\documentclass[12pt]{elsarticle}
\usepackage{graphics,graphicx,epsfig,subfigure}
\usepackage{booktabs}
\usepackage{amssymb,amsmath,amsthm,amscd,amsxtra,latexsym,amsfonts}
\usepackage{multirow}
\usepackage{xcolor}
\usepackage{enumerate}
\usepackage{float}
\usepackage{algorithmic}
\usepackage[boxed]{algorithm}
\usepackage{bm}
\usepackage{url}
\usepackage[colorlinks,linkcolor=blue,anchorcolor=black,citecolor=black]{hyperref}

\theoremstyle{definition}

\usepackage{natbib}
\biboptions{numbers,sort&compress}
\journal{Elsevier}

\begin{document}

\begin{frontmatter}

  \title{An efficient peridynamics-based statistical multiscale method for fracture in composite structure with randomly distributed particles}
  \author[nwpu,nwpuc]{Zihao Yang}
  \author[nwpu]{Shaoqi Zheng}
  \author[nwpu]{Shangkun Shen}
  \author[dlut2]{Fei Han\corref{mycorrespondingauthor}}
  \ead{hanfei@dlut.edu.cn}

  \cortext[mycorrespondingauthor]{Corresponding author}
  
  \address[nwpu]{School of Mathematics and Statistics, Northwestern Polytechnical University, Xi'an 710072, China}
  \address[nwpuc]{Innovation Center NPU Chongqing, Northwestern Polytechnical University, Chongqing 400000, China}
  \address[dlut2]{State Key Laboratory of Structural Analysis for Industrial Equipment, Department of Engineering Mechanics, Dalian University of Technology, Dalian 116023, China}
  
  \begin{abstract}
    The fracture simulation of random particle reinforced composite structures remains a challenge. Current techniques either assumed a homogeneous model, ignoring the microstructure
    characteristics of composite structures, or considered a micro-mechanical model, involving intractable computational costs.
    This paper proposes a peridynamics-based statistical multiscale (PSM) framework to simulate the macroscopic structure fracture with high efficiency.
    The heterogeneities of composites, including the shape, spatial distribution and volume fraction of particles, are characterized within the representative volume elements (RVEs), and their impact on structure failure are extracted as two types of peridynamic parameters, namely, statistical critical stretch and equivalent micromodulus.
    At the microscale level, a bond-based peridynamic (BPD) model with energy-based micromodulus correction technique is introduced to simulate the fracture in RVEs, and then the computational model of statistical critical stretch is established through micromechanical analysis.
    Moreover, based on the statistical homogenization approach, the computational model of effective elastic tensor is also established. Then, the equivalent micromodulus can be derived from the effective elastic tensor, according to the  energy density equivalence between classical continuum mechanics (CCM) and BPD models.
    At the macroscale level, a macroscale BPD model with the statistical critical stretch and the equivalent micromodulus is constructed to simulate the fracture in the macroscopic homogenized structures.
    The algorithm  framework of the PSM method is also described. Two- and three-dimensional numerical examples illustrate the validity, accuracy and efficiency of the proposed method.

  \end{abstract}
  
  
\end{frontmatter}


\section{Introduction}

Particle reinforced composite materials, such as polymer composites, concretes and short-fiber reinforced composites, have been widely used in a variety of engineering and industrial products.
The safety and reliability of composite structures are closely related to the mechanical properties of materials.
Macroscopically, the nonlinear mechanical behaviors such as fracture of composite structures depend on the microscopic heterogeneity of materials \cite{shackelford2016introduction,tang2013study}, including the shape, spatial distribution and volume fraction of particles.
Accurate modeling and efficient simulation of fracture in particle reinforced composite structures by considering the microstructure are of great significance for the optimal design of  composite structures with increased fracture toughness.
Current techniques either assumed a homogeneous model, ignoring the microstructure characteristics, or considered a micro-mechanical model, involving intractable computational costs, especially for the large-scale structures, such as concrete dam.
Therefore, it is still challenging to analyze the fracture and failure in particle reinforced composite structures from a numerical point of view.

Different techniques have been proposed to address the discontinuities of fracture. Several methods involve discontinuities of particle reinforced composites within the local continuum
framework where partial differential equations are used as the governing equations, including the finite element method \cite{hillerborg1976analysis,kassam1995finite,ayyar2006microstructure}, extended finite element method \cite{sukumar2001modeling,huynh2009extended,rostami2019xfem,spangenberger2021extended}, cohesive zone model \cite{sun2006modeling,de2013mode,sun2020prediction,quintanas2020phase}, meshfree method \cite{ghosh2013computational,bui2018analysis} and phase-field model \cite{nguyen2015phase,kuhn2016discussion,zhang2020modelling,xia2021mesoscopic}.
However, since the mesh size must be smaller than the size of particles in composite structures while employing
above numerical methods, simulating the detailed microstructure at structural scale is not computationally feasible.
A possible way to balance the accuracy and efficiency is the multiscale analysis framework, which has been proved to be effective for the linear and continuous deformation problems of composites \cite{yang2020stochastic,shu2020multiscale,yang2020novel,yang2022efficient,yang2022second,yang2017high,ma2018multi,dong2022stochastic}. And the effect of microstructure on composite structures was also taken into account.
As for the fracture analysis, attempts inspired by the multiscale framework have been made to combine a homogenized model with a small area of explicit microstructure geometry representation near the location where cracks initiate and grow \cite{canal2012intraply}.
Nevertheless, the extension to the non-linear regime, and particularly to situations involving strain localization and fracture, are much more complex and the results' correctness is not always guaranteed.
All the above methods involve discontinuities within the local continuum framework. However, the difficulty in handling discontinuities mainly arises from the basic incompatibility of
discontinuities with the partial differential equations that are used in the governing equations of the local continuum framework.

Peridynamics, which was proposed by Silling in 2000 \cite{silling2000reformulation}, is a recently developed nonlocal theory to redefine mechanical
problems by replacing the partial differential equations with integral equations, which are mathematically compatible with discontinuities.
In peridynamics, a material point is assumed to interact with surrounding points in a certain neighborhood region.
A peridynamic bond is defined between a pair of points in the neighborhood. A bond is set to break irreversibly when it is stretched beyond the critical length \cite{silling2005meshfree}. And cracks are represented by
domains that are crossed by the broken bonds. Therefore, peridynamics can effectively simulate the crack nucleation and propagation.

Some works have studied the fracture in composites with peridynamics, including particle reinforced composites.
The homogenized peridynamic models, which ignored microstructure characteristics and heterogeneities of composite materials, have been applied to predict the fracture and failure in composites.
For example, Hu et al. \cite{hu2012peridynamic} developed a homogenized peridynamics description of fiber-reinforced composites, and studied the dynamic brittle fracture and damage.
Zhou et al. \cite{zhou2017analyzing} put forth a bond-based peridynamic model to study in-plane dynamic fracture process in orthotropic composites.
Sau et al. \cite{sau2019peridynamic} introduced a numerical scheme to compute the micropolar peridynamic stress, and the stress tensor can be analyzed at the points in failure zones.
Since the influence of microstructure on nonlinear mechanics behaviors of composites cannot be ignored, peridynamic models considering the explicit representation of the microstructure of composites have been proposed.
A mesoscopic peridynamic model was proposed in \cite{li2018meso} for meso-fracture simulation of cracking process in concrete.
Dong et al. \cite{dong2021improved} proposed an improved ordinary state-based peridynamic model and employed it to study mesoscale crack initiation and propagation of concrete under uniaxial tension.
Peng et al. \cite{peng2021application} established a micro-calculation peridynamic model of concrete by the MATLAB-ABAQUS co-simulation method and calculated the mode-I fracture test and mixed-mode fracture test.
In order to reduce the cumbersome microstructural characterization and huge computation while considering the microstructures
and heterogeneities of materials as much as possible, Chen et al. \cite{mehrmashhadi2019stochastically,wu2021stochastically} proposed the intermediately
homogenized peridynamic model, which took into account some microscale information of composites, i.e., the volume fraction of particles, ignoring
the topology of component phases in microstructure, such as the shape and distribution of particles.
An effective peridynamic model considering the multiscale structure of particle reinforced composite materials to simulate fracture and failure behaviors of composite structures with high efficiency is still lacking.
This paper aims to propose a new peridynamics-based statistical multiscale framework.
It sequentially couples the peridynamics at macroscale with the peridynamics at microscale based on statistical homogenization theory to simulate the macroscopic fracture of random particle reinforced composite structures.
In the newly proposed framework (see \autoref{Fig:sec2.2-1}), the heterogeneities of composites are described as the representative volume elements (see \autoref{Fig:sec2.1-1}), and the impact of RVEs, including the shape, spatial distribution and volume fraction of particles, on structure failure are extracted as  two types of peridynamic parameters, namely, statistical critical stretch and equivalent micromodulus.
The peridynamics-based statistical multiscale framework proposed in this paper has the potential to simulate fracture in random particle reinforced composite strictures with sufficient accuracy and less computational cost, especially for large-scale structures.


The remainder of this paper is outlined as follows.
In Section 2, the multiscale representation of randomly distributed composites and the peridynamics-based statistical multiscale framework are described.
Section 3 is devoted to the proposed peridynamics-based statistical multiscale method.
Section 4 describes algorithm framework developed for the proposed approach, and a flowchart of this numerical algorithm is used to illustrate the implementation of the proposed approach.
In Section 5, two and three dimensional numerical examples are given to discuss the feasibility and effectiveness of the proposed method.
Finally, concluding remarks are presented in Section 6.

\section{Multiscale fracture of randomly distributed composites}
\subsection{Multiscale representation of randomly distributed composites}
\begin{figure}[ht]
  \centering
  \includegraphics[scale=0.55]{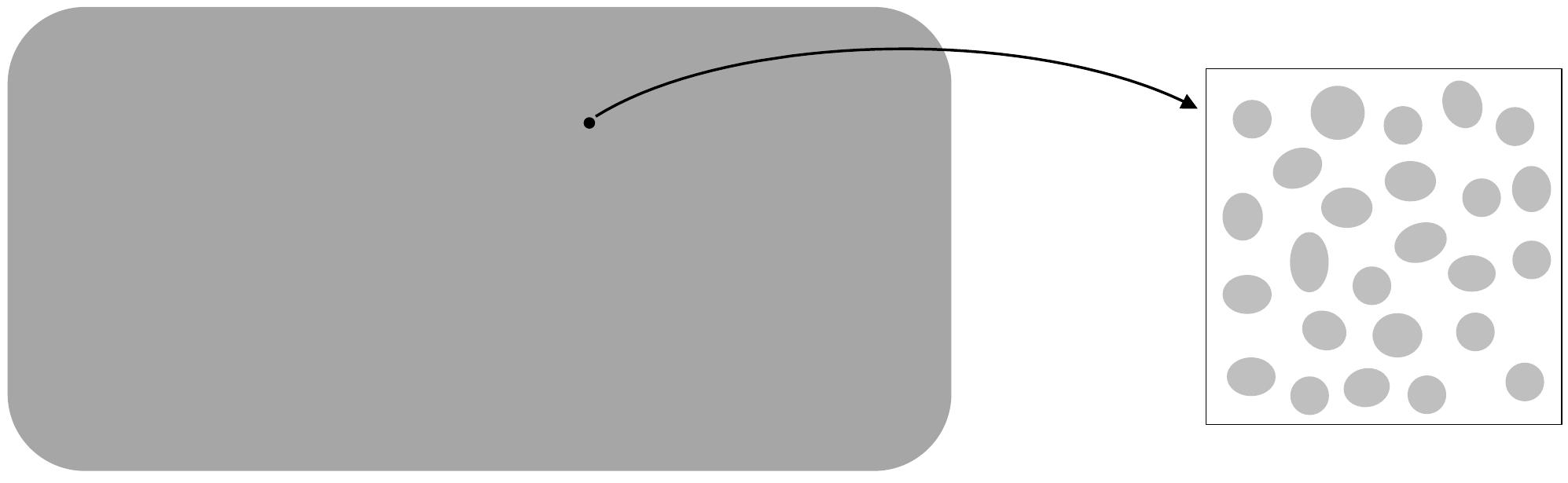}
  \caption{Composite structures $\Omega$ (left) and the microscopic RVE (right).}\label{Fig:sec2.1-1}
\end{figure}

We consider the composite structures made from the matrix and reinforced particles in this study. For the three-dimensional (3D) case, all the particles are considered as ellipsoids or polyhedrons inscribed inside the ellipsoids.
The size of each particle is denoted by the long axis $a$ of corresponding ellipsoid.
For this kind of composite structure, a multiscale representation \cite{han2010statistical}, as shown in \autoref{Fig:sec2.1-1}, is introduced as follows:


1) There exists a constant $\varepsilon$ satisfying $a\ll\varepsilon\ll L$, where $L$ denotes the size of the investigated composite domain $\Omega$ at the macroscale. Thus, the structure made from composite materials can be regarded as a set of RVEs at the microscale with the same size $\varepsilon$. The composite structure is supposed to be stationary random \cite{jikov2012homogenization}, i.e., all the RVEs have the same probability distribution of particles. Then, the composite structure has periodically random distribution of particles, and can be represented by the probability distribution $P$ inside a typical RVE with $\varepsilon$-size.



2) Each ellipsoid can be defined by nine random parameters: the size parameters, i.e., the length $a$, $b$ and $c$ of the ellipsoid in three axes, the location parameters, i.e., the coordinates $(x_{1}^0, x_{2}^0, x_{3}^0)$ of central point, the orientation parameters, i.e., the Euler angles $(\psi_1, \psi_2, \psi_3)$ of the rotations. Suppose that there are $N$ ellipsoids inside a RVE $\varepsilon Y^m$, then we can define a sample $\omega^m$ of particle distribution in a normalized RVE $Y^m$ as $\omega^m=(\bm{I}_1, \bm{I}_2\cdots,\bm{I}_N)\in P$, where $m =1,2,3, \ldots$ denotes the index of samples, and $\bm{I}_j$ denotes a sample of nine random parameters for the $j$-th ellipsoid, generated by their probability density functions.


For the two-dimensional (2D) case, all the particles are considered as ellipses defined by five random parameters, including the length $a$ and $b$ of two axes, the coordinates $(x_{1}^0, x_{2}^0)$ of the central point, and the orientation parameters $\psi$. A sample of particles distribution for a 2D RVE can be obtained similarly.


By the above multiscale representation, the investigated structure $\Omega$ is logically
composed of $\varepsilon$-size cells subjected to same probability
distribution model $P$, i.e., $\Omega=\bigcup_{(\omega^m,t\in Z)}\varepsilon(Y^m+t)$.
For the whole structure $\Omega$, define
$\omega=\{\omega^m|\bm x\in \varepsilon Y^m \subset \Omega\}$. Thus the elastic
tensor of the materials can be periodically expressed as $\mathbb{A}(\frac{\bm x}{\varepsilon},\omega\bigr)=\{a_{ijkl}\bigl(\frac{\bm x}{\varepsilon},\omega\bigr)\}$, and for a given sample $\omega^m$, the elastic tensor can be defined as follows:
\begin{equation}\label{eq:coeff}
  a_{ijkl}\bigl(\frac{\bm x}{\varepsilon},\omega^m\bigr)=
  \left\{
  \begin{aligned}
     & a^1_{ijkl}, \quad \text{if}\,\,\, \bm x\in \bigcup_{i=1}^{N}e_i,                 \\
     & a^2_{ijkl}, \quad \text{if}\,\,\, \bm x\in \varepsilon Y^m-\bigcup_{i=1}^{N}e_i,
  \end{aligned}
  \right.
\end{equation}
where $\varepsilon Y^m$ denotes the domain of a RVE belonging to $\Omega$, $e_i$
is the $i$-th ellipse in the $\varepsilon Y^m$, and $a^1_{ijkl}$ and $a^2_{ijkl}$ are the
elastic coefficients of particles and matrix, respectively.
%

We introduce a variable $\bm y=\bm x/\varepsilon \in Y^m$ which denotes local coordinates on the normalized RVE $Y^m$, and then the elastic tensor can be expressed as
\begin{equation}\label{eq.coeff1}
  \mathbb{A}(\frac{\bm x}{\varepsilon},\omega\bigr) =\mathbb{A}(\bm y,\omega\bigr).
\end{equation}
\subsection{Peridynamics-based statistical multiscale (PSM) methodological framework}

For the composite structure with randomly distributed particles, it macroscopically exhibits nonlinear mechanical behavior and size effects, which originates from its microscopic heterogeneity.
In order to explain the cracking behavior in the composite structure more reasonably, the microstructural characteristics and heterogeneities cannot be ignored. However, this will inevitably lead to a huge amount of computation, especially for the large structures made from composite materials with randomly distributed particles.
Based on the multiscale representation of randomly distributed composites and inspired by the statistical homogenization theory, an efficient PSM framework is proposed to analyze the fracture of random composite structure.

\begin{figure}[ht]
  \centering
  \includegraphics[scale=0.65]{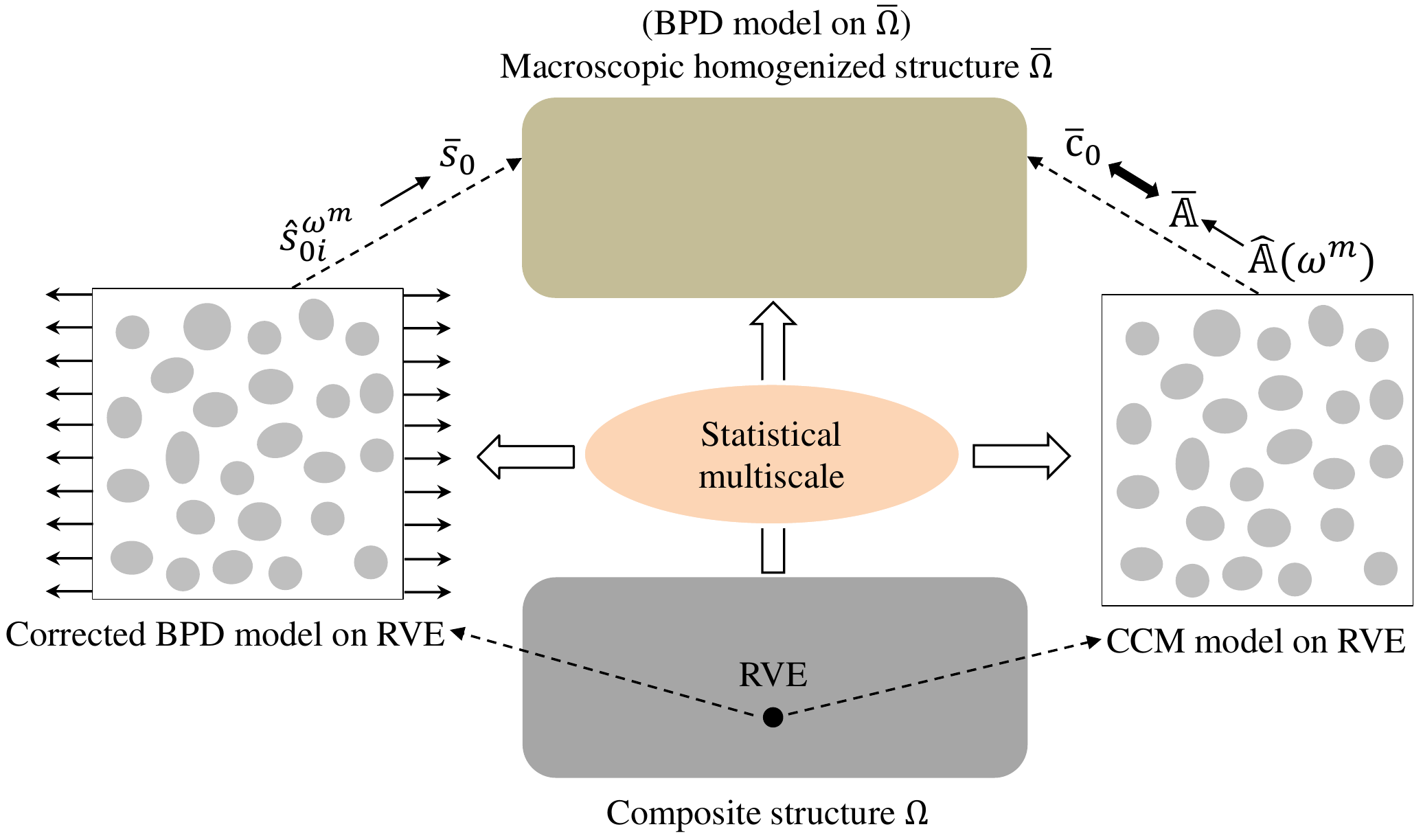}
  \caption{Methodological framework.}\label{Fig:sec2.2-1}
\end{figure}

In the newly proposed framework, as shown in \autoref{Fig:sec2.2-1}, the heterogeneities of composites, including the shape, spatial distribution and volume fraction of particles, are characterized within the RVEs, and their impacts on the macroscopic structure failure are extracted as two types of peridynamic parameters, namely, statistical critical stretch $\bar{s}_0$ and equivalent micromodulus $\bar{c}_0$.
Since the fracture of composite structure occurs at some local RVEs, which is related to the loads on the RVEs, the micromechanical computation is applied to analyze the fracture in the RVEs to obtain the statistical critical stretch $\bar{s}_0$.
In order to estimate the critical stretch of the macroscopic homogenized material, the fracture of the RVEs of composites is simulated by the BPD model, which has advantages in simulating fractures. For the fracture of the RVEs, the mechanical response near the interfaces between particles and matrix is critical since cracks often initiate near the composite interfaces. However, the nonlocal PD bond model fails to accurately characterize the properties of materials on both sides of composite interfaces, so the elastic solution consistent with the local CCM model cannot be obtained, resulting in inaccurate fracture simulation in the RVEs. An alternative approach is to modify the mechanical parameters of PD bonds near the composite interfaces based on the elastic solution of CCM model and the pointwise elastic energy density equivalence between CCM and BPD models, so as to obtain the consistent solutions with local CCM model. Consequently, the corrected BPD model can better simulate the fracture of the RVEs \cite{YANG2022CICP}.
On the other hand, the above correction method destroys the stationary randomness of PD micromodulus distribution in the RVEs of composites, since the PD micromodulus correction is carried out one by one for the RVEs with specific random samples. Therefore, the homogenization strategy for the BPD model, i.e., upscaling the microscale PD micromodulus on the RVEs into the macroscale effective elastic tensor \cite{madenci2020peridynamic}, will no longer be applicable. Fortunately, estimating the effective elastic tensor of composites does not involve discontinuous deformation, and the elastic solution obtained from the CCM model still works. Therefore, the traditional homogenization method based on the CCM model can be used to obtain the effective elastic tensor of macroscale homogenized material. Then the equivalent PD micromodulus can be derived based on the elastic energy density equivalence between CCM and BPD models. Specifically, the PSM methodology framework proposed in this paper includes following three steps:

\begin{itemize}
  \item At the microscale level, an energy-based correction method will be introduced to define the micromodulus of interface bonds in RVEs, which makes sure that for the elastic deformation of RVEs, the elastic energy density of CCM and BPD models will be equivalent.
        And an improved microscale BPD model is obtained to analyze the fracture of RVEs. We can further define the effective
        critical stretch $\hat{s}_{0i}^{\omega^m}(i=1,2,3)$ along different stretching directions as the critical tensile deformation when the RVE breaks completely.
        According to the probability theory, we establish the computational model of statistical critical stretch $\bar{s}_0$ related to $\hat{s}_{0i}^{\omega^m}$ with different samples $\omega^m(m=1,2,\ldots,M)$.
  \item At the microscale level, the statistical homogenization approach is introduced to solve the CCM model of the composite materials, and we build the computational model of homogenized elastic coefficients $\hat{\mathbb{A}}(\omega^m)$ defined as integral average over the RVEs.
        According to the probability theory, we give the computational models of effective elastic tensor $\bar{\mathbb{A}}$ related to $\hat{\mathbb{A}}(\omega^m)$ with different samples $\omega^m(m=1,2,\ldots,M)$.
        Further, according to the  energy density equivalence between CCM and BPD models of macroscopic homogenized materials, it is able to derive the equivalent micromodulus $\bar{c}_0$ from the effective elastic tensor $\bar{\mathbb{A}}$.
  \item At the macroscale level, based on above homogenization models, the macroscale BPD model with statistical critical stretch $\bar{s}_0$ and equivalent micromodulus $\bar{c}_0$ is then constructed to analyze the fracture in the  macroscopic homogenized structures $\bar{\Omega}$.
\end{itemize}

\section{Peridynamics-based statistical multiscale (PSM) method}
\subsection{Microscopic peridynamic modeling for RVEs of composites}
\subsubsection{Micromodulus tensors for the RVEs}

\begin{figure}[ht]
  \centering
  \includegraphics[scale=0.6]{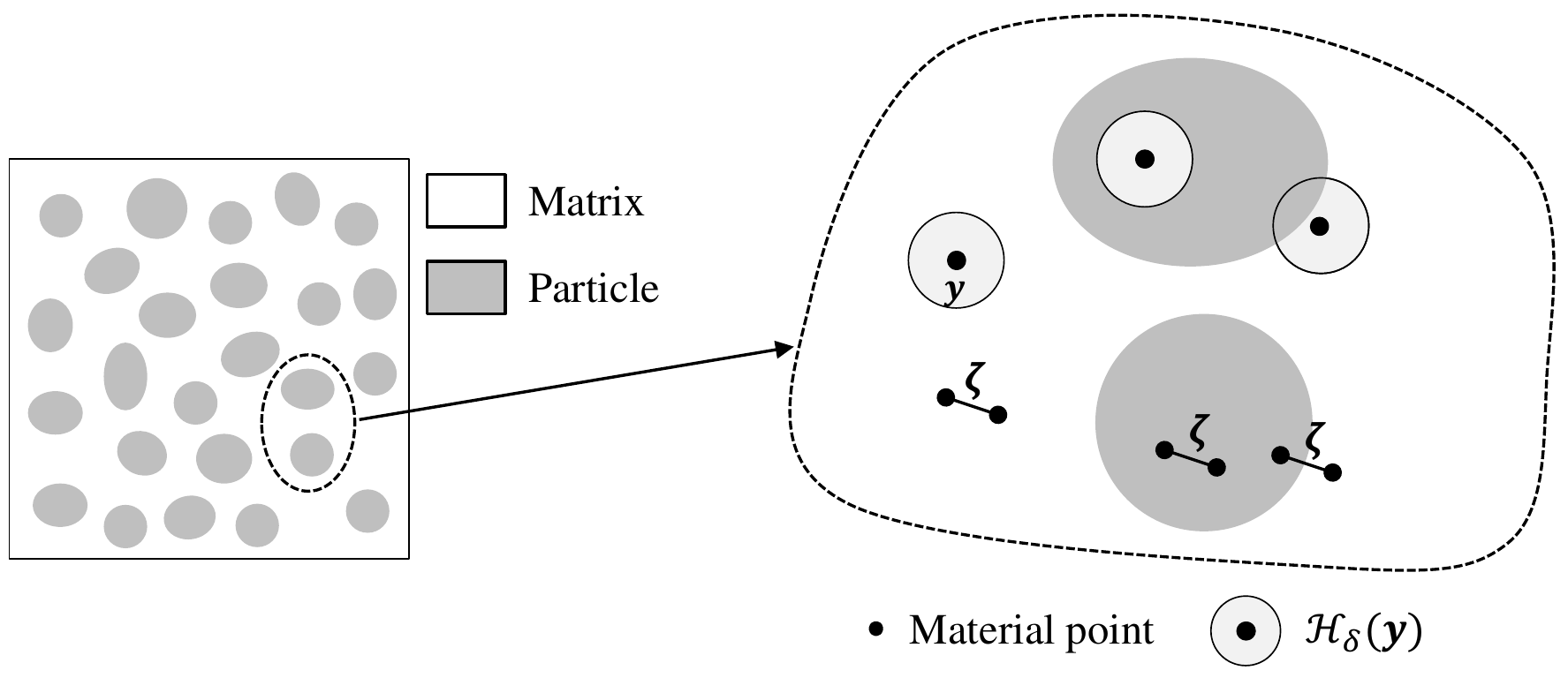}
  \caption{Three kinds of PD bonds in the RVE $Y^m$.}\label{Fig:sec3.2.1-1}
\end{figure}

For the two-phase composite materials, there exist three types of PD bonds $\bm \zeta$ in the RVE, including the particle bonds, matrix bonds and interface bonds, as shown in \autoref{Fig:sec3.2.1-1}. For a specified sample $\omega^m$, the PD bonds in the RVE can be defined as $\bm \zeta=\bm y'-\bm y(\bm y,\bm y'\in Y^m)$,
and corresponding micromodulus can be expressed by $c_0(\bm y,\bm \zeta,\omega^m)$.
Since the composite interfaces play an important role on the material failure, it is very important to accurately define the micromodulus of interface bonds in the BPD model.
In our previous work \cite{YANG2022CICP}, an energy-based micromodulus correction method was proposed  to determine the mechanical properties of interface bonds in the BPD model for the composite structure with randomly distributed particles.


A stretch boundary condition along the positive direction of $y_i$-axis ($i=1,2,3$) is specified on one side of the RVE $Y^m$ with the other sides being simply supported and a middle point being fixed (see \autoref{Fig:sec3.2.2-1}(a) for the example of $y_1$-axis). According to the CCM theory, the displacement field $\bm u^i_{\text{ccm}}(\bm y,\omega^m)(i=1,2,3)$ can be obtained through the equations of elasticity defined on the RVE $Y^m$.
For the infinitesimal elastic deformation, the elastic energy density of CCM and BPD models with the same distribution of displacement field should be equivalent.
Thus, the elastic energy density of CCM and BPD models with the displacement field $\bm u^i_{\text{ccm}}(\bm y,\omega^m)(i=1,2,3)$, i.e., $W^{i}_{\text{ccm}}(\bm y,\omega^m)$ and $W^i_{\text{pd}}(\bm y,\omega^m)(i=1,2,3)$, can be further obtained as
\begin{equation}\label{eq3.1.1}
  W^{i}_{\text{ccm}}(\bm y,\omega^m)=\frac{1}{2}\bm \epsilon^i(\bm y,\omega^m):\mathbb{A}(\bm y,\omega^m):\bm \epsilon^i(\bm y,\omega^m),
\end{equation}
\begin{equation}\label{eq3.1.2}
  W^i_{\text{pd}}(\bm y,\omega^m)=\frac{1}{4}\int_{\mathcal{H}_{\delta}(\bm y)} \frac{c_0(\bm y,\bm \zeta,\omega^m)+c_0(\bm y',\bm \zeta,\omega^m)}{2}\bigl(\bm \zeta\cdot(\bm u^i_{\text{ccm}}(\bm y,\omega^m)-\bm u^i_{\text{ccm}}(\bm y',\omega^m)) \bigr)^2dV_{\bm y'},
\end{equation}
where $\bm \epsilon^i(\bm y,\omega^m)(i=1,2,3)$ denotes the strain tensor defined as
\begin{equation*}
  \bm \epsilon^i(\bm y,\omega^m)=\frac{1}{2}\bigl(\nabla\bm u^i_{\text{ccm}}(\bm y,\omega^m)+\nabla ^\text{T} \bm u^i_{\text{ccm}}(\bm y,\omega^m)\bigr).
\end{equation*}
Based on the energy equivalence between CCM and BPD model for elastic continuous deformation of the RVE, the correction parameters $\bm \alpha(\bm y)=(\alpha_1(\bm y),\alpha_2(\bm y),\alpha_3(\bm y))$ and $\bm \beta(\bm y)=(\beta_1(\bm y'),\beta_2(\bm y'),\beta_3(\bm y'))$ are defined to modify the micromodulus at the points $\bm y$ and $\bm y'$ in the RVE $Y^m$, respectively
\begin{equation}\label{eq3.1.3}
  \alpha_i(\bm{y}) = \frac{W^i_{\text{ccm}}(\bm{y},\omega^m)}{W^i_{\text{pd}}(\bm{y},\omega^m)},\quad \beta_i(\bm{y}') = \frac{W^i_{\text{ccm}}(\bm{y}',\omega^m)}{W^i_{\text{pd}}(\bm{y}',\omega^m)}.
\end{equation}
Thus, the micromodulus at the point $\bm y$ and point $\bm y'$ are modified as
\begin{equation}\label{eq3.1.4}
  \hat{c}_i^{\omega^m}(\bm{y},\bm \zeta) = \alpha_i(\bm y)c(\bm y,\bm \zeta,\omega^m),\quad \hat{c}_i^{\omega^m}(\bm{y}',\bm \zeta) =\beta_i(\bm y')c(\bm y',\bm \zeta,\omega^m).
\end{equation}
Further, the micromodulus of a bond $\bm \zeta$ can be defined by the harmonic averaging of $\hat{c}_i^{\omega^m}(\bm y,\bm \zeta)$ for $\bm y$ and $\hat{c}_i^{\omega^m}(\bm{y}',\bm \zeta)$ for $\bm y'$, namely
\begin{equation}\label{eq3.1.5}
  k_i^{\omega^m}(\bm{y},\bm \zeta) = \frac{2\hat{c}_i^{\omega^m}(\bm{y},\bm \zeta)\hat{c}_i^{\omega^m}(\bm{y}',\bm \zeta)}{\hat{c}_i^{\omega^m}(\bm{y},\bm \zeta)+\hat{c}_i^{\omega^m}(\bm{y}',\bm \zeta)}.
\end{equation}
Finally, the scalar scaling is applied to define the micromodulus for the bond $\bm \zeta$
\begin{equation}\label{eq3.1.8-c}
  \tilde{c}_0^{\omega^m}(\bm{y}, \bm \zeta) = \frac{1}{\sqrt{\bigl(\frac{n_1}{k_1^{\omega^m}(\bm{y},\bm \zeta)}\bigr)^{2} + \bigl(\frac{n_2}{k_2^{\omega^m}(\bm{y},\bm \zeta)}\bigr)^{2} + \bigl(\frac{n_3}{k_3^{\omega^m}(\bm{y},\bm \zeta)}\bigr)^{2}}},
\end{equation}
where $\bm \zeta/ |\bm \zeta| = (n_1, n_2, n_3)^{\text{T}}$ denotes the unit direction vector of bond $\bm \zeta$ in the RVE.

\subsubsection{Statistical critical stretch: microscopic peridynamic simulation}
Based the definition of micromodulus by energy-based correction method,
the improved BPD equilibrium equation, defined on the RVE $Y^m$, is written as follows:
\begin{equation}\label{eq:sec3.2.2-1}
  \int_{\mathcal{H}_{\delta}(\bm y)} \bm f^{\omega^m}(\bm y', \bm y) dV_{\bm y'} + \bm b(\bm y) = 0, \quad \bm y,\bm y'\in Y^m.
\end{equation}
We denote by $\mathcal{H}_{\delta}(\bm y)$ the neighborhood in the RVE, and the pairwise force function $\bm f^{\omega^m}(\bm y', \bm y)$ can be defined as
\begin{equation}\label{eq:sec3.2.2-2}
  \bm f^{\omega^m}(\bm y', \bm y) =\frac12 \tilde{c}_0^{\omega^m}(\bm y, \bm \zeta) \frac{\bm \zeta \otimes \bm \zeta}{|\bm \zeta|^2}\cdot\bigl( \bm u(\bm y') - \bm u(\bm y) \bigr)\hat{\kappa}(t,\bm \zeta),
\end{equation}
where $\bm u$ denotes the displacement in the RVE, and $\hat{\kappa}(t,\bm \zeta)$ denotes the history-dependent scalar-valued functions defined by Eq. \eqref{eq:sec3.2.2-5}.

\begin{figure}[ht]
  \centering
  \includegraphics[scale=0.8]{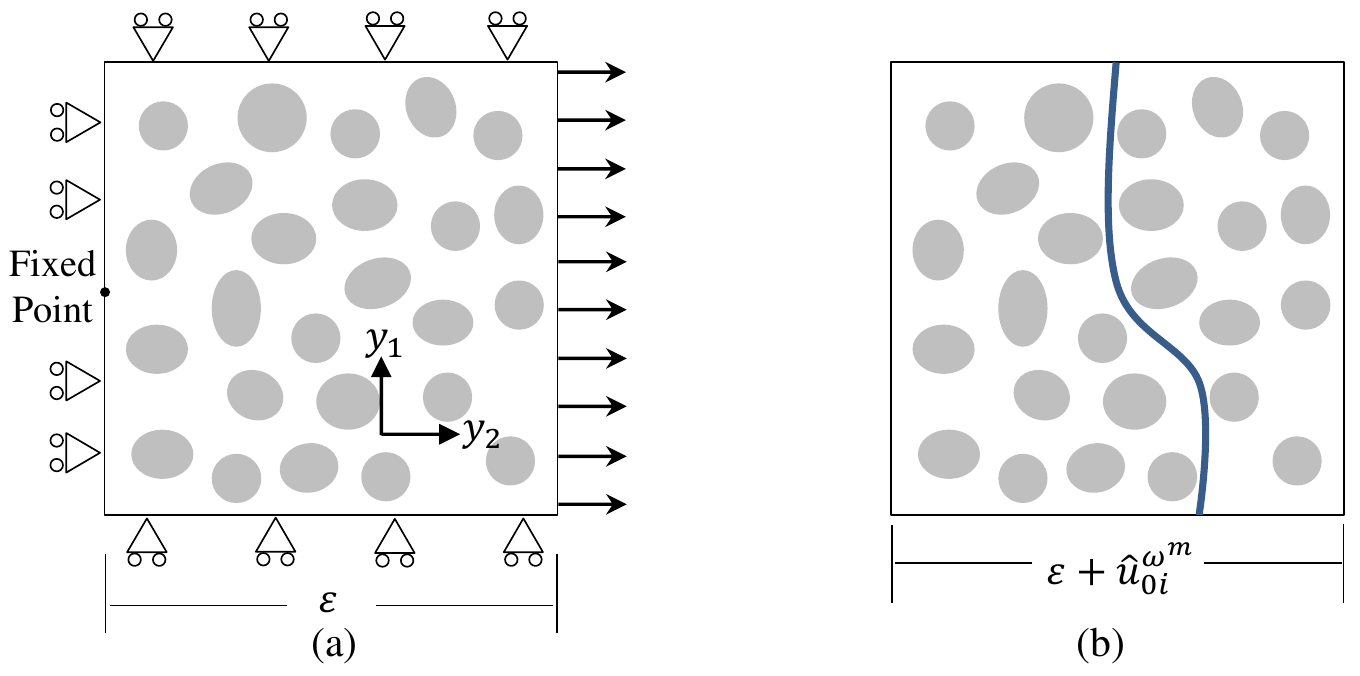}
  \caption{The  critical bond stretch of the RVE $Y^m$ with a stretch boundary condition along the positive direction of $y_1$-axis.}\label{Fig:sec3.2.2-1}
\end{figure}

The bond stretch $s_{\bm \zeta}$ in the RVE $Y^m$ can be defined as follows
\begin{align}\label{eq:sec3.2.2-4}
  s_{\bm \zeta}=\frac{|\bm \zeta+\boldsymbol{u}(\boldsymbol{y}+\bm \zeta)-\boldsymbol{u}(\boldsymbol{y})|-|\bm \zeta|}{|\bm \zeta|}.
\end{align}
In this study, since the two-phase particle reinforced composite structures are considered, the failure law of the RVE $Y^m$ in composite structure is implemented by defining different history-dependent scalar-valued functions $\kappa$ for the matrix bonds, particle bonds and interface bonds as follows
\begin{align}\label{eq:sec3.2.2-5}
  \hat{\kappa}(t,\bm \zeta)=
  \left\{
  \begin{array}{lr}
    1, \, \text{if}\,\,\, s_{\bm \zeta}<s_0(\bm y,\bm \zeta) \quad \forall \ 0\leq t'\leq t, \\
    0, \, \text{otherwise},
  \end{array}
  \right.
\end{align}
where $t'$ and $t$ denote the computational steps, $\bm \zeta$ denotes the three kinds of PD bonds, and $s_0(\bm y,\bm \zeta)$ denotes critical  stretch value corresponding to the particle bonds, matrix bonds and interface bonds.
A stretch boundary condition along the positive direction of $y_i$-axis ($i=1,2,3$) is specified on one right side of the RVE $Y^m$ with the other sides being simply supported and a middle point being fixed (see \autoref{Fig:sec3.2.2-1}(a) for the example of $y_1$-axis).
According to the improved BPD model and failure criterion of PD bonds, a critical tensile displacement $\hat{u}_{0i}^{\omega^m}$ along the direction of $y_i$-axis $(i=1,2,3)$ is reached until the RVE breaks completely (see \autoref{Fig:sec3.2.2-1}(b)).
Since the RVE corresponds to a material point $\bm x$ in the macro-structure,
the effective critical bond stretch $\hat{s}_{0i}^{\omega^m}(i=1,2,3)$ along the $x_i$-direction $(i=1,2,3)$ can be defined as follows
\begin{equation}\label{eq:sec3.2.2-6}
  \hat{s}_{0i}^{\omega^m}=\frac{\varepsilon+\tilde{u}_{0i}^{\omega^m}-\varepsilon}{\varepsilon}=\frac{\tilde{u}_{0i}^{{\omega^m}}}{\varepsilon}.
\end{equation}

From the Kolmogorov strong law of large number, we can evaluate the expected critical stretch $\bar{s}_{0i}(i=1,2,3)$ along $x_i$-direction by taking different samples $\omega^{m}(m=1,2,\cdots,M)$
\begin{equation}\label{eq:s0i}
  \bar{s}_{0i}=\lim_{M\rightarrow\infty}\frac{\displaystyle\sum_{m=1}^{M}\hat{s}_{0i}^{\omega^{m}}}{M}.
\end{equation}
The scalar scaling is then applied to define the statistical critical stretch $\bar{s}_0$ corresponding to the bond $\bm \xi$ as follows
\begin{equation}\label{eq:s0}
  \bar{s}_0(\bm \xi)= \frac{1}{\sqrt{\bigl(\frac{v_1}{\bar{s}_{01}}\bigr)^{2} + \bigl(\frac{v_2}{\bar{s}_{02}}\bigr)^{2} + \bigl(\frac{v_3}{\bar{s}_{03}}\bigr)^{2}}},
\end{equation}
where $\bm \xi/ |\bm \xi| = (v_1, v_2, v_3)^{\text{T}}$ denotes the unit direction vector of bond $\bm \xi=\bm x-\bm x'(\bm x,\bm x'\in \bar{\Omega})$ in the macroscale homogenized structure $\bar{\Omega}$, which will be discussed in Section 3.2.

\subsubsection{Equivalent  micromodulus: statistical homogenization approach}
The CCM model for description of elastic mechanical behavior in randomly distributed composite structure $\Omega$ is given by
\begin{eqnarray}\label{eq:sec3.2.3-1}
  \vspace{2mm}\displaystyle  -\frac{\partial}{\partial
    x_j}\left[a_{ijkl}({\frac{\bm x}{\varepsilon},\omega})\frac{1}{2}\left(\frac{\partial u^{\varepsilon}_k}{\partial
      x_l}+\frac{\partial u^{\varepsilon}_l}{\partial x_k}\right)\right]=q_i(\bm x),\ \ \text{in}\ \Omega,
\end{eqnarray}
where $\bm {u}^{\varepsilon}({\bm x},\omega)$ is the displacement vector.
As mentioned in Section 2.1, $\bm x$ and $\bm y$ denote the coordinate system of the macroscale defined on the structure $\Omega$ and the microscale defined on the normalized RVE, respectively, they are related to each other as $\bm y=\bm x/{\varepsilon}$.
Differentiation with respect to $\bm x$ can be defined as
\begin{eqnarray}\label{eq:sec3.2.3-3}
  \frac{\partial }{\partial x_i}=\frac{\partial }{\partial
    x_i} + \varepsilon^{-1} \frac{\partial }{\partial y_i},\
  (i=1,2,3).
\end{eqnarray}
Inspired by the classical mathematical homogenization theory, one supposes that $\bm {u}^{\varepsilon}$ can be expressed as following asymptotic forms for a specified sample $\omega^m$ \cite{yang2020stochastic}
\begin{eqnarray}\label{eq:sec3.2.3-4}
  \bm u^{\varepsilon}(\bm x,\omega)=\bm u^{0}(\bm x,\bm y,\omega^m)+\varepsilon \bm{u}^{1}(\bm x,\bm y,\omega^m)+\varepsilon^2 \bm u^{2}(\bm x,\bm y,\omega^m)+ \cdots.
\end{eqnarray}

Introducing Eqs. \eqref{eq:sec3.2.3-3} and \eqref{eq:sec3.2.3-4} into Eq. \eqref{eq:sec3.2.3-1} and matching the terms of the same order of $\varepsilon$,
%
%
and a series of equations are obtained if the coefficients of $\varepsilon^i(i=-2,-1,0,\ldots)$ from both sides of above equation are required to be equal to each other.
From the coefficients of $\varepsilon^{-2}$ following equation can be obtained $\bm u^{0}=\bm u^{0}(\bm x)$.
From the coefficients of $\varepsilon^{-1}$ following equation can be obtained
\begin{eqnarray*}\label{eq:sec3.3.1-8}
  \bm {u}(\bm x,\bm y,\omega^m)= N_{\alpha}(\bm y,\omega^m)\frac{\partial \bm {u}^0}{\partial x_{\alpha}}(\bm x)+\bm {\breve{u}}(\bm x),
\end{eqnarray*}
where $\bm {\breve{u}}^1$ is a function which is independent of $\bm y$, and ${\bm N}_{\alpha}(\bm y,\omega^m)$ is a matrix-valued function defined in the RVE $Y^m$ as follows
\begin{equation}\label{eq:sec3.2.3-5}
  \frac{\partial }{{\partial {y_j}}}\left[ {{a_{ijkl}}(\bm y,\omega^m)\frac{1}{2}\left( {\frac{{\partial {N_{\alpha km}}}}{{\partial {y_l}}} + \frac{{\partial {N_{\alpha lm}}}}{{\partial {y_k}}}} \right)} \right] =  - \frac{{\partial {a_{ijm\alpha }}(\bm y,\omega^m)}}{{\partial {y_j}}},\quad \bm y \in Y^m,
\end{equation}
with following boundary condition
\begin{equation}\label{eq:sec3.2.3-6}
  \bm {N}_{\alpha m}(\bm y,\omega^m)=0, \quad \forall \bm y\in\partial Y^m.
\end{equation}
%
%
From the coefficients of $\varepsilon^{0}$ following equation defined in the macroscale homogenized structure $\hat{\Omega}$ corresponding to the sample $\omega^m$ can be obtained
\begin{eqnarray}\label{eq:sec3.2.3-7}
  - \frac{\partial }{{\partial {x_j}}}\left( {\hat{a}_{ijkl}(\omega^m)\frac{1}{2}\left( {\frac{{\partial u_k^{\text{0}}({\bm x})}}{{\partial {x_l}}} + \frac{{\partial u_l^{\text{0}}({\bm x})}}{{\partial {x_k}}}} \right)} \right) = {q_i}({\bm x}),\quad \bm x \in \hat{\Omega},
\end{eqnarray}
where $\hat{a}_{ijkl}(\omega^m)$ is the component of the homogenized elastic tensor $\hat{\mathbb{A}}(\omega^m)$ defined as
\begin{equation}\label{eq:sec3.2.3-8}
  \hat{a}_{ijkl}(\omega^m) = \frac{1}{{\left| Y^m \right|}}\int_{Y^m} {\left( {a_{ijkl}^{}({\bm y},\omega^m) + a_{ijpq}^{}({\bm y},\omega^m)\frac{1}{2}\left( {\frac{{\partial N_{kpl}^{}({\bm y},\omega^m)}}{{\partial {y_q}}} + \frac{{\partial N_{kql}^{}({\bm y},\omega^m)}}{{\partial {y_p}}}} \right)} \right)} dy.
\end{equation}
From Kolmogorov strong law of large number, we can evaluate the effective elastic tensor $\bar{a}_{ijkl}$ by taking different samples $\omega^{m}(m=1,2,\cdots,M)$
\begin{equation}\label{eq:aijkl}
  \bar{a}_{ijkl} =\lim_{M\rightarrow\infty}\frac{\displaystyle\sum_{m=1}^{M}\hat{a}_{ijkl}(\omega^m)}{M}.
\end{equation}

The equivalent micromodulus $\bar c_0(\bm x,\bm x')$ for the BPD model on homogenized structure $\bar{\Omega}$ can be defined through the effective elastic tensor $\bar{\mathbb{A}}=\{\bar{a}_{ijkl}\}$ as follows \cite{azdoud2013morphing}
\begin{equation}\label{eq:eff-microc}
  \bar{\mathbb{A}}=\frac{1}{2}\int_{\mathcal{H}_{\delta}(\bm x)} \bar c_0(\bm x,\bm \xi) \frac{\bm \xi \otimes \bm \xi \otimes \bm \xi \otimes \bm \xi}{|\bm \xi|^2} dV_{\bm x'}.
\end{equation}
\subsection{Statistical peridynamic modeling for macroscopic homogenized structure}

According to the statistical homogenization theory, the effect of microstructure, including the shape, spatial distribution and volume fraction of particles, on the  structure failure is homogenized as two types of peridynamic parameters, namely, the statistical critical stretch $\bar{s}_0$ and equivalent micromodulus $\bar{c}_0$. Therefore, the BPD equilibrium equation, defined on homogenized structure $\bar{\Omega}$, is written as follows:

\begin{equation}\label{eq:sec3.3-1}
  \int_{\mathcal{H}_{\delta}(\bm x)} \bm f(\bm x', \bm x) dV_{\bm x'} + \bm b(\bm x) = 0, \quad \bm x,\bm x'\in \bar{\Omega}.
\end{equation}
We denote by $\mathcal{H}_{\delta}(\bm x)$ the neighborhood in macroscale homogenized structure $\bar{\Omega}$, and the pairwise force function $\bm f(\bm x', \bm x)$ can be defined as
\begin{equation}\label{eq:sec3.3-2}
  \bm f (\bm x', \bm x) = \bar{c}_0(\bm x,\bm \xi) \frac{\bm \xi \otimes \bm \xi}{|\bm \xi|^2}\cdot \bigl( \bm u(\bm x') - \bm u(\bm x) \bigr)\bar{\kappa}(t,\boldsymbol{\xi}),
\end{equation}
where $\bm u$ denotes the displacement in the homogenized structure $\bar{\Omega}$, and $\bar{\kappa}(t,\bm \xi)$ denotes the history-dependent scalar-valued functions defined by Eq. \eqref{eq:sec3.3-5}.

The bond stretch $s_{\bm \xi}$ in the homogenized structure $\bar{\Omega}$ can be defined as follows
\begin{align}\label{eq:sec3.3-4}
  s_{\bm \xi}=\frac{|\bm \xi+\boldsymbol{u}(\boldsymbol{x}+\boldsymbol{\xi})-\boldsymbol{u}(\boldsymbol{x})|-|\boldsymbol{\xi}|}{|\boldsymbol{\xi}|}.
\end{align}
The failure law of the homogenized structure $\bar{\Omega}$ is implemented by defining history-dependent scalar-valued functions as follows
\begin{align}\label{eq:sec3.3-5}
  \bar{\kappa}(t,\boldsymbol{\xi})=
  \left\{
  \begin{array}{lr}
    1, \, \text{if}\,\,\, s_{\bm \xi}<\bar{s}_0(\bm \xi) \quad \forall \ 0\leq t'\leq t, \\
    0, \, \text{otherwise},
  \end{array}
  \right.
\end{align}
where $t'$ and $t$ denote the computational steps.

\section{Peridynamics-based statistical multiscale algorithm}
The proposed PSM framework includes some key models. The energy-based micromodulus correction model consists of Eqs. \eqref{eq3.1.1}-\eqref{eq3.1.8-c}, the statistical critical stretch model consists of Eqs. \eqref{eq:sec3.2.2-6}-\eqref{eq:s0}, the equivalent micromodulus model consists of Eqs. \eqref{eq:sec3.2.3-8}-\eqref{eq:eff-microc}, and the statistical peridynamic model for macroscopic homogenized structure consists of Eqs. \eqref{eq:sec3.3-1}-\eqref{eq:sec3.3-5}, which comprise a closed solving system that can be used for performing
objective simulation of composite structure with randomly distributed particles.

\begin{algorithm}[H]
  \renewcommand{\thealgorithm}{of PSM}
  \caption{computation steps of the proposed PSM method}
  
  {\bf Algorithm I: Definition of corrected micromodulus $\tilde{c}_0^{\omega^m}(\bm{y}, \bm \zeta)$}\\
  {\bf Input:} geometry of RVE $Y^m$ with sample $\omega^m$, elastic properties of particles and matrix\\
  {\bf Output:} Micromodulus $\tilde{c}_0^{\omega^m}(\bm{y}, \bm \zeta)$ of PD bonds in $Y^m$
  \begin{algorithmic}[1]
    \STATE{Compute displacement field of the RVE through CCM theory.}
    \STATE{Compute elastic energy density of the CCM and BPD models by Eqs. \eqref{eq3.1.1} and \eqref{eq3.1.2}.}
    \STATE{Compute correction factor along different directions similarly by Eq. \eqref{eq3.1.3}.}
    \STATE{Compute $\tilde{c}_0^{\omega^m}(\bm{y}, \bm \zeta)$ by harmonic average and scalar scaling in Eqs. \eqref{eq3.1.5}-\eqref{eq3.1.8-c}.}
  \end{algorithmic}
  
  \vspace{6 pt}
  
  {\bf Algorithm II: Computation of   statistical critical stretch $\bar{s}_{0}$}\\
  {\bf Input:} geometry of RVE $Y^m$, micromodulus $\tilde{c}_0^{\omega^m}(\bm{y}, \bm \zeta)$\\
  {\bf Output:} statistical critical stretch $\bar{s}_{0}$
  \begin{algorithmic}[1]
    \STATE{Compute displacement of the RVE by Eq. \eqref{eq:sec3.2.2-1} with a $y_i$-direction boundary condition.}
    \STATE{Compute effective critical stretch $\hat{s}_{0i}^{\omega^m}$ related to sample $\omega^m$ by Eq. \eqref{eq:sec3.2.2-6}.}
    \STATE{Compute statistical critical stretch $\bar{s}_{0}(\bm \xi)$ by Eqs. \eqref{eq:s0i} and \eqref{eq:s0}.}
  \end{algorithmic}
  
  \vspace{6 pt}
  
  {\bf Algorithm III: Computation of  equivalent micromodulus $\bar{c}_0$}\\
  {\bf Input:} geometry of RVE $Y^m$, elastic properties of particles and matrix\\
  {\bf Output:} equivalent micromodulus $\bar{c}_0$
  \begin{algorithmic}[1]
    \STATE{Compute cell function ${\bm N}_{\alpha}(\bm y,\omega^m)$ by solving Eqs. \eqref{eq:sec3.2.3-5} and \eqref{eq:sec3.2.3-6}.}
    \STATE{Compute homogenized elastic tensor $\{\hat{a}_{ijkl}(\omega^{m})\}$ related to sample $\omega^m$ by Eq. \eqref{eq:sec3.2.3-8}.}
    \STATE{Compute effective elastic tensor $\bar{\mathbb{A}}=\{\bar{a}_{ijkl}\}$ by Eq. \eqref{eq:aijkl}.}
    \STATE{Compute micromodulus $\bar c_0(\bm x,\bm \xi)$ by Eq. \eqref{eq:eff-microc}.}
  \end{algorithmic}
  
  \vspace{6pt}

  {\bf Algorithm IV: Computation of  displacement filed $\bm u(\bm x)$}\\
  {\bf Input:} geometry of $\bar{\Omega}$, statistical critical stretch $\bar{s}_{0}$, equivalent micromodulus $\bar{c}_0$\\
  {\bf Output:} displacement filed $\bm u(\bm x)$ in homogenized structure $\bar{\Omega}$
  \begin{algorithmic}[1]
    \STATE{Compute displacement of the homogenized structure $\bar{\Omega}$ by Eq. \eqref{eq:sec3.3-1}.}
    \STATE{Compute bond stretch by Eq. \eqref{eq:sec3.3-4} and judge bond failure by Eq. \eqref{eq:sec3.3-5}.}
  \end{algorithmic}
\end{algorithm}

\begin{figure}
  \centering
  \includegraphics[width=0.8\linewidth]{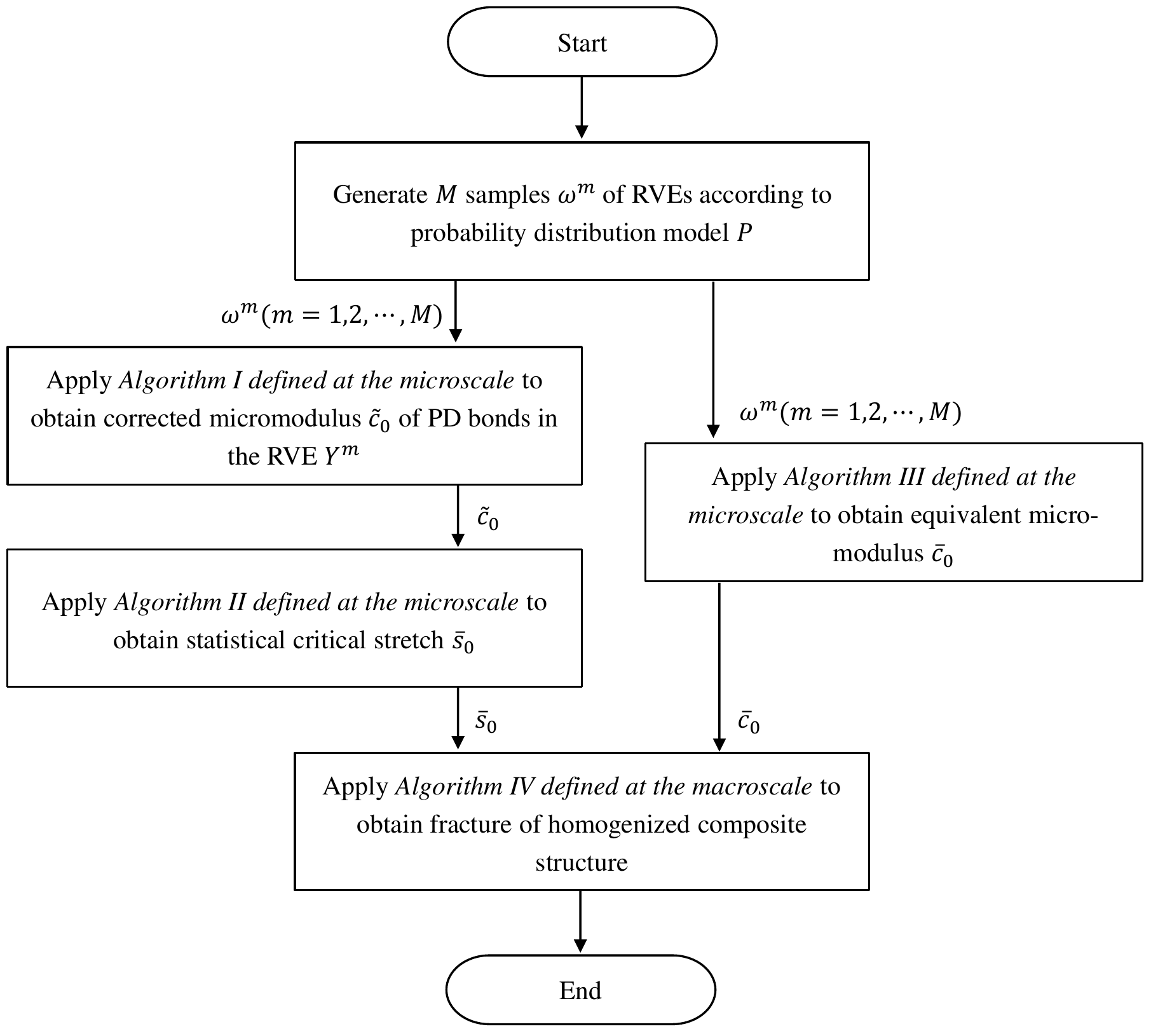}
  \caption{Flowchart for computing fracture of random composite structure.}
  \label{Fig:sec4.3-1}
\end{figure}

A detailed implementation procedure of the proposed PSM method for the fracture simulation of composite structure reinforced by randomly distributed particles is described in Algorithm of PSM.
The models are numerically implemented based on finite element discretization, in
which both the continuous elements (CE) and discrete elements (DE) are applied \cite{li2023peridynamics}.
The displacement field in the crack area belongs to $L^2$ space, and for the other area of the structure, the displacement field is in $H^1$ space.
According to the regularity of the structure, we apply the DEs in the crack area and CEs in the remainder of the structure.
Based on the computation of displacement, the bond stretch is obtained. Further, the initiation and propagation of crack are driven by the bond breaking according to bond failure criterion.
The update of crack area leads to an update of the meshes, which make some CEs transfer to DEs.
Besides,
\autoref{Fig:sec4.3-1} shows a flowchart of the proposed adaptive algorithm for fracture simulation.

\section{Numerical examples}
The feasibility and effectiveness of the proposed PSM method are illustrated with the 2D and 3D examples in this section.
We consider the brittle composites with Young's modulus of $E_M$ = 71.7 GPa for the matrix, and $E_P$ = 427 GPa for the particles.
Poisson's ratio for both matrix and particles are set to be 1/3 in 2D examples and 1/4 in 3D examples, respectively.
All 2D examples are meshed by quadrilateral elements and 3D cases are meshed by hexahedral elements.
The horizon, $\delta$, of neighborhoods in the PD model is defined as three times of average grid size to ensure the integral accuracy of nonlocal effects in numerical computations.
The critical stretches of particle bonds, matrix bonds and interface bonds are set to 0.00338, 0.01161 and 0.007495, respectively in Section 5.1, Section 5.2.1 and Section 5.2.2, and they are set to 0.00338, 0.01161 and 0.00387, respectively in Section 5.2.3, and 0.01, 0.01 and 0.01, respectively in Section 5.3.
2D and 3D RVEs are chosen as the square and cube domains with side length of 1mm, respectively, in all the examples.
The micromodulus of both matrix and particles in the RVE with a given sample $\omega^m$ is assumed to be an exponential function \cite{azdoud2013morphing}
\begin{equation}\label{eq:mirco-exa}
  c_0(\bm y, \bm \zeta, \omega^m) =\bigl(a_0+a_1\cos(2\theta)+a_2\cos(4\theta)\bigr) e^{-|\bm{\zeta}|/l}
\end{equation}
where $a_0$,  $a_1$ and $a_2$ can be estimated by the Poisson's ratio and Young's modulus of matrix and particles \cite{azdoud2013morphing}, $\theta$ denotes the angle between the bond $\bm \zeta$ and $y_1$-axis, and $l$ is a characteristic length set to be one-third of average gird size of the RVE.
Similarly, after the statistical homogenization, the equivalent micromodulus can also be defined as
\begin{equation}\label{eq:emirco-exa}
  \bar{c}_0(\bm x, \bm \xi) =\bigl(\bar{a}_0+\bar{a}_1\cos(2\varphi)+\bar{a}_2\cos(4\varphi)\bigr) e^{-|\bm{\xi}|/\bar{l}}
\end{equation}
where $\bar{a}_0$,  $\bar{a}_1$ and $\bar{a}_2$ can be estimated by the effective Poisson's ratio and Young's modulus of composite materials, $\varphi$ denotes the angle between the bond $\bm \xi$ and $x_1$-axis, and $\bar{l}$ is a characteristic length set to be one-third of average gird size of macroscale homogenized structures.
It should be pointed out that the homogenized materials in Section 5.2.1 and 5.2.3 are isotropic and $\bar{a}_1=\bar{a}_2=0$ \cite{azdoud2013morphing}.
Moreover, our numerical experiments are performed on a desktop workstation with 96G memory and 2.20GHz Xeon 5220R CPU.

\subsection{Validation of  accuracy and efficiency for the PSM method}

\begin{figure}[ht]
  \centering
  \includegraphics[scale=0.4]{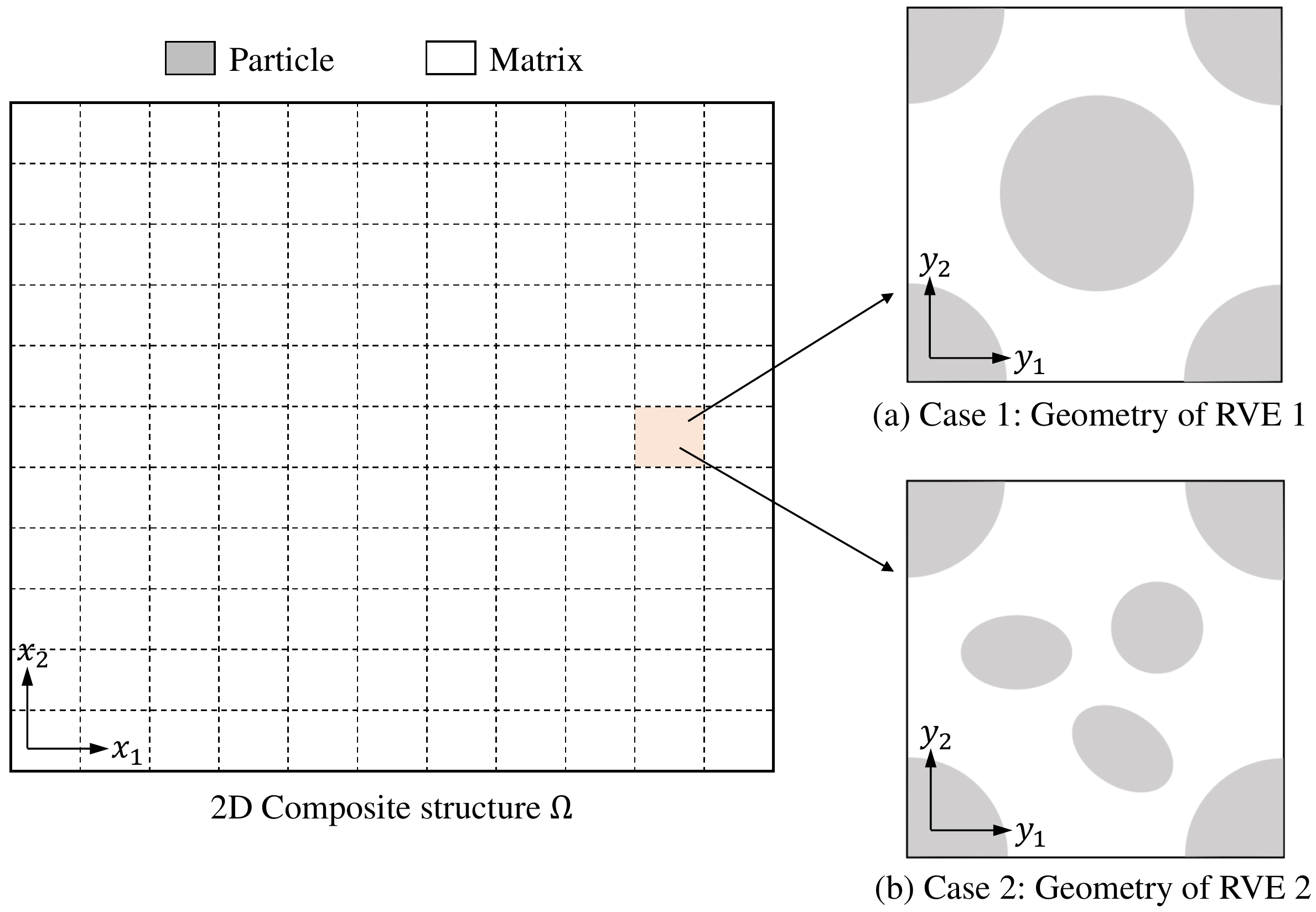}
  \caption{2D Macroscopic structure and related RVEs.}\label{Fig:sec5.1-1}
\end{figure}

\begin{figure}[H]
  \centering
  \subfigure[]{
    \includegraphics[width=0.29\textwidth]{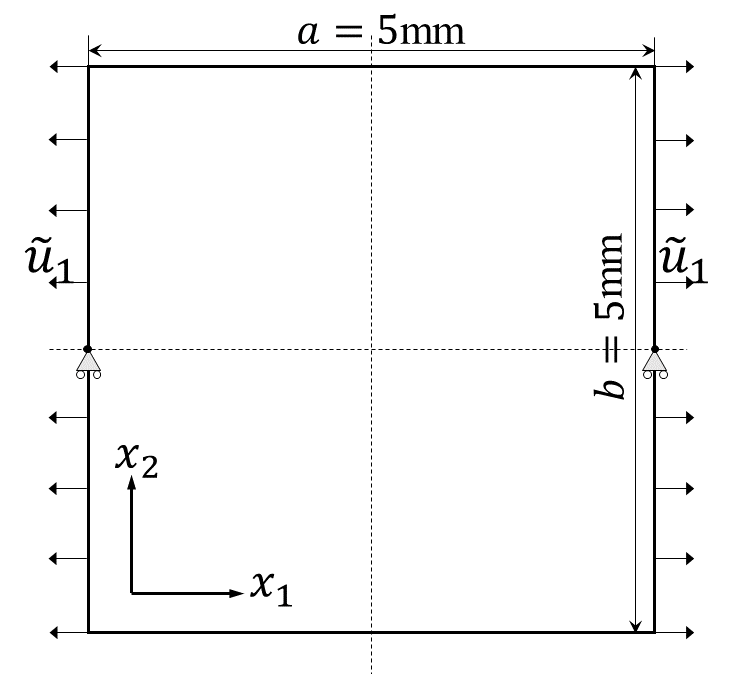}\label{fig:macro-5-1}}
  \subfigure[]{
    \includegraphics[width=0.29\textwidth]{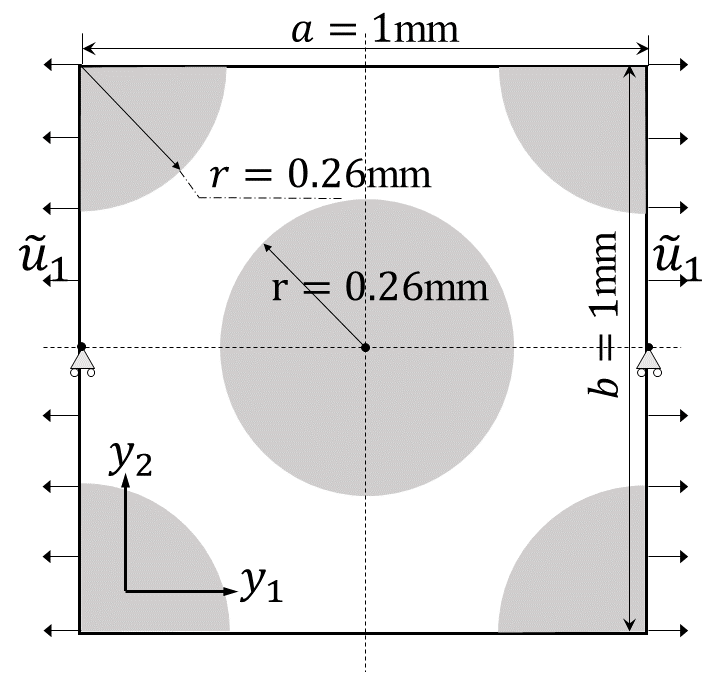}\label{fig:rve1_geo}}
  \subfigure[]{
    \includegraphics[width=0.29\textwidth]{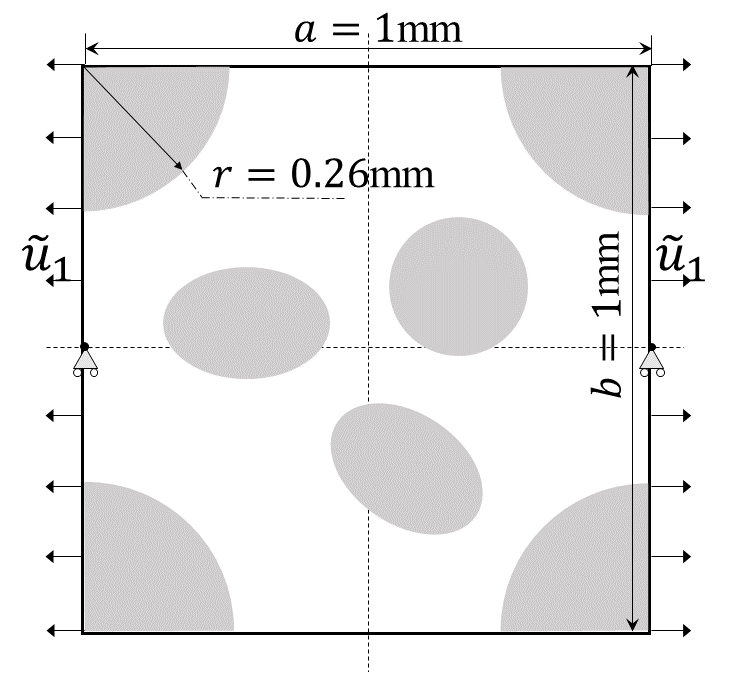}\label{fig:rve2_geo}}
  \caption{(a) geometry and boundary conditions for the macroscopic structure, (b) and (c) geometry and boundary conditions for RVE 1 and RVE 2.}\label{fig:sec5-1}
\end{figure}

To validate the proposed PSM method, we consider two kinds of periodic composite structures, namely, the square structures with side length of 5mm and composed of RVE 1 as illustrated in \autoref{Fig:sec5.1-1}(a), and RVE 2 as illustrated in \autoref{Fig:sec5.1-1}(b). \autoref{fig:sec5-1} shows the geometry and boundary conditions of the composite structures and corresponding RVEs.
The biaxial tensile boundary condition with $\tilde{u}_1=0.02$ mm along $y_1$-axis is specified on the left and right sides of the RVEs, and the middle points on the left and right sides of RVEs are fixed along $y_2$-axis. Similar biaxial tensile boundary condition with $\tilde{u}_1=0.03$ mm along $x_1$-axis and fixed boundary condition along $x_2$-axis are specified for the macroscopic structures.
In order to demonstrate the accuracy and efficiency of the present method, we take the single-scale direct BPD fracture analyses of composite structures as reference.
The BPD simulations of the composite structures, RVEs and homogenized structures were implemented through 100 displacement increments (steps).

\begin{figure}[H]
  \centering
  \subfigure[]{
    \includegraphics[width=0.33\textwidth]{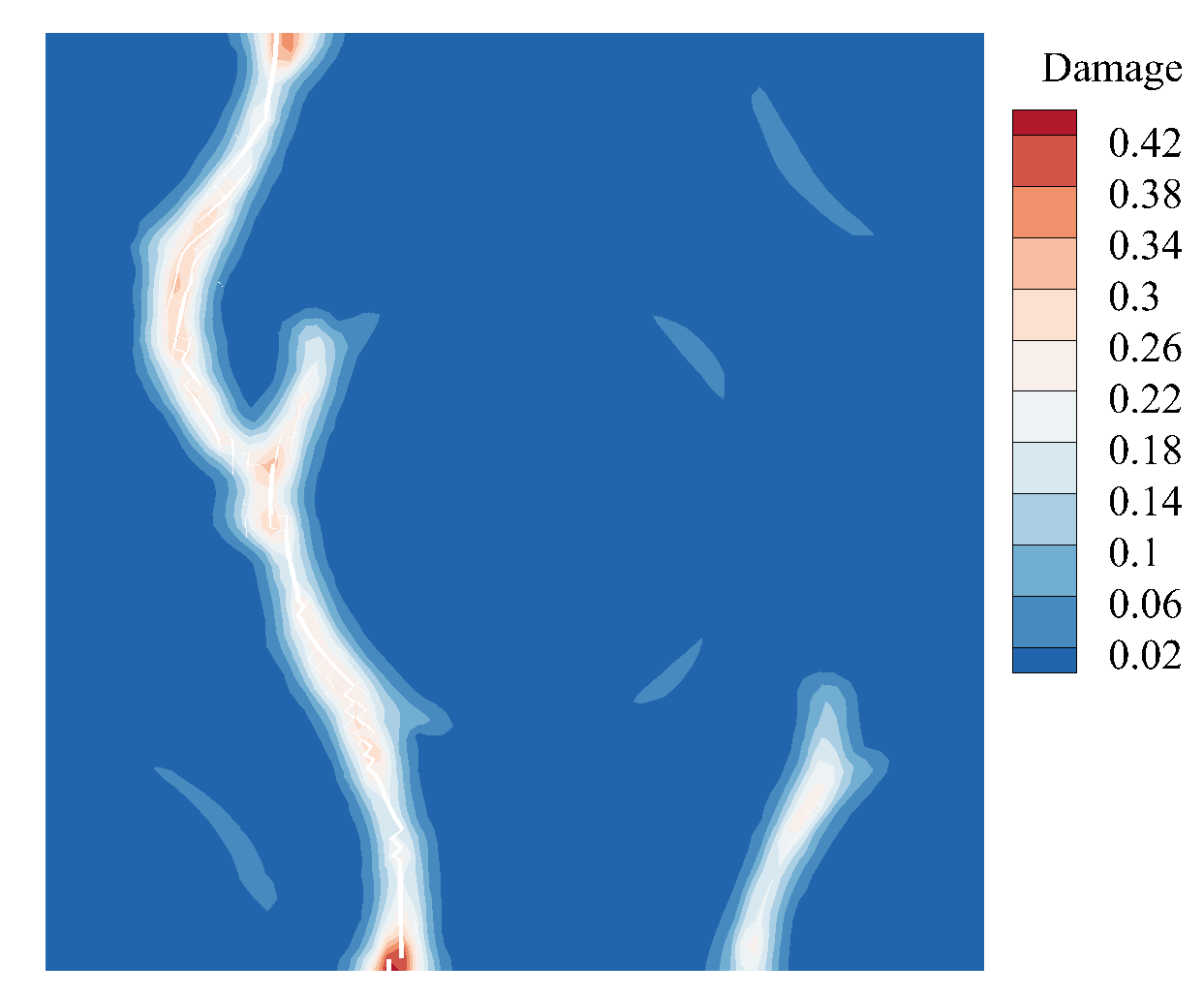}\label{fig:rve1_frac_cell}}
  \subfigure[]{
    \includegraphics[width=0.33\textwidth]{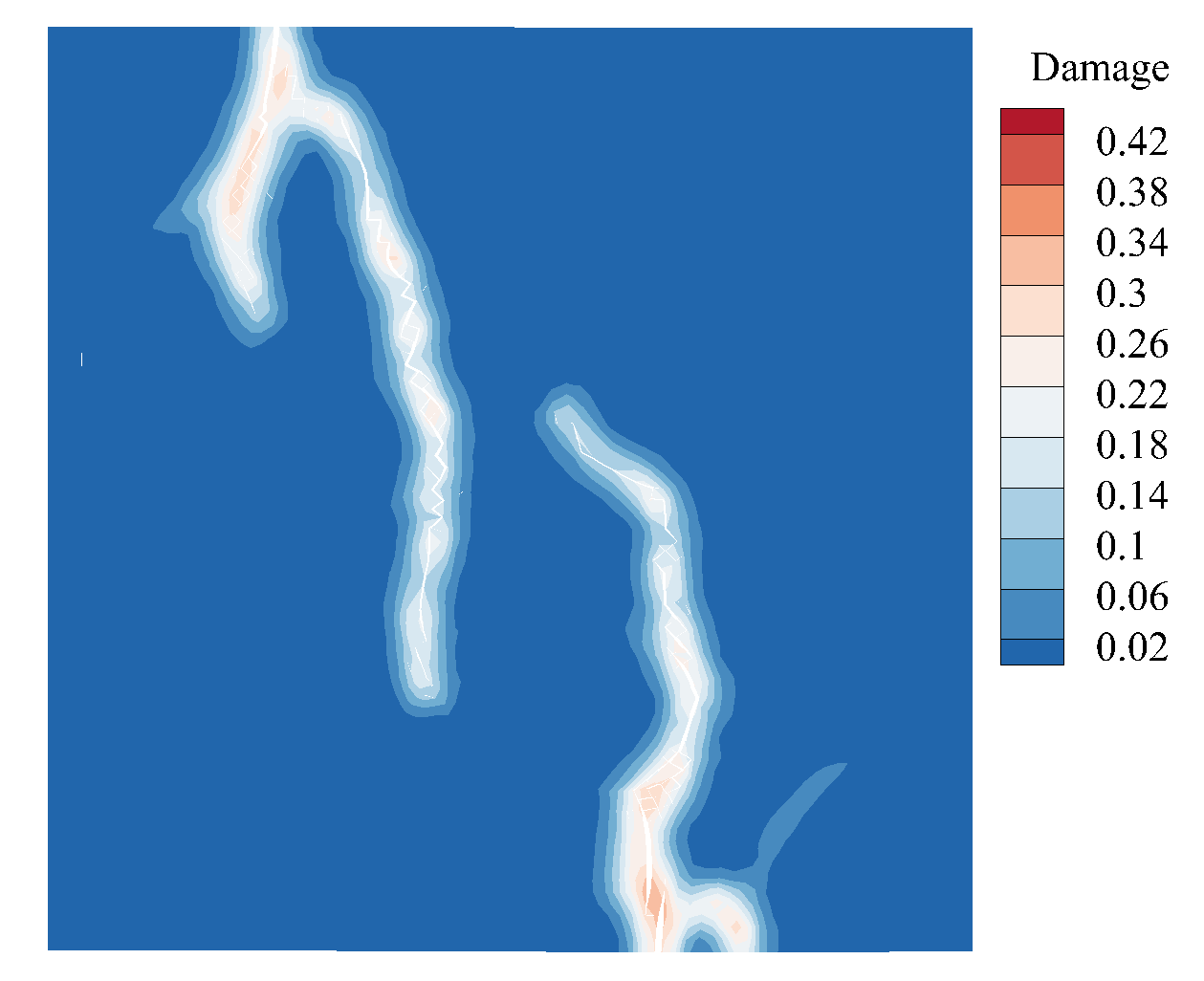}\label{fig:rve2_frac_cell}}
  \subfigure[]{
    \includegraphics[width=0.29\textwidth]{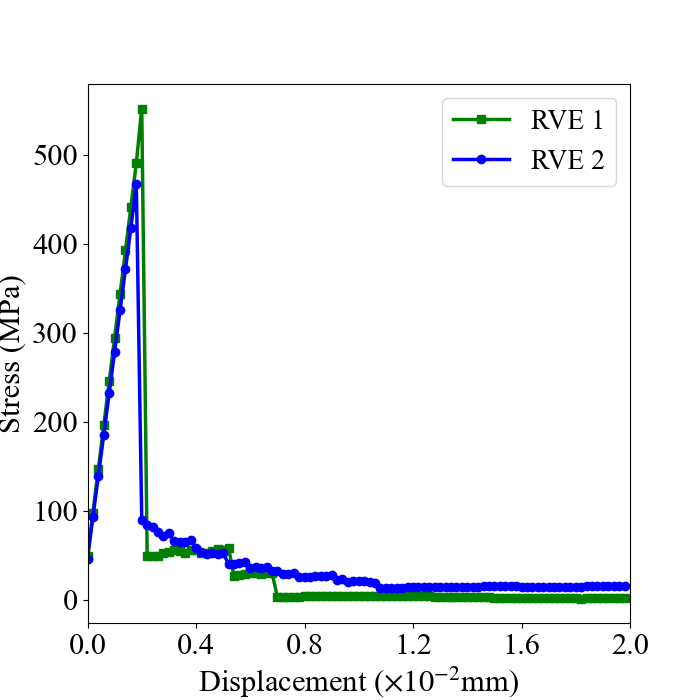}\label{fig:rve1_2_stress_cell}}
  \caption{Crack path for (a) 2D RVE 1 and (b) 2D RVE 2 at imposed displacement step 14 and 11, respectively; (c) stress versus imposed displacement for two kinds of RVEs.}\label{fig:rve_frac}
\end{figure}

\begin{figure}[H]
  \centering
  \subfigure[]{
    \includegraphics[width=0.33\textwidth]{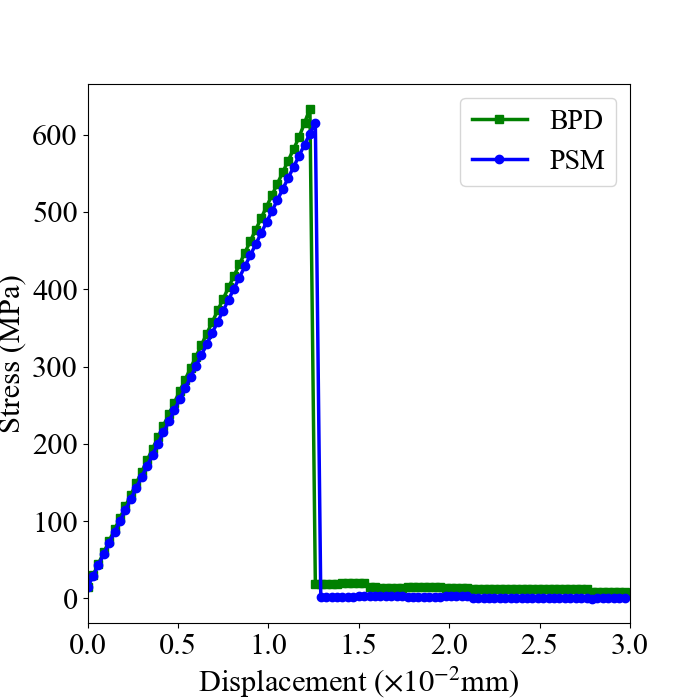}}
  \subfigure[]{
    \includegraphics[width=0.33\textwidth]{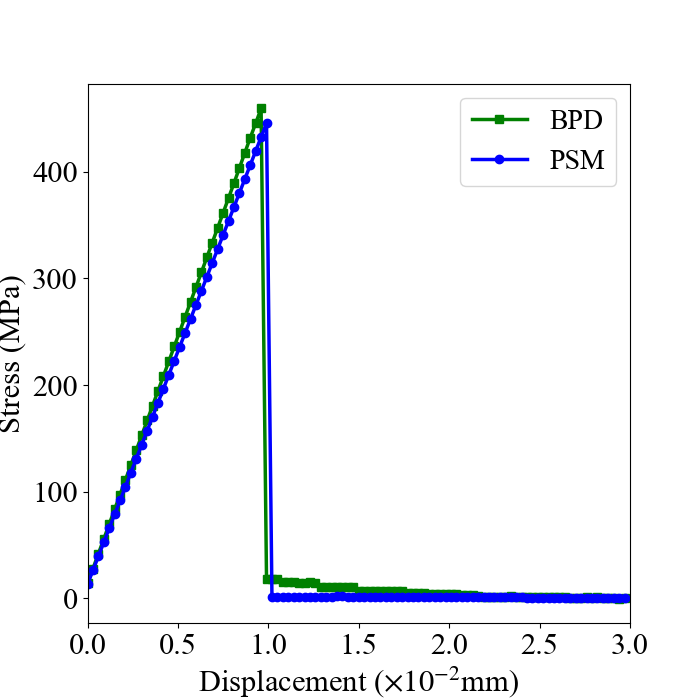}}
  \caption{Stress versus imposed displacement for composite structure with (a) RVE 1, $\varepsilon=1/5$ and (b) RVE 2, $\varepsilon=1/5$ simulated by PSM method and single-scale direct BPD method, respectively.}\label{fig:sec5_macro}
\end{figure}

\autoref{fig:rve_frac} shows the crack paths and the stress in the $y_1$ direction versus displacement of two RVEs. It can be found from \autoref{fig:rve_frac}(a) and \autoref{fig:rve_frac}(b) that the microstructure can significantly affect the crack propagation path, and RVE 1 and RVE 2 break completely at step 14 and 11, respectively. Besides, we can obtain from \autoref{fig:rve_frac}(c) that the effective critical stretches of two RVEs (see Eq. \eqref{eq:sec3.2.2-6}) are 0.0056 and 0.0044, respectively. Based on the homogenization method, the equivalent Young's modulus (see Eq. \eqref{eq:sec3.2.3-8}) are obtained as 118.286 GPa and 110.354 GPa, respectively. Thus, the fracture of homogenized macroscopic structure can be simulated by the BPD model with effective critical bond stretch and micromodulus calculated by equivalent Young's modulus.
\autoref{fig:sec5_macro} shows the stress versus displacement for two kinds of composite structures simulated by the PSM method and single-scale direct BPD analyses.
From \autoref{fig:sec5_macro}, it can be seen that the failure process simulated by the PSM method is in good agreement with that by the single-scale direct BPD analyses.
The information of CE and DE meshes for the composite structures, RVEs and homogenized structures is listed in \autoref{tab:time_compare},
which shows the advantage of the proposed PSM method by saving almost 80\% computational time and it is very important in engineering computation.
As a result, the PSM method is effective to analyze the fracture behaviors of particle reinforced composite structures.

\begin{table}[H]
  \setlength{\abovecaptionskip}{0cm}
  \setlength{\belowcaptionskip}{0.3cm}
  \centering
  \caption{Comparison of computational time for 2D composite structures, RVEs and homogenized structures.}\label{tab:time_compare}
  \scalebox{0.8}{
    \begin{tabular}{cccccccc}
      \toprule  
      \multirow{2}*{} & \multicolumn{3}{c}{Composite with RVE 1} & \quad  & \multicolumn{3}{c}{Composite with RVE 2}                                                      \\ \cline{2-4} \cline{6-8}
                      & Composite                                & RVE 1  & Homogenized structure                    & \quad & Composite & RVE 2  & Homogenized structure \\ \hline
      Elements        & 77,296                                   & 5,820  & 28,224                                   & \quad & 91,128    & 5,770  & 28,224                \\
      CE Nodes        & 77,857                                   & 5,973  & 28,561                                   & \quad & 91,729    & 5,923  & 28,561                \\
      DE Nodes        & 309,184                                  & 23,280 & 112,896                                  & \quad & 364,512   & 23,080 & 112,896               \\
      Time(s)         & 1,191,714                                & 17,047 & 242,212                                  & \quad & 1,263,596 & 52,318 & 346,950               \\
      \bottomrule  
    \end{tabular}
  }
\end{table}

\subsection{Fracture of 2D random composite structures by the PSM method}

\subsubsection{Composite structures with uniform distribution of particles}

\begin{figure}[H]
  \centering
  \includegraphics[scale=0.48]{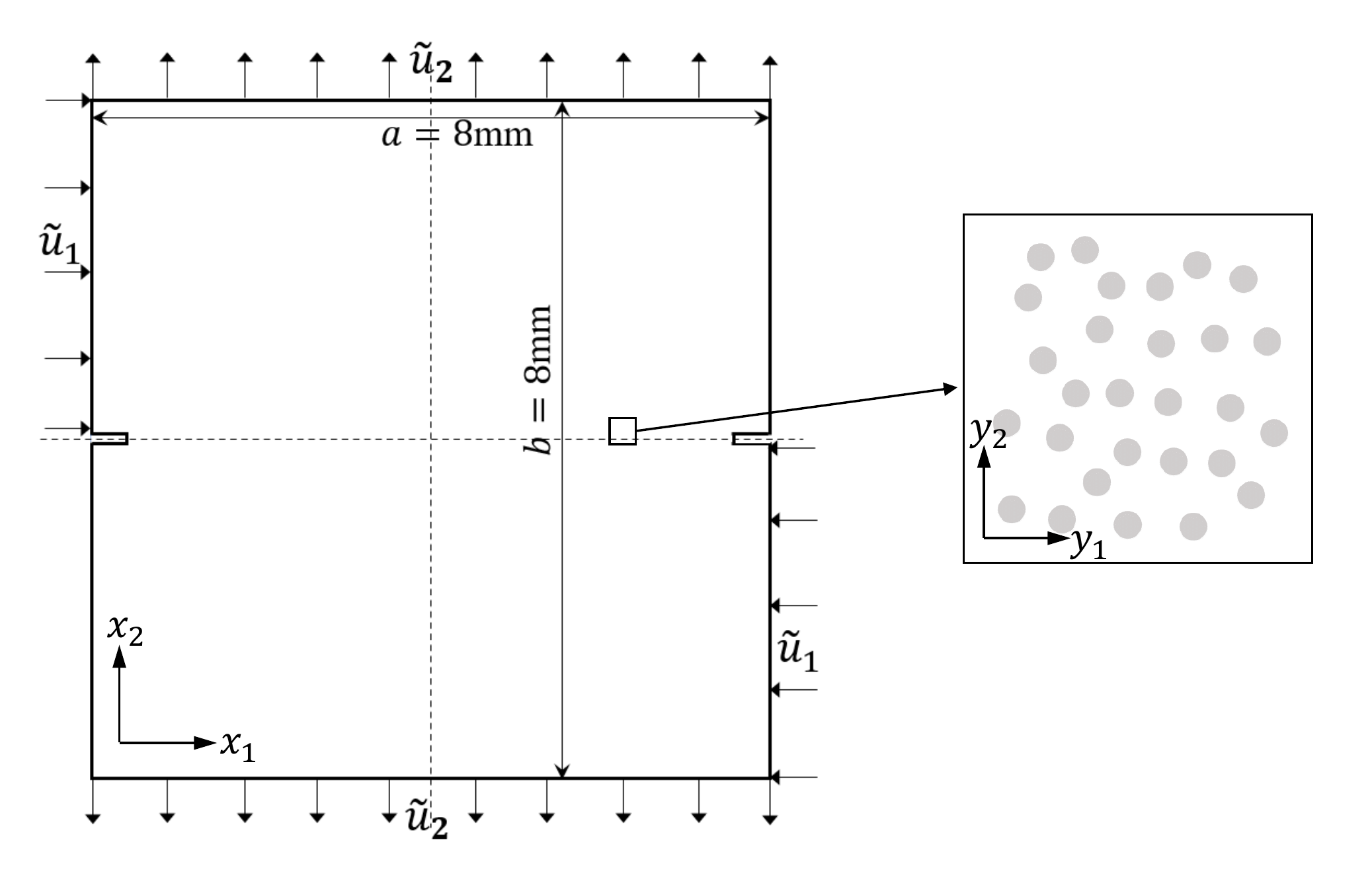}
  \caption{Geometry and boundary conditions of the structure with uniform distribution of  circular particles.}\label{fig:sec5.2-1}
\end{figure}

We consider the composite structure with side length of 8mm and reinforced by uniform distribution of circular particles, as shown in \autoref{fig:sec5.2-1}, where the notches with width of 0.2mm and length of 1mm are preset at the left and right sides of the structure, and the shear and tensile boundary conditions with $\tilde{u}_1=\tilde{u}_2=0.1$mm are given.
The volume fraction of particles in the RVE is 14\%. The size and boundary condition of the RVE are the same as those in Section 5.1.
The BPD simulations of the RVE and homogenized structures were implemented through 100 displacement increments (steps).

\begin{figure}[H]
  \centering
  \subfigure[]{
    \includegraphics[width=0.4\textwidth]{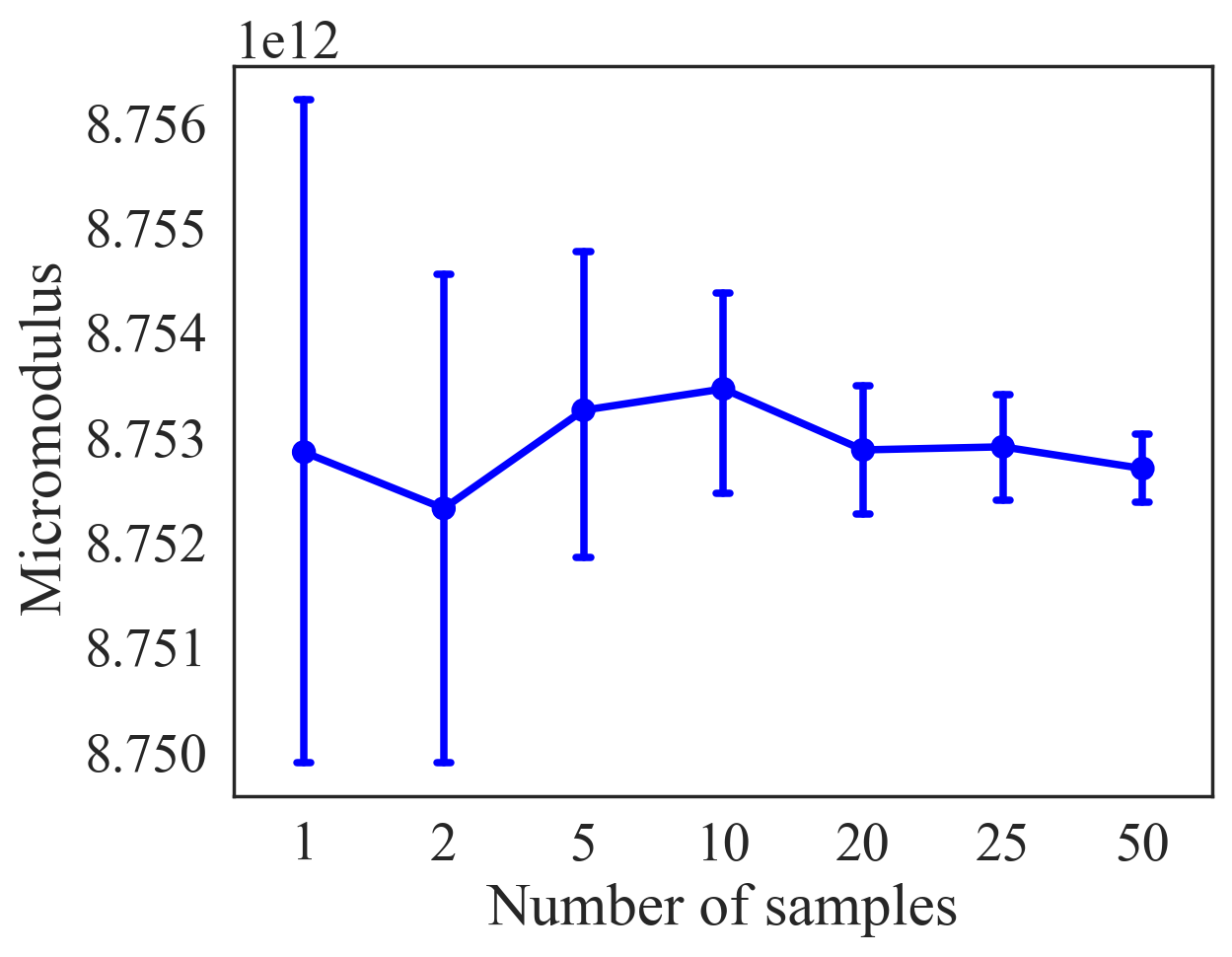}\label{fig:uniform_plot_c}}
  \subfigure[]{
    \includegraphics[width=0.4\textwidth]{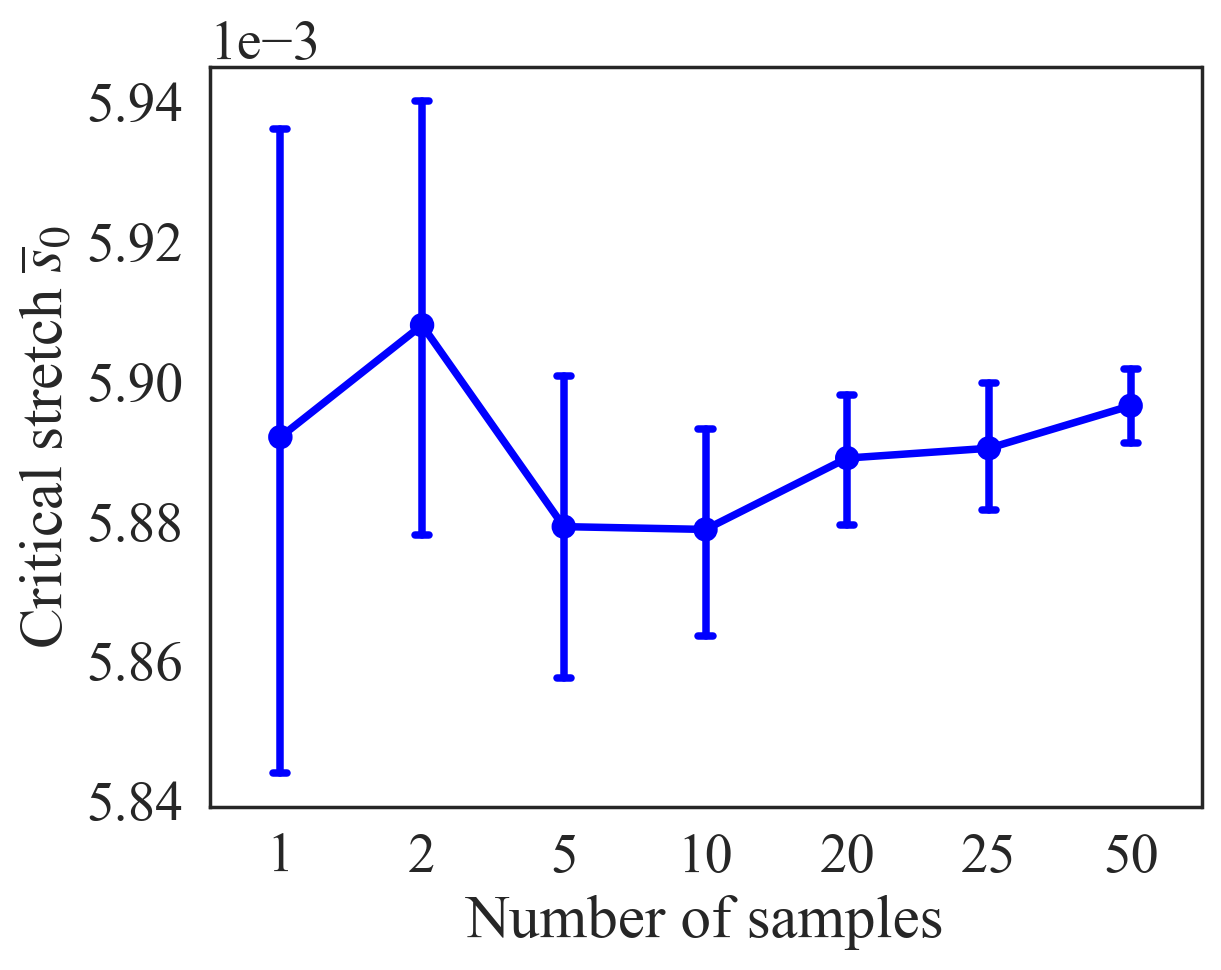}\label{fig:uniform_plot_t}}
  \caption{The expected value of (a) $\bar{a}_0$ in micromodulus defined by Eq. \eqref{eq:emirco-exa}, (b) critical stretch with different number of samples.}\label{fig:uniform_plot}
\end{figure}

\begin{figure}[H]
  \centering
  \subfigure[Step 12]{
    \includegraphics[width=0.32\textwidth]{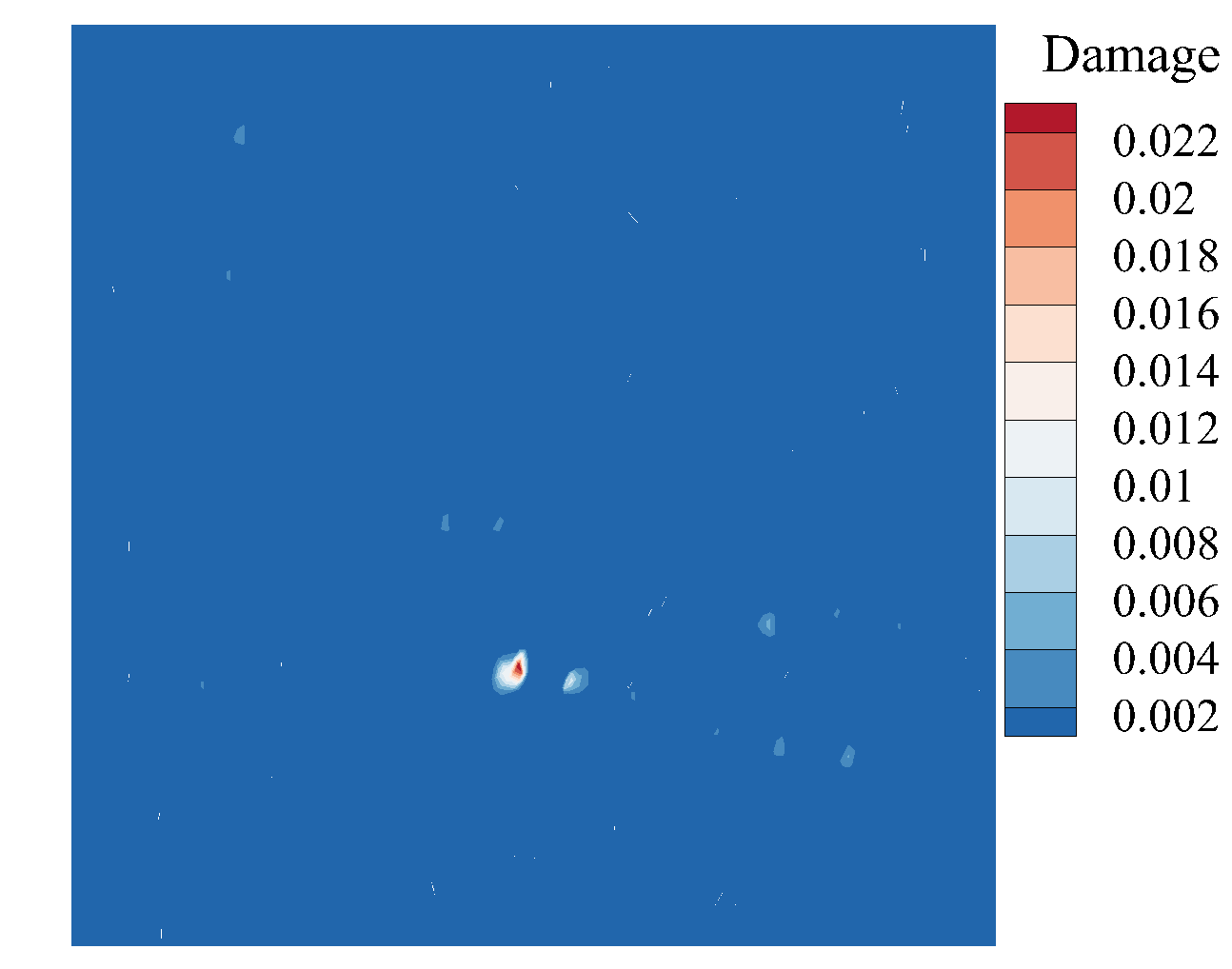}}
  \subfigure[Step 16]{
    \includegraphics[width=0.32\textwidth]{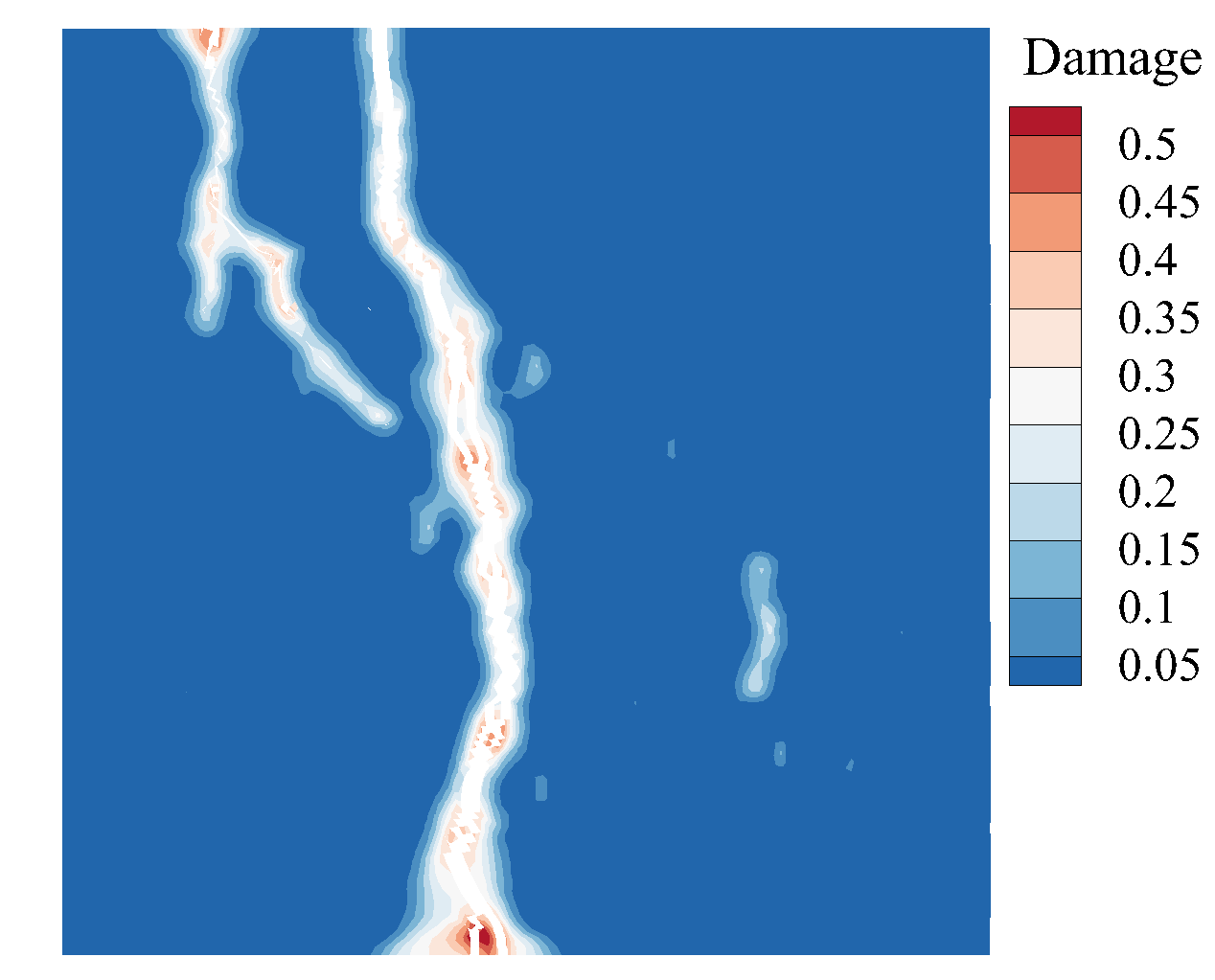}}
  \subfigure[Step 100]{
    \includegraphics[width=0.32\textwidth]{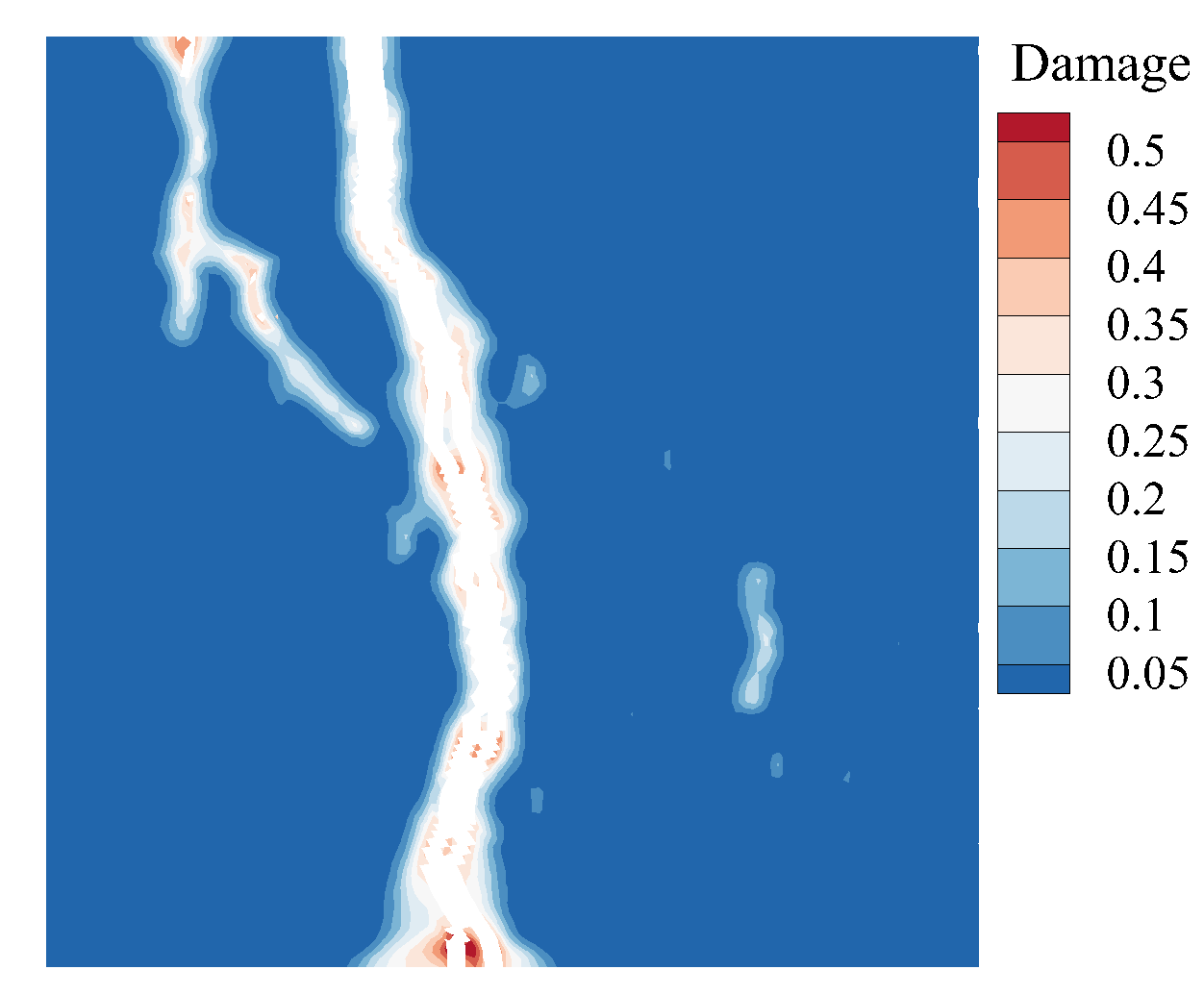}}
  \caption{Crack paths of the RVE at different steps.}\label{fig:frac_uniform}
\end{figure}

\begin{figure}[H]
  \centering
  \subfigure[Step 8]{
    \includegraphics[width=0.32\textwidth]{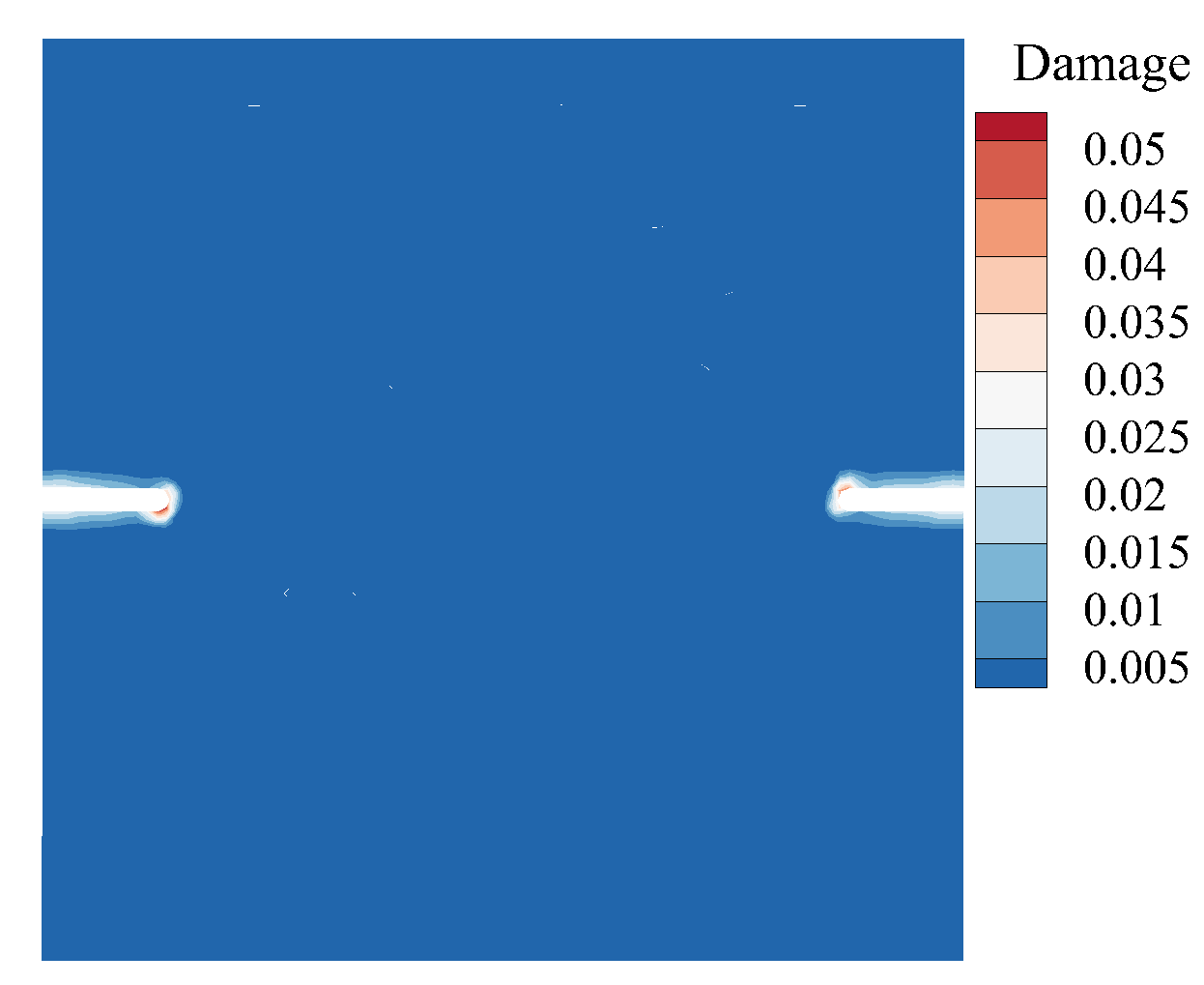}}
  \subfigure[Step 10]{
    \includegraphics[width=0.32\textwidth]{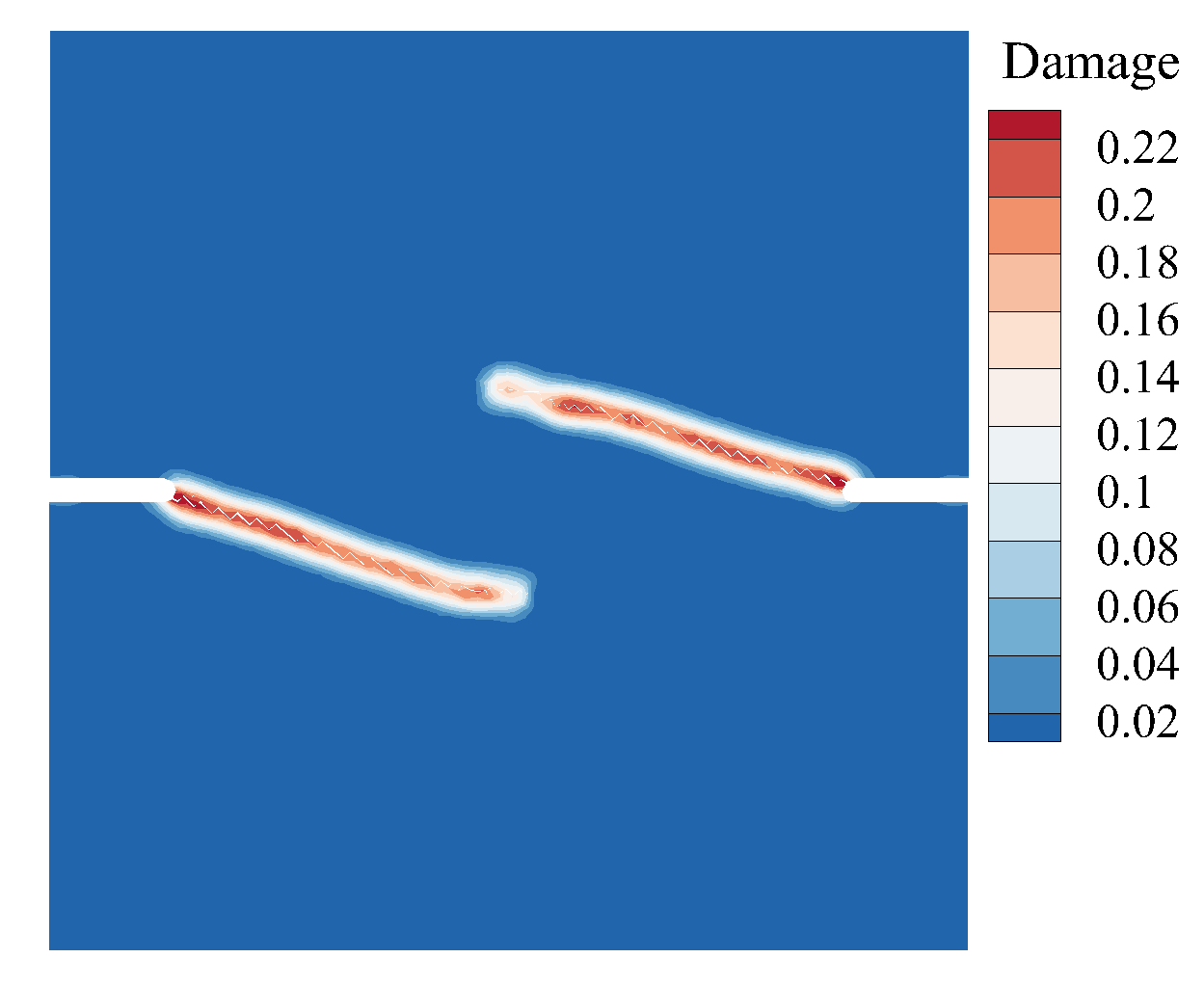}}
  \subfigure[Step 100]{
    \includegraphics[width=0.32\textwidth]{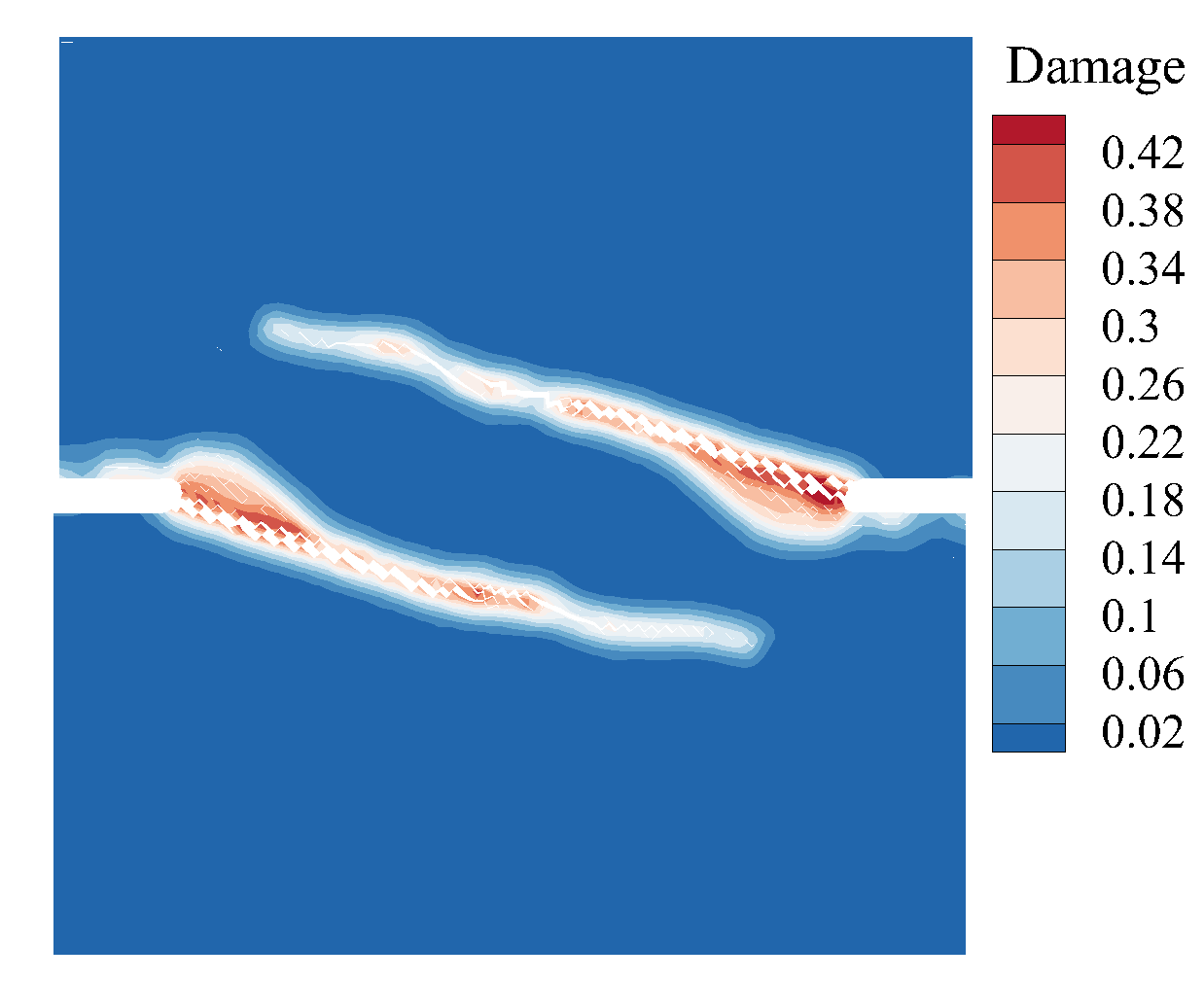}}
  \caption{Crack paths of the homogenized structure under biaxial tension and shear at different steps.}\label{fig:frac_biaxial}
\end{figure}

\autoref{fig:uniform_plot} displays the expected values of $\bar{a}_0$ in micromodulus defined by Eq. \eqref{eq:emirco-exa} and critical stretch of the random composite material.
Statistically, different samples
have different results. But accompanied by the increasing number of samples with the same statistical characteristic, the mathematical expectation of computational results
should converge. Obviously, as shown in \autoref{fig:uniform_plot}, the scatter of data
decreases with the increasing number of samples.
Therefore, 25 samples were taken in this study to avoid an unacceptable scatter of numerical results.
\autoref{fig:frac_uniform} shows the crack paths in a specified sample RVE simulated by the microscopic BPD model.
\autoref{fig:frac_biaxial}  show the crack paths for the homogenized structure under biaxial tension and shear conditions simulated by the macroscopic BPD model with the expectations of micromodulus and critical stretch, which is consistent with the experimental results \cite{nooru1993experimental}.
Based on the proposed PSM framework, a single-scale direct BPD simulation of the composite structure can be replaced by the microscopic BPD simulation of RVEs and the macroscopic BPD simulation of homogenized structure with high efficiency.

\subsubsection{Composite structures reinforced by particles with different inclination angles}
\begin{figure}[H]
  \centering
  \includegraphics[scale=0.35]{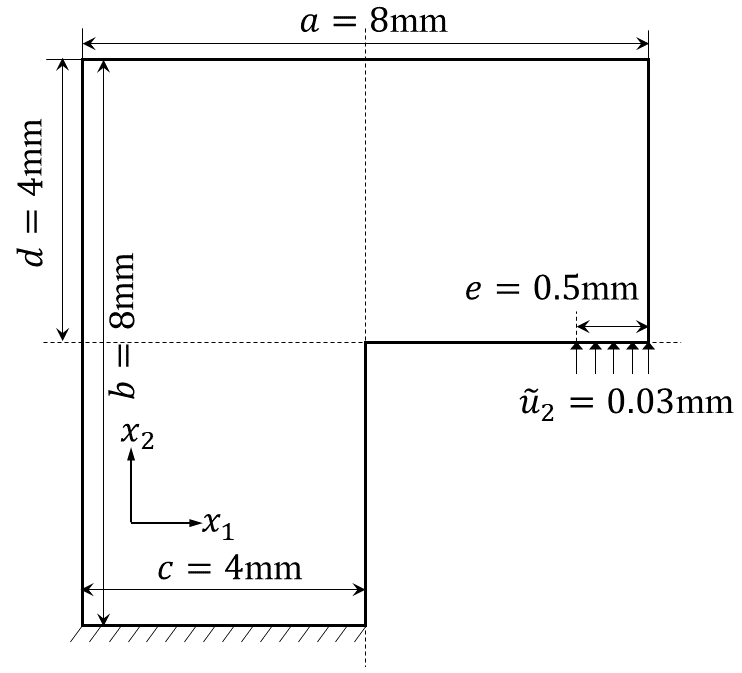}
  \caption{Geometry and boundary conditions of the L-shaped composite panel.}\label{fig:sec5.2-2}
\end{figure}

\begin{figure}[H]
  \centering
  \subfigure[$15^\circ$]{
    \includegraphics[width=0.3\textwidth]{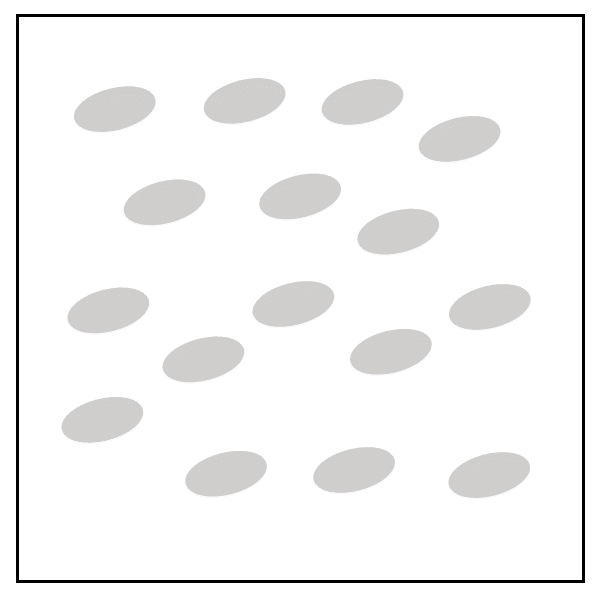}\label{fig:angle_15}}
  \subfigure[$45^\circ$]{
    \includegraphics[width=0.3\textwidth]{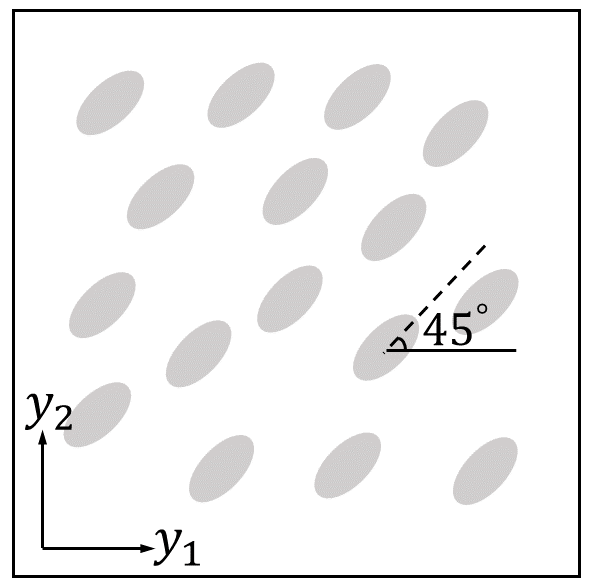}\label{fig:angle_45}}
  \subfigure[$60^\circ$]{
    \includegraphics[width=0.3\textwidth]{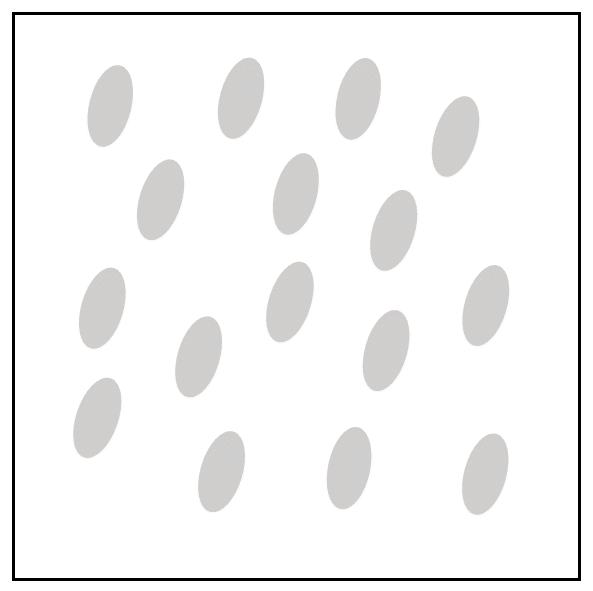}\label{fig:angle_75}}
  \caption{Geometry of RVEs with uniform distribution of ellipsoidal particles with different inclination angles.}\label{fig:angle_geo}
\end{figure}

We consider the L-shaped composite panel reinforced by uniform distribution of elliptical particles, and the boundary conditions
and dimensions are presented in \autoref{fig:sec5.2-2}.
Five RVEs containing 14\% of elliptical particles with different inclination angles between the long axis of ellipse and $y_1$-axis, including $15^\circ$, $30^\circ$, $45^\circ$, $60^\circ$ and $75^\circ$ are taken into account as shown in \autoref{fig:angle_geo}. The size and boundary condition of the RVEs are the same as those in Section 5.1.
The simulations of the RVEs and homogenized structures were implemented by 100 displacement increments (steps).

\begin{figure}[H]
  \centering
  \subfigure[$15^\circ$]{
    \includegraphics[width=0.32\textwidth]{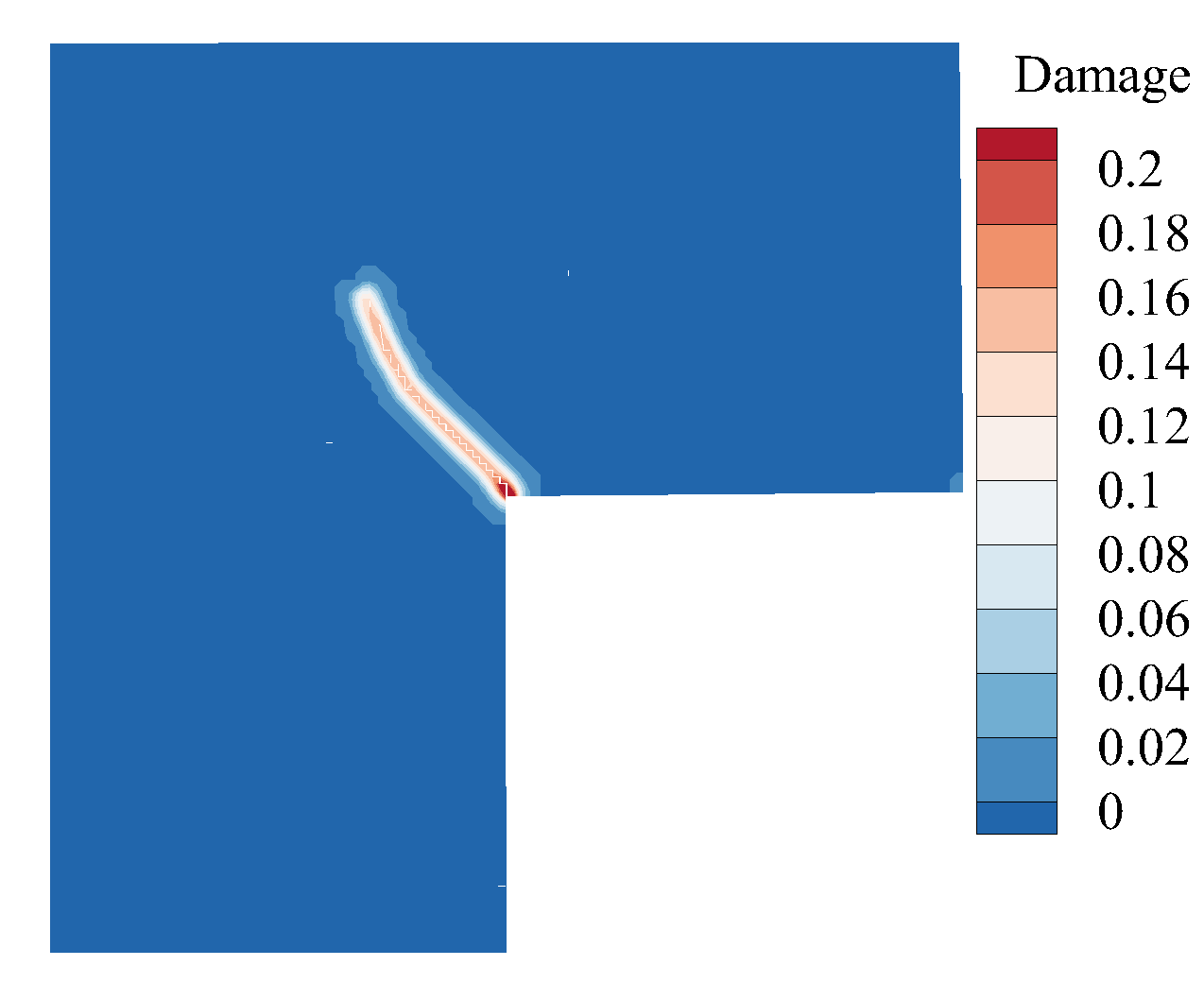}\label{fig:L_15}}
  \subfigure[$30^\circ$]{
    \includegraphics[width=0.32\textwidth]{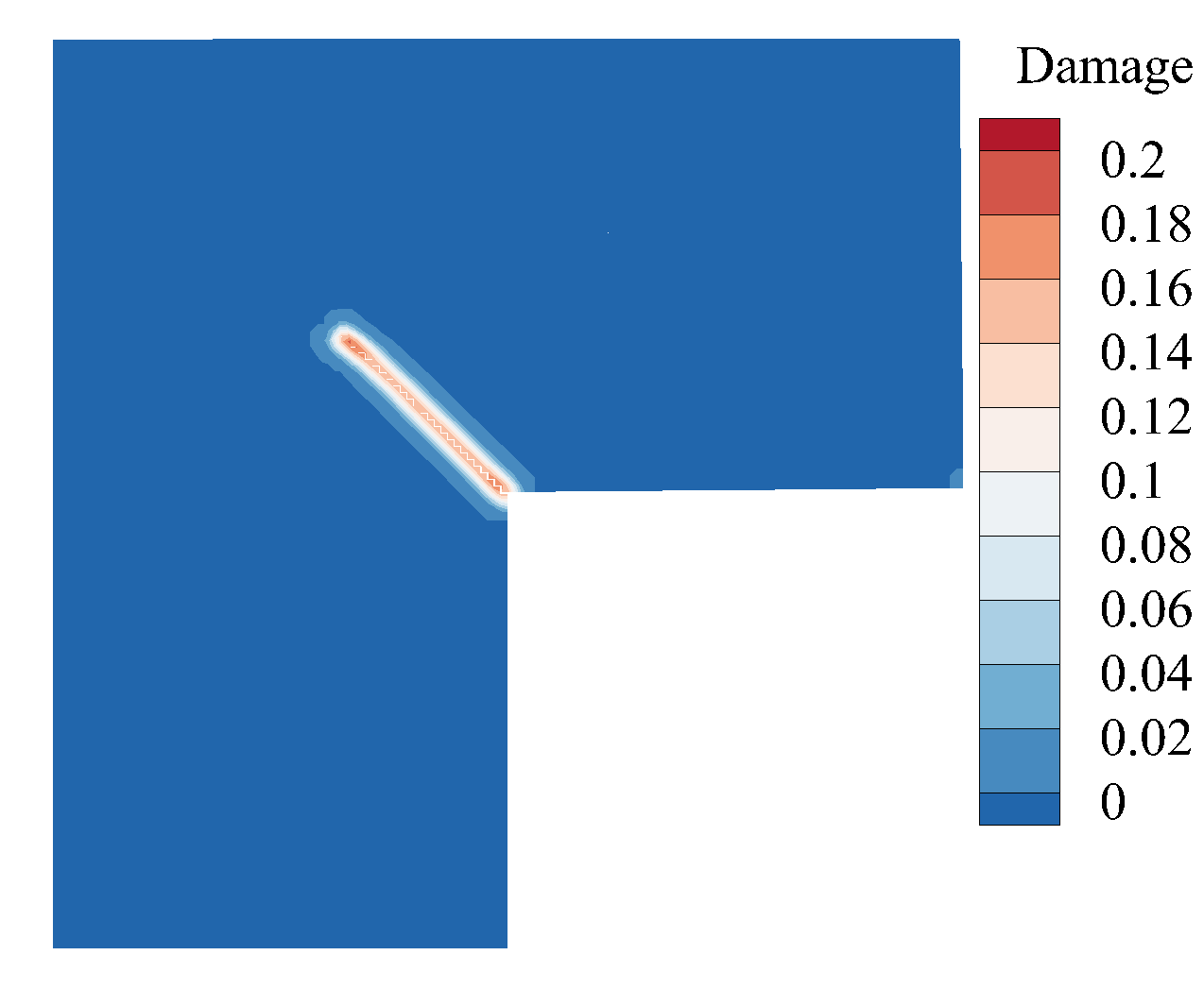}\label{fig:L_30}}
  \subfigure[$45^\circ$]{
    \includegraphics[width=0.32\textwidth]{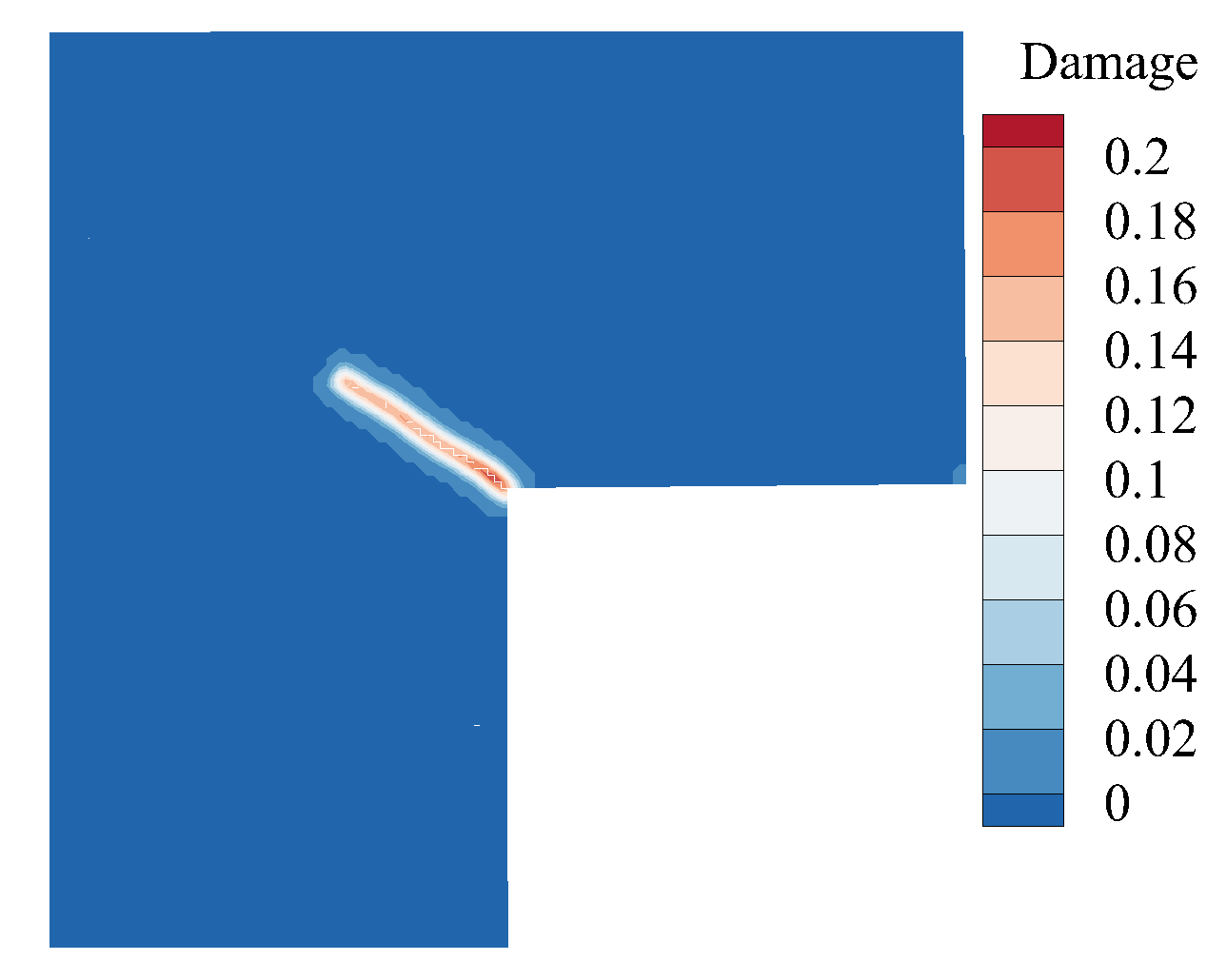}\label{fig:L_45}}
  \subfigure[$60^\circ$]{
    \includegraphics[width=0.32\textwidth]{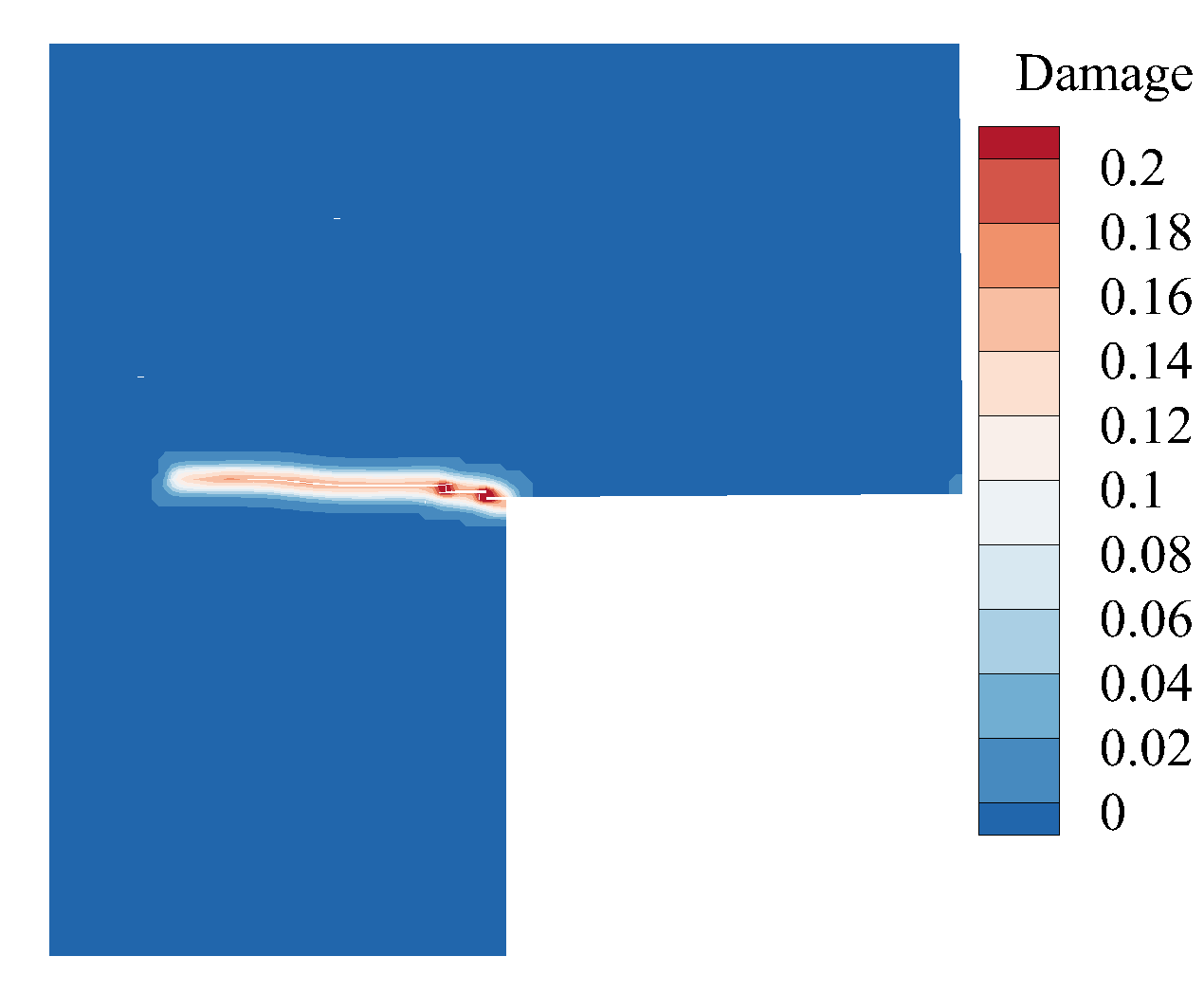}\label{fig:L_60}}
  \subfigure[$75^\circ$]{
    \includegraphics[width=0.32\textwidth]{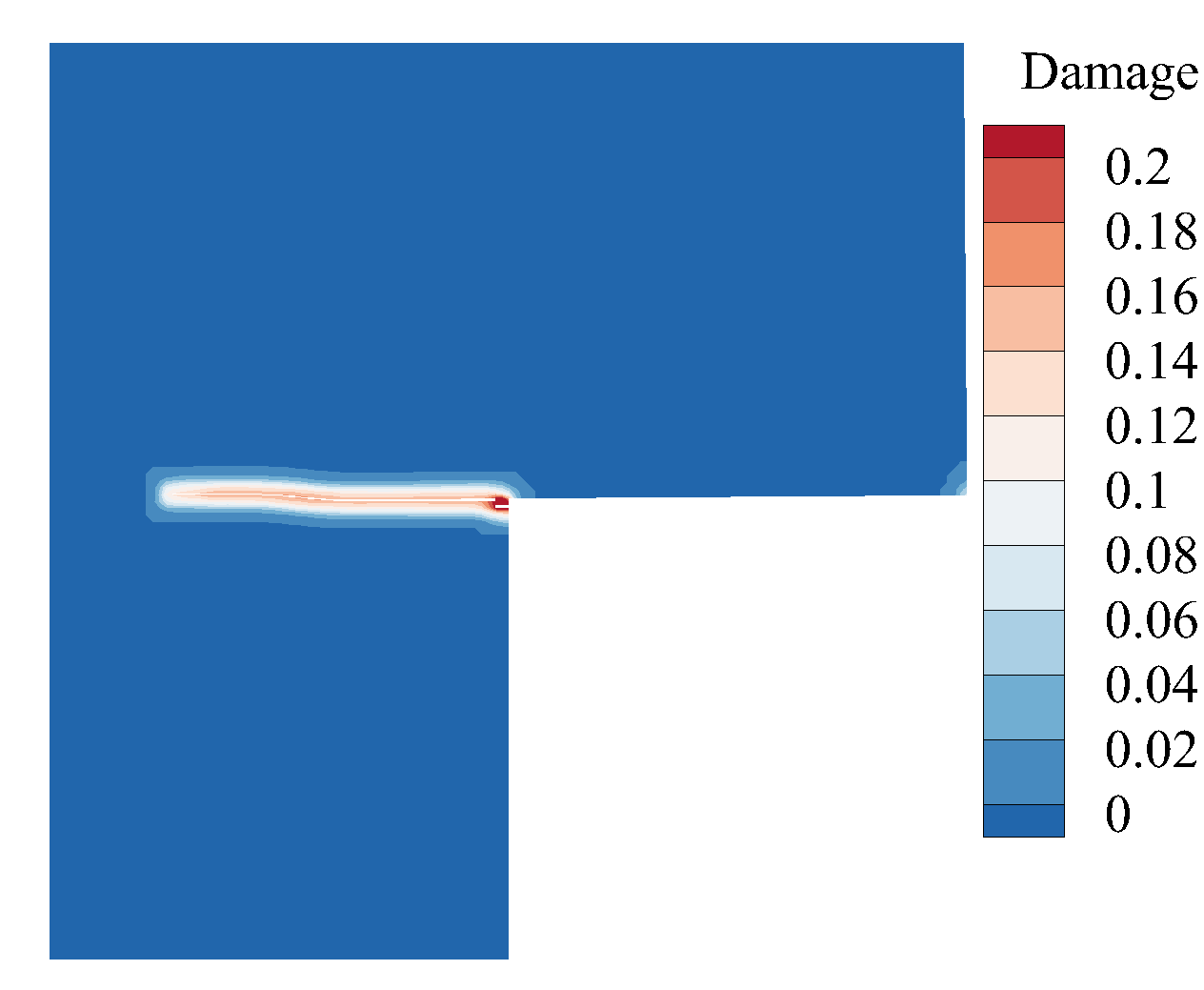}\label{fig:L_75}}
  \caption{Crack paths of homogenized L-shaped panel structures corresponding to RVEs with different inclination angles.}\label{fig:L_frac}
\end{figure}

\autoref{tab:angle_rve} presents the expected values of micromodulus defined by Eq. \eqref{eq:emirco-exa} and critical stretch of elliptical particle reinforced composites with inclination angle from $15^\circ$ to $75^\circ$.
From \autoref{tab:angle_rve}, it can be seen that the critical stretch $\bar{s}_{01}$ increases with the increase of inclination angle.
Since the biaxial tensile direction is parallel to $y_1$-axis, the smaller the inclination angle between the long axis of ellipse and $y_1$-axis is, the easier the stress concentration around ellipsoidal particles in the RVE is. Therefore, for the case of  inclination angle $15^\circ$, the cracks are most likely to initiate, which leads to the smallest critical stretch. The similar results can be obtained for the critical stretch $\bar{s}_{02}$ from \autoref{tab:angle_rve}.
The crack paths of the homogenized L-shaped panels corresponding to the composites with five RVEs are shown in \autoref{fig:L_frac}, and the  loading steps when crack initiates are 61, 67, 65, 65 and 63, respectively, which demonstrates that the shapes of particles in RVEs have the significant impact on the fracture of macroscopic structures.
Besides, the crack paths are similar to the results described in \cite{winkler2001experimental,narayan2019gradient}. Therefore, this example illustrates the effectiveness of the
PSM method.

\begin{table}[H]
  \setlength{\abovecaptionskip}{0cm}
  \setlength{\belowcaptionskip}{0.3cm}
  \centering
  \caption{The coefficients of equivalent micromodulus defined by Eq. \eqref{eq:emirco-exa} and critical stretch of ellipsoidal particle reinforced composites with different inclination angles.}\label{tab:angle_rve}
  \scalebox{0.8}{
    \begin{tabular}{cccccc}
      \toprule  
                     & $15^\circ$            & $30^\circ$            & $45^\circ$            & $60^\circ$            & $75^\circ$            \\ \hline
      $\bar{a}_{0}$  & $1.91\times 10^{18}$  & $1.88\times 10^{18}$  & $1.86\times 10^{18}$  & $1.84\times 10^{18}$  & $1.82\times 10^{18}$  \\
      $\bar{a}_{1}$  & $1.72\times 10^{15}$  & $1.69\times 10^{15}$  & $1.67\times 10^{15}$  & $1.65\times 10^{15}$  & $1.63\times 10^{15}$  \\
      $\bar{a}_{2}$  & $-4.38\times 10^{15}$ & $-4.31\times 10^{15}$ & $-4.25\times 10^{15}$ & $-4.21\times 10^{15}$ & $-4.16\times 10^{15}$ \\
      $\bar{s}_{01}$ & 0.0050                & 0.0056                & 0.0058                & 0.0062                & 0.0068                \\
      $\bar{s}_{02}$ & 0.0070                & 0.0066                & 0.0060                & 0.0056                & 0.0048                \\
      \bottomrule  
    \end{tabular}
  }
\end{table}

\subsubsection{Composite structures with different volume fractions of particles}
\begin{figure}[H]
  \centering
  \includegraphics[width=0.48\textwidth]{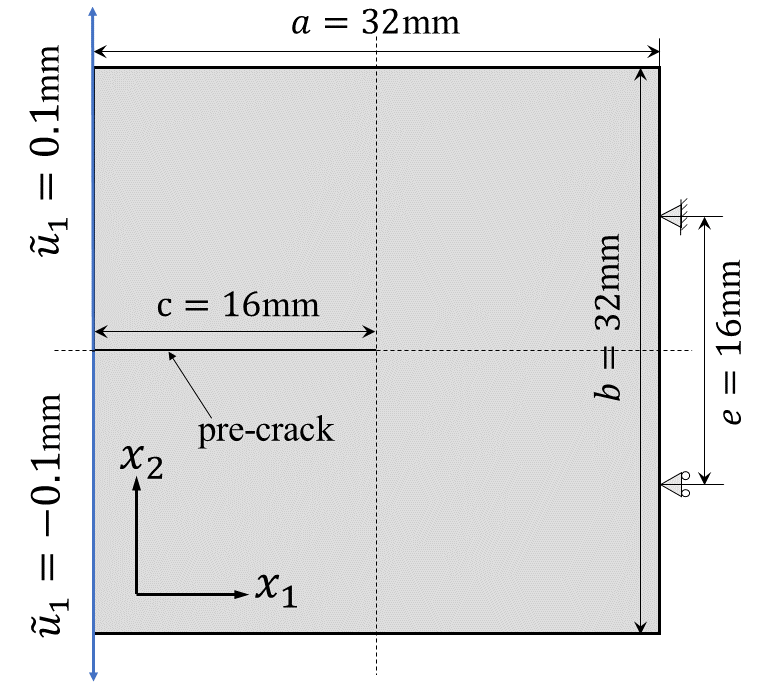}
  \caption{Geometry and boundary conditions of the structure with different volume fractions of particles.}\label{fig:volume_macro}
\end{figure}

\begin{figure}[H]
  \centering
  \subfigure[10\%]{
    \includegraphics[width=0.35\textwidth]{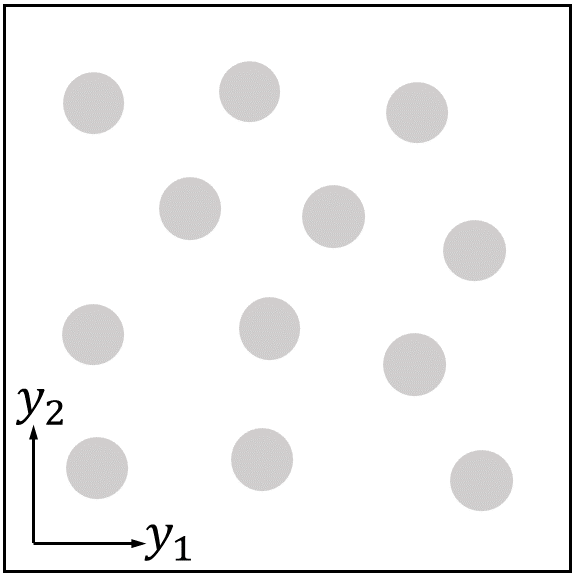}\label{fig:volume_10}}
  \subfigure[24\%]{
    \includegraphics[width=0.35\textwidth]{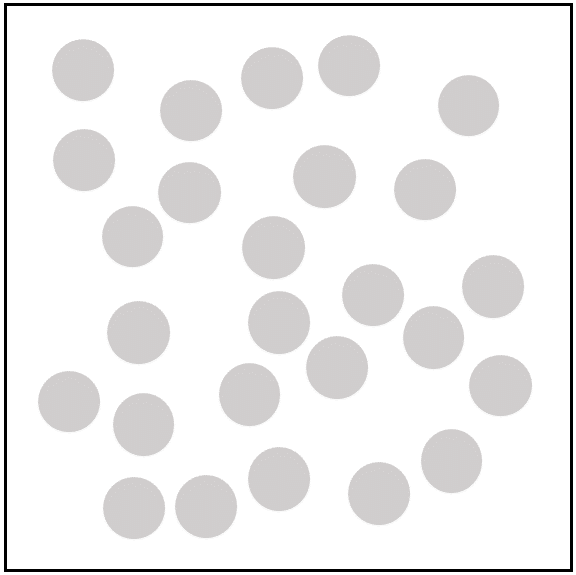}\label{fig:volume_22}}
  \caption{Geometry of RVEs with uniform distribution of circular particles with different volume fractions.}\label{fig:volume_geo}
\end{figure}

We consider the composite structures with the geometry and boundary conditions as shown in \autoref{fig:volume_macro}.
The composite structures reinforced by uniform distribution of circular particles with the radius of 0.0522mm and ten volume fractions of particles are studied, including from 6\% to 24\% with an interval of 2\%, as shown in \autoref{fig:volume_geo}. Twenty samples are selected for the RVEs of each volume fraction.
The biaxial tensile boundary condition with $\tilde{u}_1=0.003$ mm along $y_1$-axis is preset on the left and right sides of RVEs, and the middle points on the left and right sides of RVEs are fixed along $y_2$-axis.
The simulations of the RVEs and homogenized structures were implemented by 100 and 166 displacement increments (steps), respectively.

\begin{figure}[H]
  \centering
  \subfigure[]{
    \includegraphics[width=0.45\textwidth]{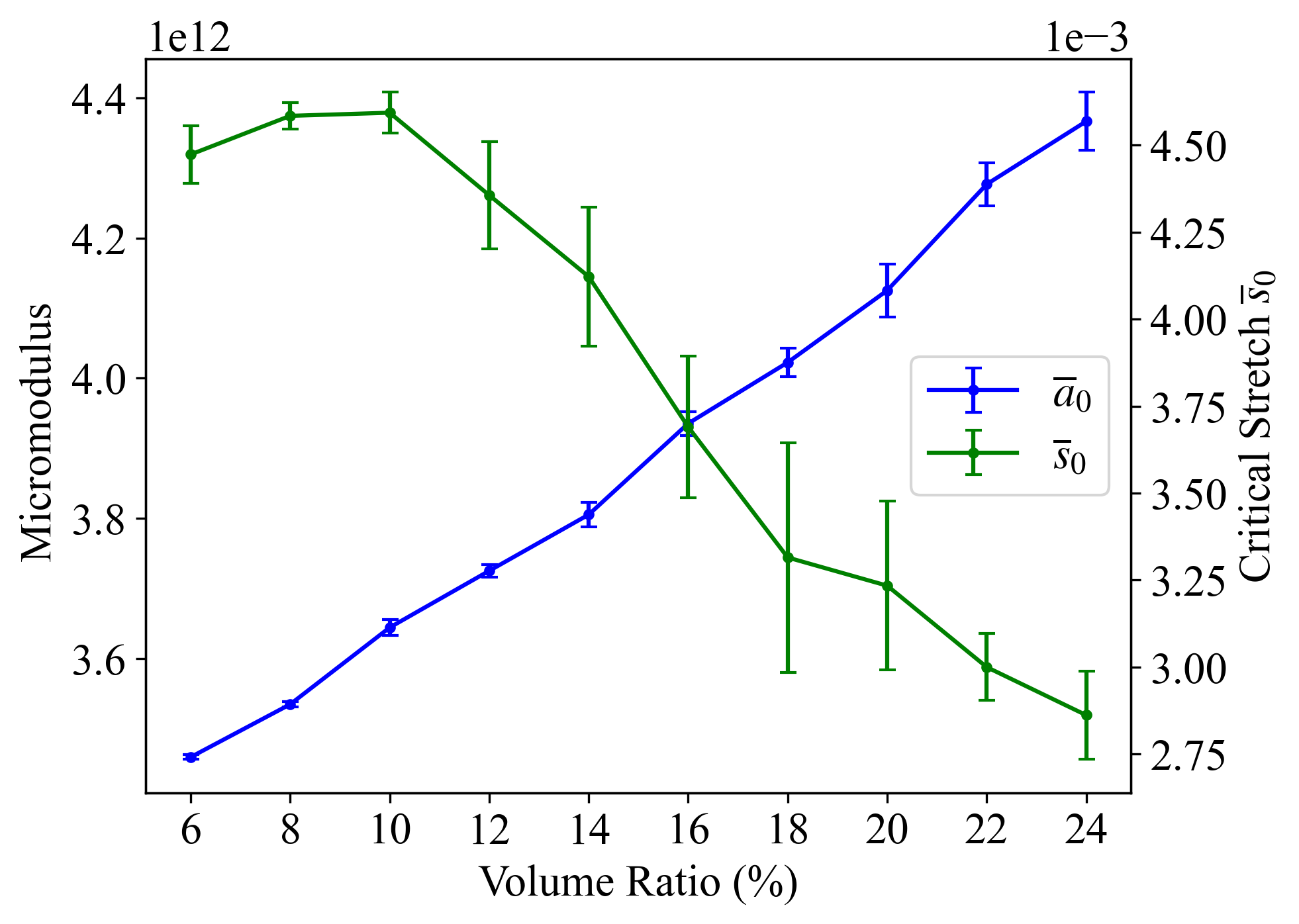}}
  \subfigure[]{
    \includegraphics[width=0.52\textwidth]{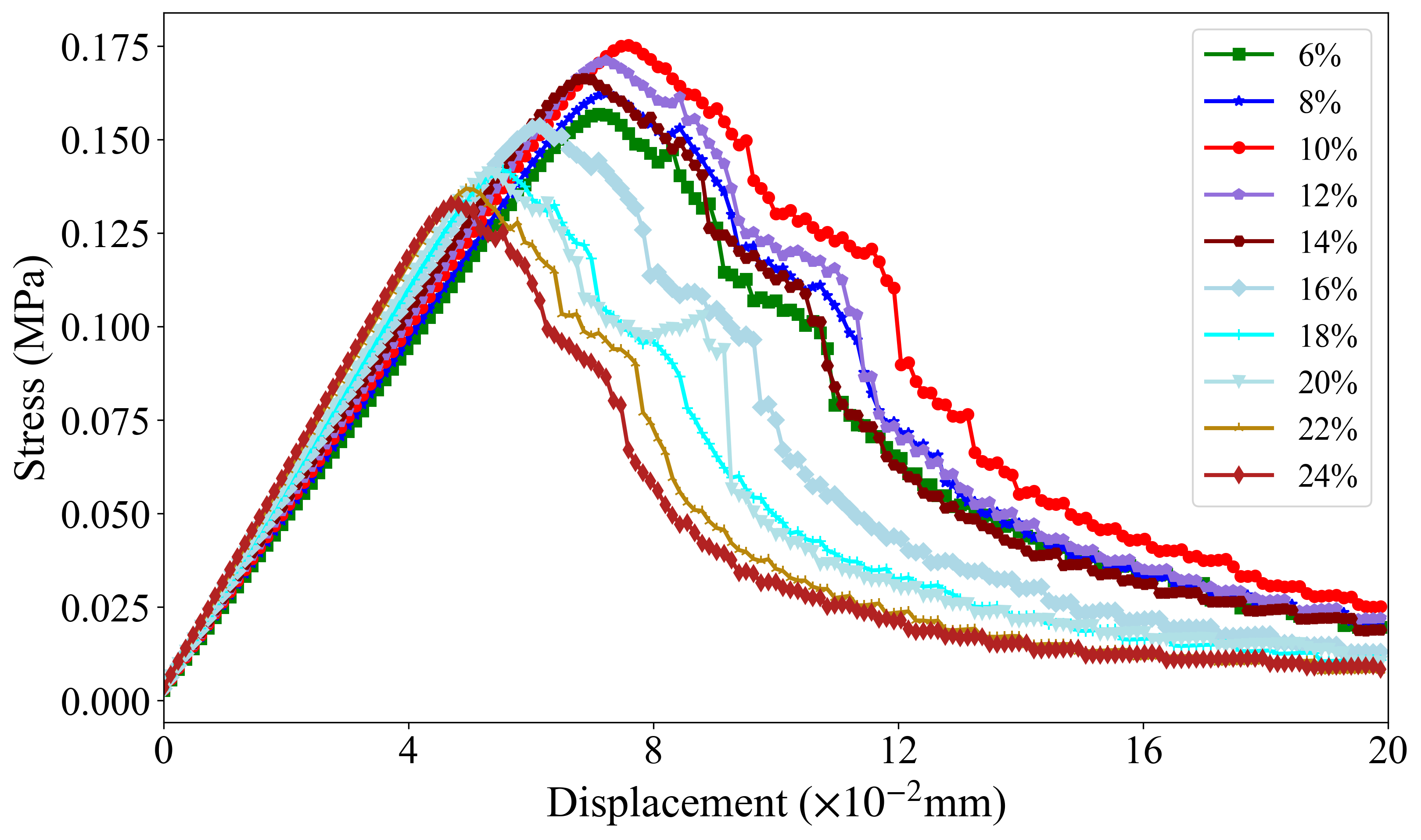}}
  \caption{(a) The coefficient $\bar{a}_0$ of equivalent micromodulus defined by Eq. \eqref{eq:emirco-exa} and critical stretch versus volume fraction of particles, (b) stress versus imposed displacement for homogenized structures corresponding to RVEs with different volume fractions of particles.}\label{fig:volume_curve}
\end{figure}

\begin{figure}[H]
  \centering
  \subfigure[10\%]{
    \includegraphics[width=0.43\textwidth]{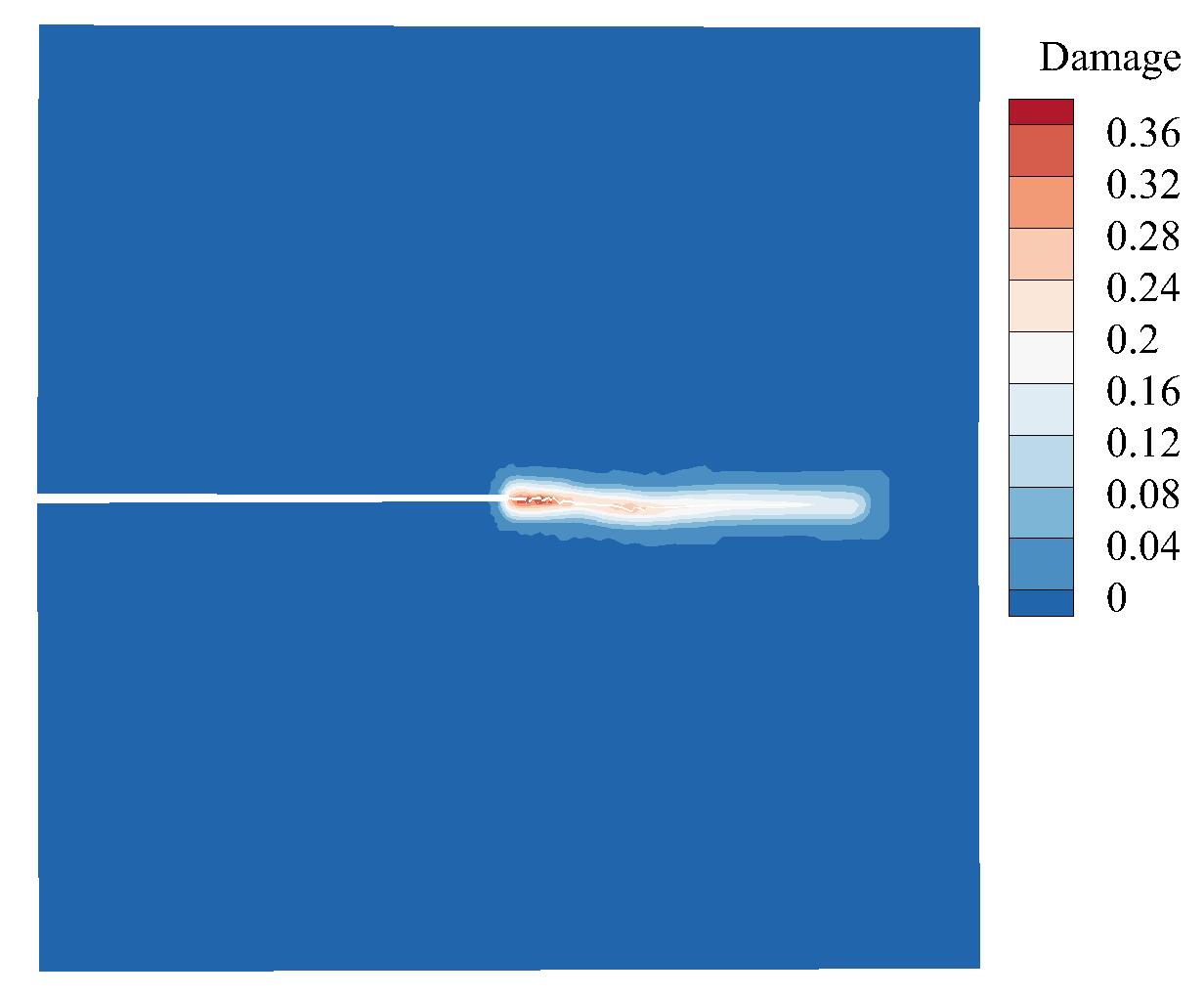}\label{fig:volume_10_frac}}
  \subfigure[24\%]{
    \includegraphics[width=0.43\textwidth]{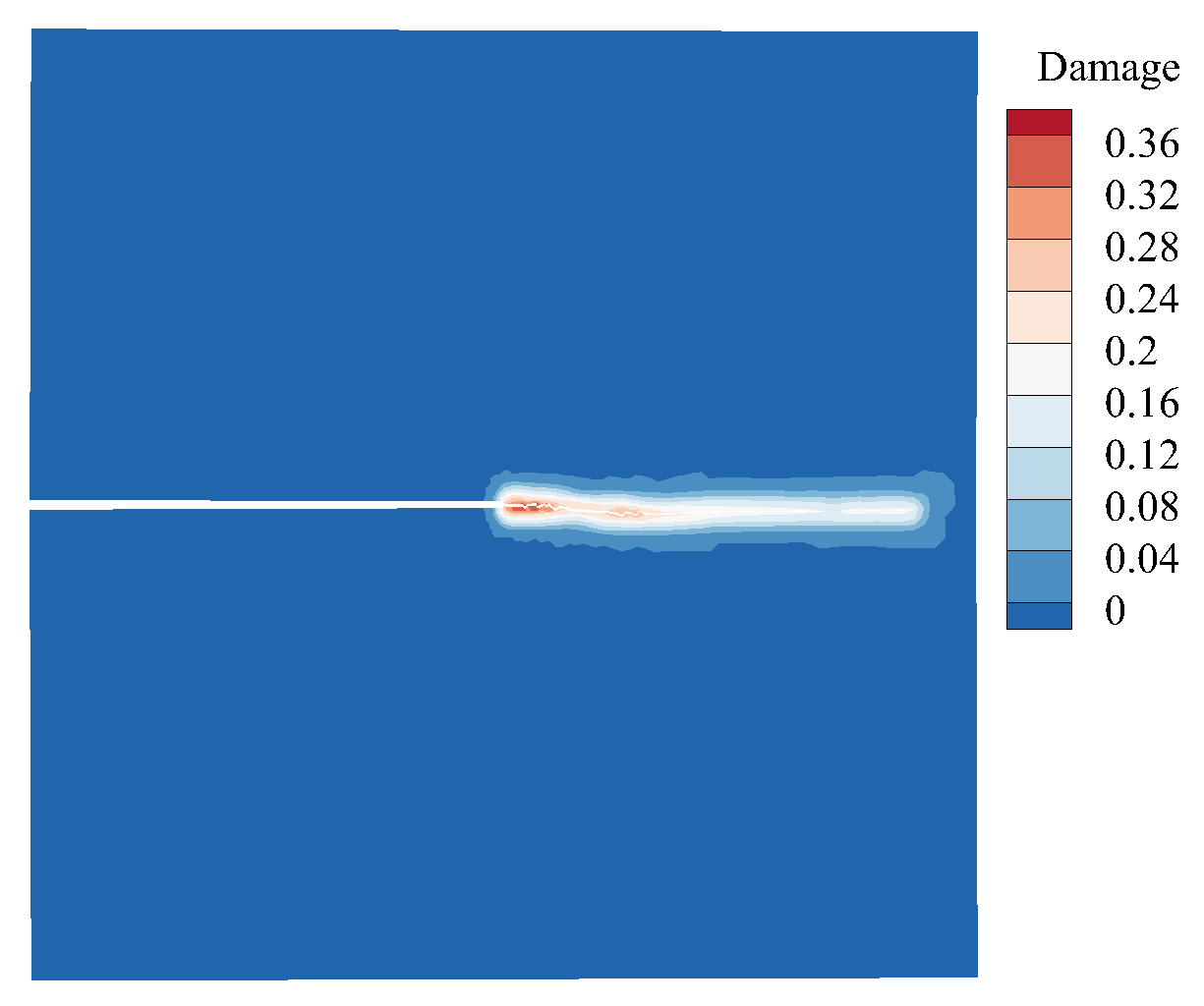}\label{fig:volume_24_frac}}
  \caption{Crack paths of homogenized structures corresponding to RVEs with different volume fractions of particles at Step 166.}\label{fig:volume_frac_macro}
\end{figure}

\autoref{fig:volume_curve}(a) displays the expected values of $\bar{a}_0$ in micromodulus defined by Eq. \eqref{eq:emirco-exa} and critical stretch of the RVEs with different volume fractions of particles.
Statistically, different samples should
have different results. But accompanied by increasing number of samples with the same statistical characteristic, the mathematical expectation of computational results should converge.
Therefore, some samples were taken in this study to avoid an unacceptable scatter of numerical results.
From \autoref{fig:volume_curve}(a), it can be seen that with the increase of particle volume fraction, the equivalent micromodulus increases continuously, while the statistical critical stretch first increases slightly and then decreases continuously. This shows that although the particles play a reinforcing role in the composites, it does not mean that the more particles, the greater the statistical critical stretch of the composites. There exists an optimal reinforcement volume fraction \cite{jamshaid2022natural,ali2014seismic,duxiaoqi}, that is, in this example, the statistical critical stretch reaches the maximum value when the particle volume fraction is $10\%$, and when the particle volume fraction is more than $10\%$, the more particles, the smaller the statistical critical stretch of the composites.
\autoref{fig:volume_curve}(b) displays the average stress at the load boundary versus imposed displacement curve for the homogenized structures corresponding to the RVEs with different volume fractions of particles.
From \autoref{fig:volume_curve}(b), it can be found that the average stress versus imposed displacement curves for different particle fractions have similar morphology, that is, with the increase of imposed displacement, the stress first increases approximately linearly, and then decreases gradually after reaching the peak value.
The final crack path is shown in \autoref{fig:volume_frac_macro}. It can be found from \autoref{fig:volume_frac_macro} that the crack growth paths of the homogenized structures corresponding to the RVEs with $10\%$ and $24\%$ particle fractions basically remain horizontal, and the crack for the $24\%$ case is longer than that for the $10\%$ case.

By comparing \autoref{fig:volume_curve}(a) and \autoref{fig:volume_curve}(b), it can be found that 10 groups of results can be roughly divided into left and right parts (see \autoref{fig:volume_curve}(a)) and upper and lower parts (see \autoref{fig:volume_curve}(b)), respectively, with the data point in \autoref{fig:volume_curve}(a) or the curve in \autoref{fig:volume_curve}(b) corresponding to $16\%$ particle fraction as the boundary.
From the left part of \autoref{fig:volume_curve}(a), it is easy to see that the curve of statistical critical stretch is on the top of the curve of equivalent micromodulus when the particle fraction is less than $16\%$, and the result is opposite when the particle fraction is greater than $16\%$ in the right part of \autoref{fig:volume_curve}(b). We know that under the same deformation condition, the smaller the critical stretch, the easier the bond will break. On the other hand, under the condition of bearing the same bond force, the higher the bond micromodulus, the harder the bond will elongate and the harder the bond will fracture.
Therefore, there exists a competition between the statistical critical stretch and the equivalent micromodulus of composite materials for structural failure.
This conclusion can also be drawn in \autoref{fig:volume_curve}(b). For the cases that the particle fraction is less than $16\%$ (the upper curves), when the statistical critical stretch reaches the maximum (corresponding to the case of particle fraction $10\%$), the curve reaches the maximum stress peak value and imposed displacement value.
In addition, the peak stress of the curves corresponding to $6\%$ and $8\%$ particle fraction is smaller than that of the curve corresponding to $12\%$ and $14\%$ particle fraction, respectively. Although the statistical critical stretch of the former is much larger than that of the latter, the equivalent micromodulus of the former is much smaller than that of the latter, resulting in a lower stress peak value. This demonstrates that when the equivalent micromodulus is small enough, even if the statistical critical stretch is large, it cannot play a good role in retarding structural damage.
On the other hand, when the particle fraction is greater than $16\%$ (the lower curves), the statistical critical stretch plays a decisive role. It can be seen from \autoref{fig:volume_curve}(b) that even if the equivalent micromodulus gradually increases, the peak stress of the corresponding curve gradually decreases with the gradual decrease of the statistical critical stretch. We can conclude from the above analysis that the equivalent micromodulus and statistical critical stretch compete with each other and jointly affect the fracture of the macroscale structure. Therefore, one can design the microstructure of composite materials to obtain composite structures with specific fracture performance.

\subsection{Applications of the PSM method to fracture in 3D random composite structures}

\begin{figure}[H]
  \centering
  \includegraphics[scale=0.7]{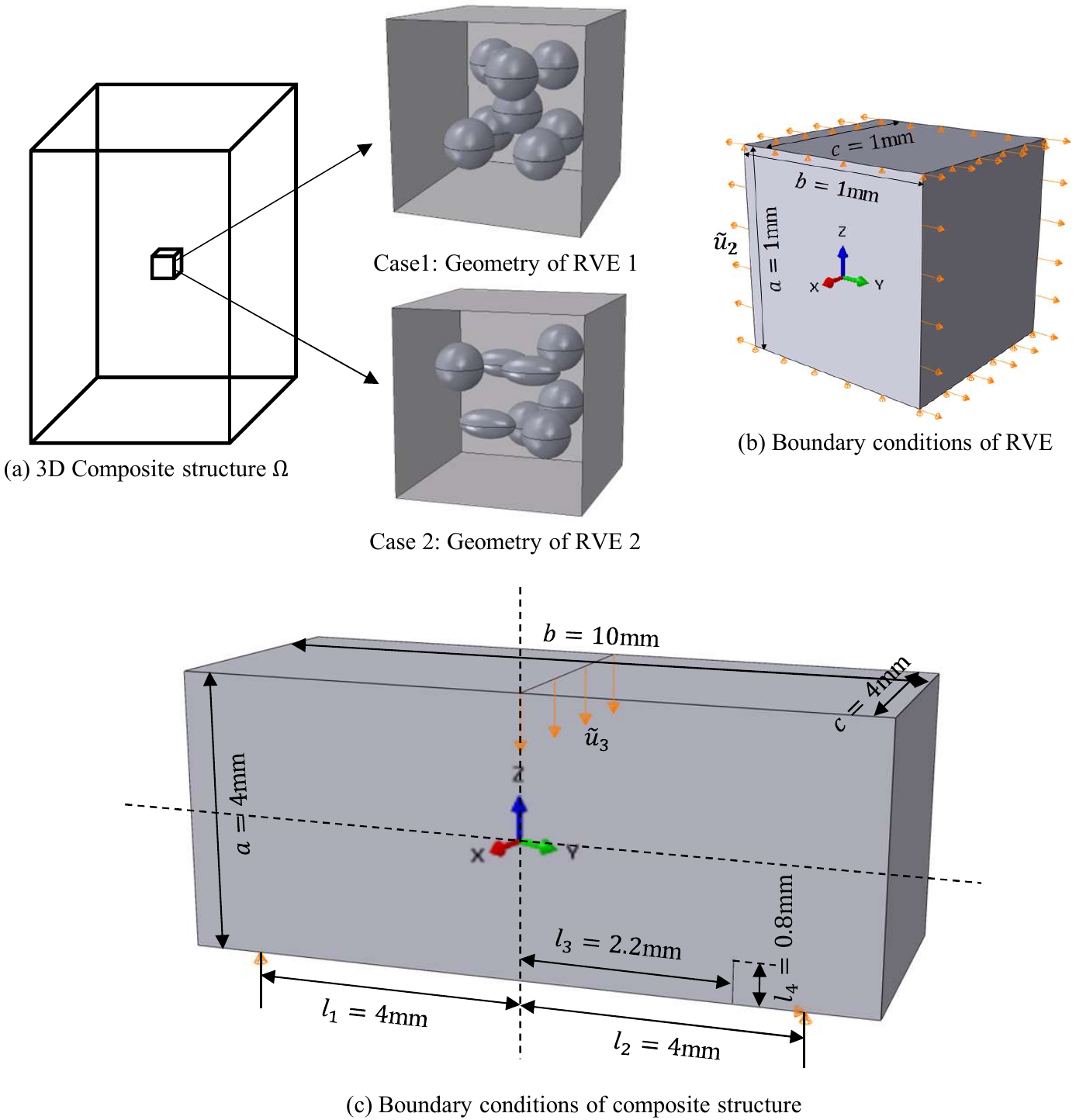}
  \caption{Geometry and boundary conditions of 3D composite structure and related RVEs.}\label{Fig:sec5.3-1}
\end{figure}

In this section, we apply the PSM method to simulate the fracture in 3D random composite structures.
\autoref{Fig:sec5.3-1}(c) shows the geometry and boundary condition with $\tilde{u}_3=0.1$mm of the composite structures, where a notch with length of 4 mm, width of 0.02 mm and height of 0.8mm is preset. Two kinds of RVEs reinforced by the uniform distribution of spherical and ellipsoidal particles, whose volume fraction is taken as 7.5\%, are considered, as shown in \autoref{Fig:sec5.3-1}(a).
The biaxial tensile boundary condition with $\tilde{u}_2=0.02$ mm along $y_2$-axis is specified on the left and right surfaces of the RVEs, and the top and bottom surfaces are fixed in the $y_3$ direction, as shown in \autoref{Fig:sec5.3-1}(b).
The BPD simulations of the RVE and homogenized structures were implemented through 50 displacement increments (steps).



\begin{figure}[H]
  \centering
  \subfigure[]{
    \includegraphics[width=0.33\textwidth]{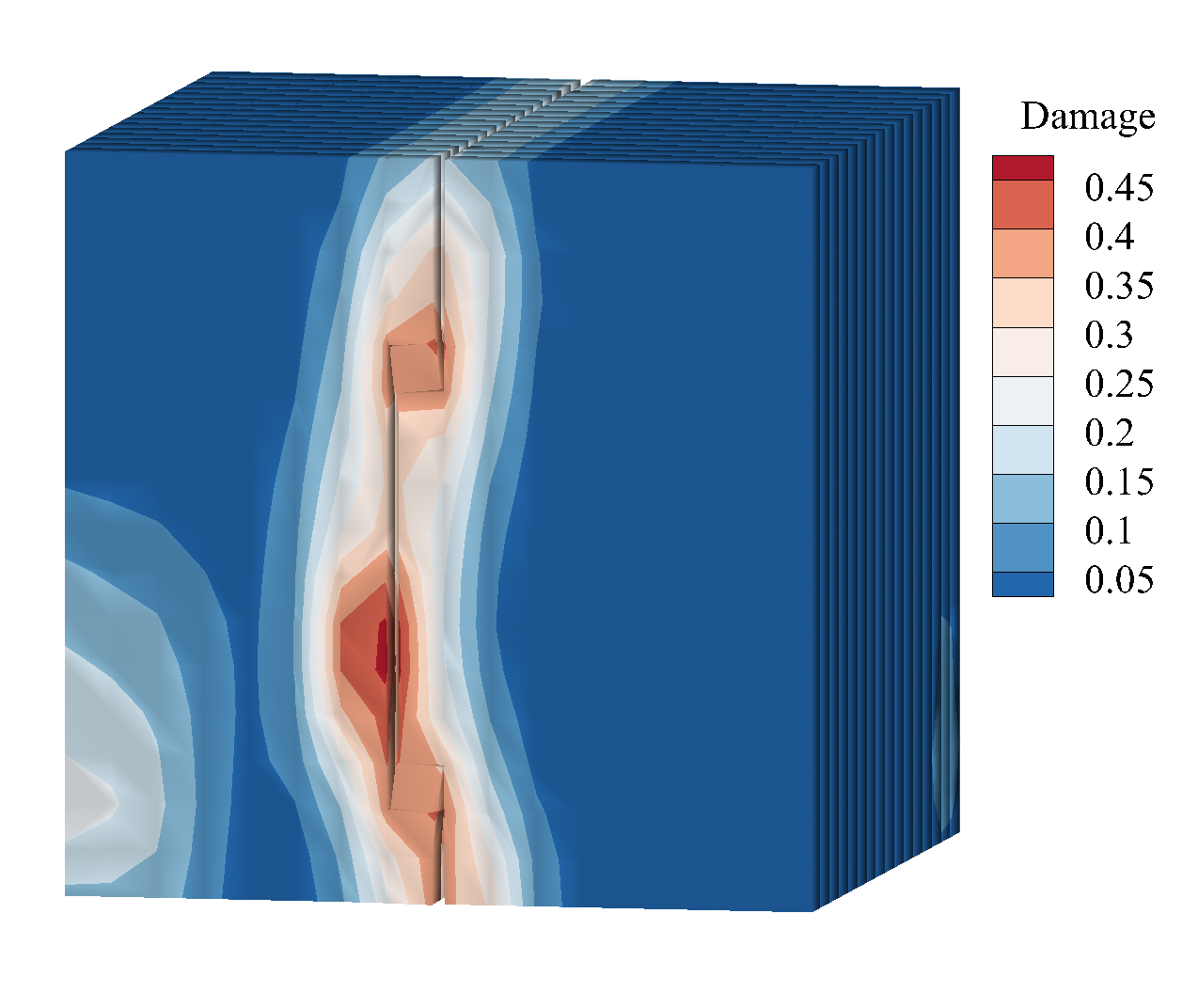}\label{fig:frac_3d_1}}
  \subfigure[]{
    \includegraphics[width=0.33\textwidth]{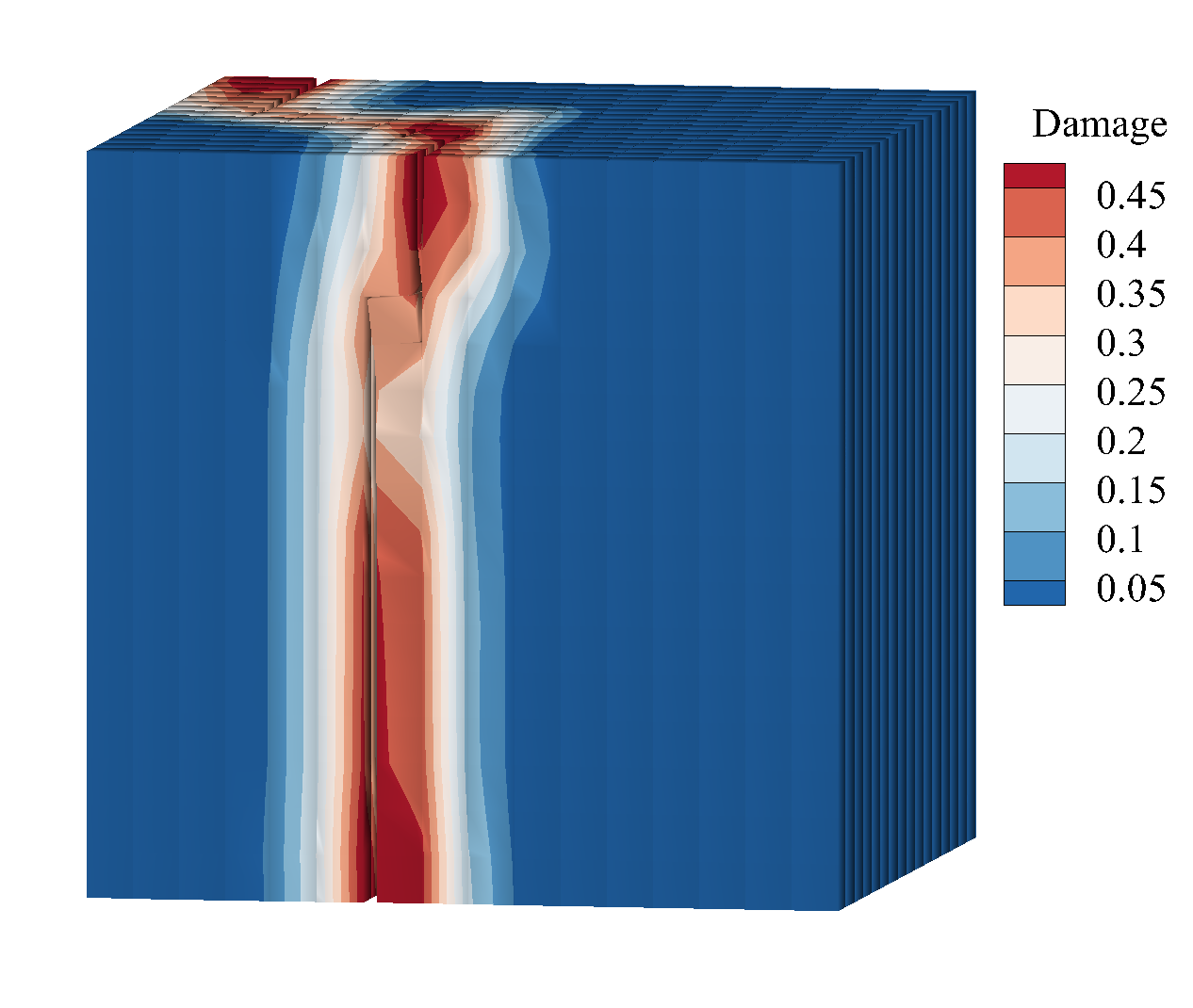}\label{fig:frac_3d_3}}
  \subfigure[]{
    \includegraphics[width=0.28\textwidth]{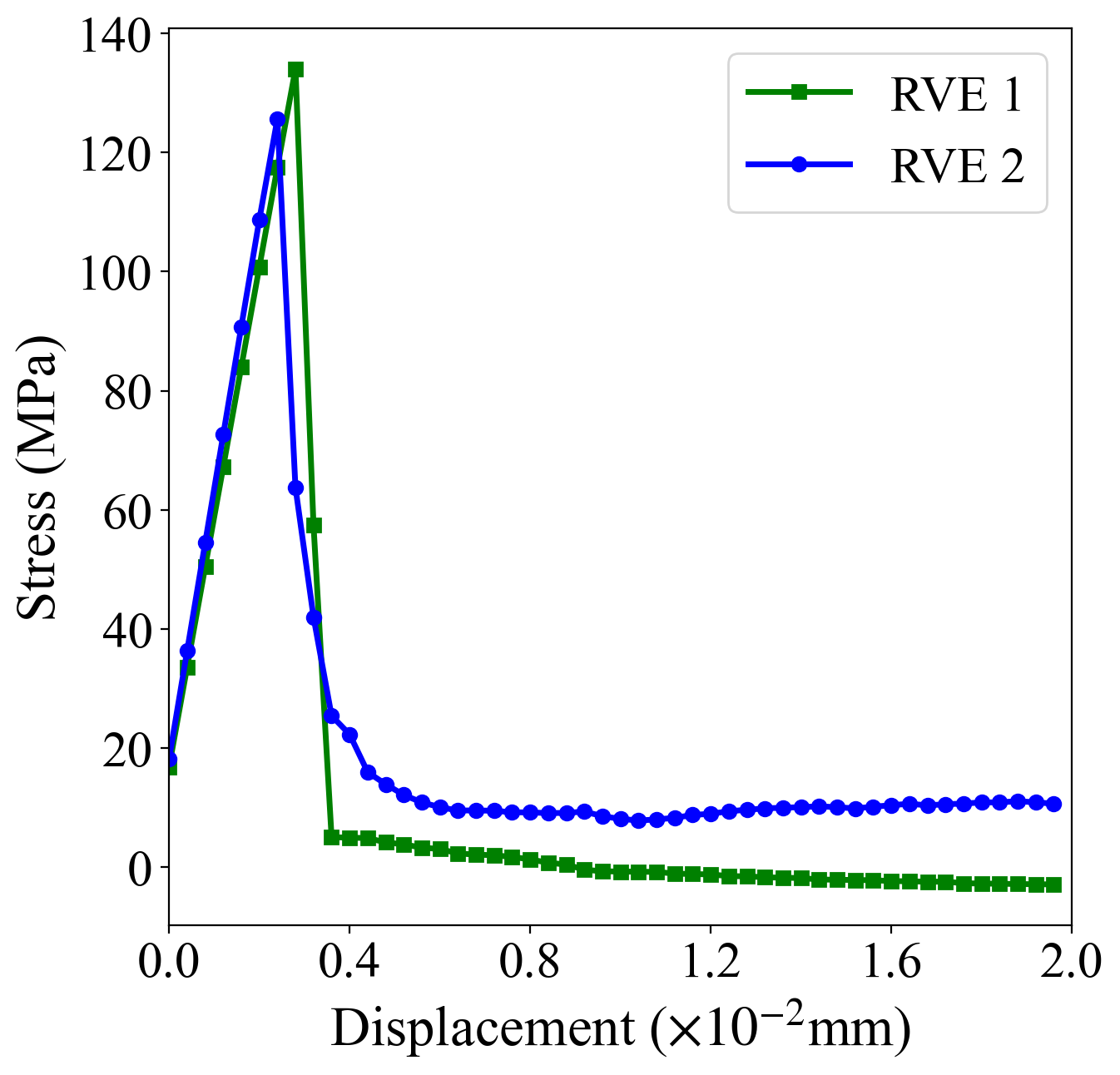}\label{fig:3d_1_2_curve_mic}}
  \caption{Crack path for (a) 3D RVE 1 and (b) 3D RVE 2 at imposed displacement load step 10 and 8; (c) stress versus imposed displacement for two kinds of RVEs.}\label{fig:rve_frac_3d}
\end{figure}

\begin{figure}[H]
  \centering
  \subfigure[]{
    \includegraphics[width=0.42\textwidth]{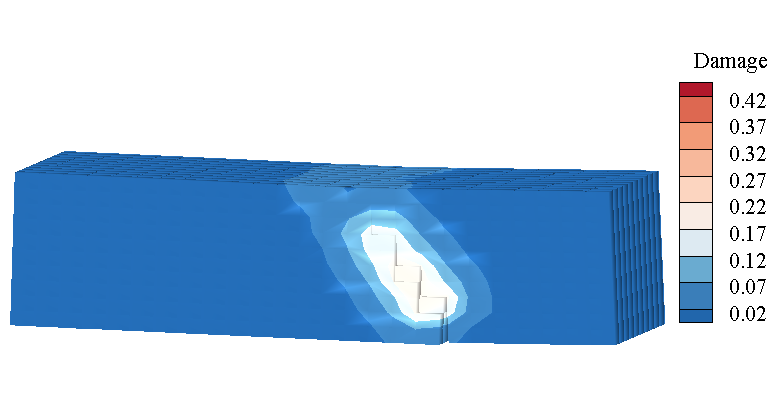}}
  \subfigure[]{
    \includegraphics[width=0.42\textwidth]{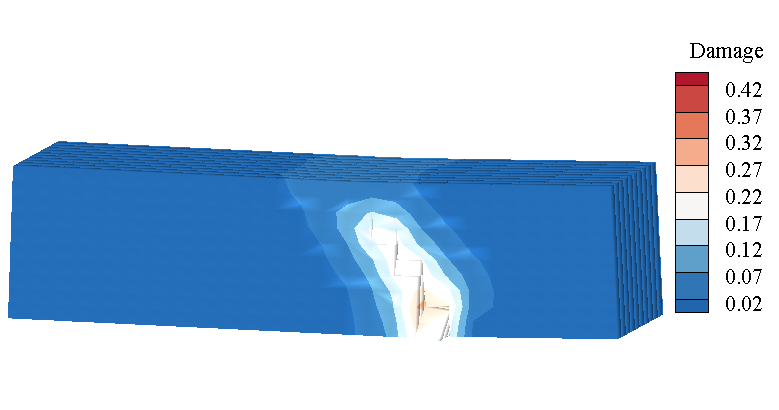}}
  \caption{Crack path in the homogenized structures of composites composed of (a) RVE1 and (b) RVE 2 with three point bending at step 50.}\label{fig:frac_3d_point}
\end{figure}

\begin{figure}[H]
  \centering
  
  \subfigure[]{
    \includegraphics[width=0.33\textwidth]{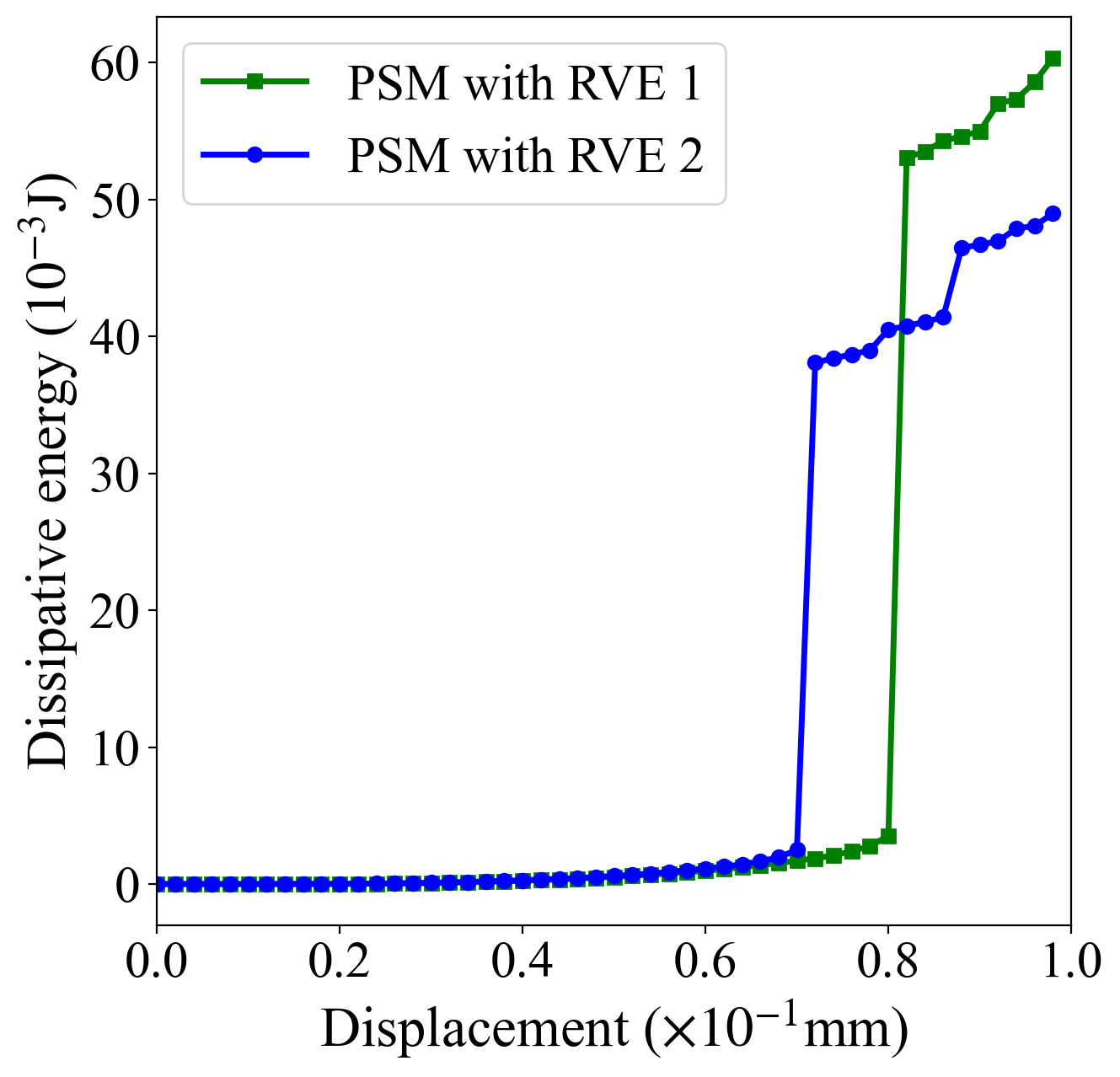}}
  \subfigure[]{
    \includegraphics[width=0.33\textwidth]{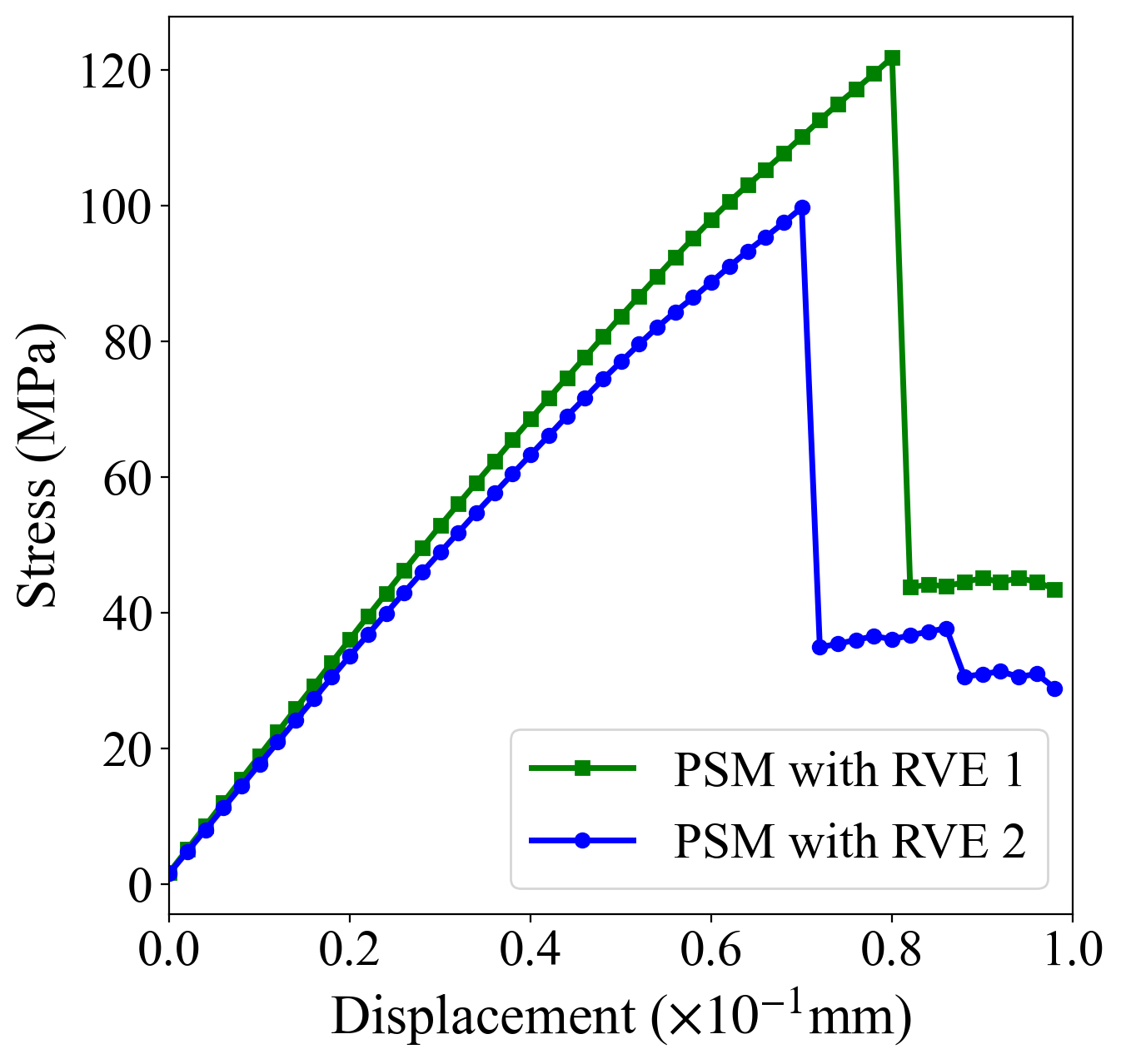}}
  \caption{(a) Dissipative energy and (b) stress versus imposed displacement for homogenized structures of composites composed of RVE 1 and RVE 2.}\label{fig:3d_curve}
\end{figure}

Since the 3D composite structure is composed of 160 RVEs, the number of elements and nodes of composite structure are much more than that of a single RVE and the homogenized structure, as shown in \autoref{tab:time_compare_3d}. Thus, the single-scale direct BPD fracture simulation of the 3D composite structures is almost impossible.
According the proposed PSM framework, we first take the fracture simulation of RVEs by the microscopic BPD model to determine the equivalent critical stretch. The crack path and stress curve versus imposed load steps for two kinds of RVEs are shown in \autoref{fig:rve_frac_3d}. Then, we apply the statistical homogenization method to obtain the equivalent micromodulus.
Finally, the composite structures can be homogenized to the macroscopic homogeneous structures, which can be simulated by the macroscopic PBD model with expectation of micromodulus and critical stretch.
\autoref{fig:frac_3d_point} shows the crack path for the 3D homogenized structures under three-point bending condition, and corresponding dissipative energy \cite{azdoud2014morphing} and stress curves are displayed in \autoref{fig:3d_curve}.


\begin{table}[H]
  \setlength{\abovecaptionskip}{0cm}
  \setlength{\belowcaptionskip}{0.3cm}
  \centering
  \caption{Comparison of computational time for 3D composite structures, RVEs and homogenized structures.}\label{tab:time_compare_3d}
  \scalebox{0.8}{
    \begin{tabular}{cccccccc}
      \toprule  
      \multirow{2}*{} & \multicolumn{3}{c}{Composite with RVE 1} & \quad   & \multicolumn{3}{c}{Composite with RVE 2}                                                       \\ \cline{2-4} \cline{6-8}
                      & Composite                                & RVE 1   & Homogenized structure                    & \quad & Composite & RVE 2   & Homogenized structure \\ \hline
      Elements        & 655,360                                  & 4,096   & 2,744                                    & \quad & 655,360   & 4,096   & 2,744                 \\
      CE Nodes        & 786,080                                  & 4,913   & 3,375                                    & \quad & 786,080   & 4,913   & 3,375                 \\
      DE Nodes        & 5,242,880                                & 32,768  & 21,952                                   & \quad & 5,242,880 & 32,768  & 21,952                \\
      Time(s)         & -                                        & 560,010 & 367,668                                  & \quad & -         & 571,430 & 351,020               \\
      \bottomrule  
    \end{tabular}
  }
\end{table}

\section{Conclusions}
The PD models are promising for the simulation of fracture, as they allow discontinuities in the displacement field. However, they are still computationally expensive, especially for fracture simulation of large-scale composite structures with randomly distributed of particles.
Then, a novel PSM method is presented to predict
the fracture of composite structures with the heterogeneities at microscale being taken into account, which makes it possible to efficiently simulate the fracture of composite structures.
A key point of the approach is that the scale separation-based statistical multiscale strategy is used to define a macroscale BPD model with an equivalent micromodulus and a statistical critical stretch that are evaluated based on microscale RVEs.
Consequently, the proposed method
can take into account the microstructural characteristics of composites when performing macroscopic fracture analysis, and at the same time save the computational cost.
The validity and efficiency of the PSM method and its CE/DE numerical algorithm have been verified by comparison with the single-scale direct PD fracture analyses.
Moreover, the equivalent micromodulus and statistical critical stretch compete with each other and jointly impact the fracture performance of macroscale structure.

\section{Acknowledgements}
This research was supported by the National Natural Science Foundation of China (12272082, 51739007),  Strategic Priority Research Program of Chinese Academy of Sciences (XDC06030102), Special Scientific Research Program of Shaanxi Provincial Department of Education (22JK0586) and Natural Science Foundation of Chongqing (CSTB2022NSCQ-MSX0296).

\bibliography{ref}

\begin{thebibliography}{10}
\expandafter\ifx\csname url\endcsname\relax
  \def\url#1{\texttt{#1}}\fi
\expandafter\ifx\csname urlprefix\endcsname\relax\def\urlprefix{URL }\fi
\expandafter\ifx\csname href\endcsname\relax
  \def\href#1#2{#2} \def\path#1{#1}\fi

\bibitem{shackelford2016introduction}
J.~F. Shackelford, Introduction to materials science for engineers, Pearson
  Upper Saddle River, 2016.

\bibitem{tang2013study}
X.~Tang, Y.~Zhou, C.~Zhang, et~al., Study on the heterogeneity of concrete and
  its failure behavior using the equivalent probabilistic model, in: Seismic
  Safety Evaluation of Concrete Dams, Elsevier, 2013, pp. 541--570.

\bibitem{hillerborg1976analysis}
A.~Hillerborg, M.~Mod{\'e}er, P.~Petersson, Analysis of crack formation and
  crack growth in concrete by means of fracture mechanics and finite elements,
  Cement. Concrete. Res. 6~(6) (1976) 773--781.

\bibitem{kassam1995finite}
Z.~H. Kassam, R.~J. Zhang, Z.~Wang, Finite element simulation to investigate
  interaction between crack and particulate reinforcements in metal-matrix
  composites, Mat. Sci. Eng. A-Struct. 203~(1-2) (1995) 286--299.

\bibitem{ayyar2006microstructure}
A.~Ayyar, N.~Chawla, Microstructure-based modeling of crack growth in particle
  reinforced composites, Compos. Sci. Technol. 66~(13) (2006) 1980--1994.

\bibitem{sukumar2001modeling}
N.~Sukumar, D.~L. Chopp, N.~Mo{\"e}s, et~al., Modeling holes and inclusions by
  level sets in the extended finite-element method, Comput. Methods Appl. Mech.
  Engrg. 190~(46-47) (2001) 6183--6200.

\bibitem{huynh2009extended}
D.~Huynh, T.~Belytschko, The extended finite element method for fracture in
  composite materials, Int. J. Numer. Meth. Eng. 77~(2) (2009) 214--239.

\bibitem{rostami2019xfem}
R.~B. Rostami, S.~Schmauder, {XFEM} simulation of fatigue crack growth in
  aluminum zirconia reinforced composites, Int. J. Multiscale. Com. 17~(5)
  (2019) 469--481.

\bibitem{spangenberger2021extended}
A.~G. Spangenberger, D.~A. Lados, Extended finite element modeling of fatigue
  crack growth microstructural mechanisms in alloys with secondary/reinforcing
  phases: model development and validation, Comput. Mech. 67~(1) (2021)
  87--105.

\bibitem{sun2006modeling}
C.~Sun, Z.~Jin, Modeling of composite fracture using cohesive zone and bridging
  models, Compos. Sci. Technol. 66~(10) (2006) 1297--1302.

\bibitem{de2013mode}
A.~De~Morais, Mode {I} cohesive zone model for delamination in composite beams,
  Eng. Fract. Mech. 109 (2013) 236--245.

\bibitem{sun2020prediction}
L.~Sun, Y.~Tie, Y.~Hou, et~al., Prediction of failure behavior of adhesively
  bonded {CFRP} scarf joints using a cohesive zone model, Eng. Fract. Mech. 228
  (2020) 106897.

\bibitem{quintanas2020phase}
A.~Quintanas-Corominas, A.~Turon, J.~Reinoso, et~al., A phase field approach
  enhanced with a cohesive zone model for modeling delamination induced by
  matrix cracking, Comput. Methods Appl. Mech. Engrg. 358 (2020) 112618.

\bibitem{ghosh2013computational}
A.~Ghosh, P.~Chaudhuri, Computational modeling of fracture in concrete using a
  meshfree meso-macro-multiscale method, Comp. Mater. Sci. 69 (2013) 204--215.

\bibitem{bui2018analysis}
T.~Q. Bui, N.~T. Nguyen, M.~N. Nguyen, et~al., Analysis of transient dynamic
  fracture parameters of cracked functionally graded composites by improved
  meshfree methods, Theor. Appl. Fract. Mec. 96 (2018) 642--657.

\bibitem{nguyen2015phase}
T.~T. Nguyen, J.~Yvonnet, Q.~Zhu, et~al., A phase field method to simulate
  crack nucleation and propagation in strongly heterogeneous materials from
  direct imaging of their microstructure, Eng. Fract. Mech. 139 (2015) 18--39.

\bibitem{kuhn2016discussion}
C.~Kuhn, R.~M{\"u}ller, A discussion of fracture mechanisms in heterogeneous
  materials by means of configurational forces in a phase field fracture model,
  Comput. Methods Appl. Mech. Engrg. 312 (2016) 95--116.

\bibitem{zhang2020modelling}
P.~Zhang, Y.~Feng, T.~Q. Bui, et~al., Modelling distinct failure mechanisms in
  composite materials by a combined phase field method, Compos. Struct. 232
  (2020) 111551.

\bibitem{xia2021mesoscopic}
Y.~Xia, W.~Wu, Y.~Yang, et~al., Mesoscopic study of concrete with random
  aggregate model using phase field method, Constr. Build. Mater. 310 (2021)
  125199.

\bibitem{yang2020stochastic}
Z.~Yang, X.~Guan, J.~Cui, et~al., Stochastic multiscale heat transfer analysis
  of heterogeneous materials with multiple random configurations, Commun.
  Comput. Phys. 27~(2) (2020) 431--459.

\bibitem{shu2020multiscale}
W.~Shu, I.~Stanciulescu, Multiscale homogenization method for the prediction of
  elastic properties of fiber-reinforced composites, Int. J. Solids. Struct.
  203 (2020) 249--263.

\bibitem{yang2020novel}
Z.~Yang, Y.~Liu, Y.~Sun, et~al., A novel second-order reduced homogenization
  approach for nonlinear thermo-mechanical problems of axisymmetric structures
  with periodic micro-configurations, Comput. Methods Appl. Mech. Engrg. 368
  (2020) 113126.

\bibitem{yang2022efficient}
Z.~Yang, J.~Huang, X.~Feng, et~al., An efficient multi-modes monte carlo
  homogenization method for random materials, SIAM. J. Sci. Comput. 44~(3)
  (2022) A1752--A1774.

\bibitem{yang2022second}
Z.~Yang, Y.~Liu, Y.~Sun, et~al., A second-order reduced multiscale method for
  nonlinear shell structures with orthogonal periodic configurations, Int. J.
  Numer. Meth. Eng. 123~(1) (2022) 128--157.

\bibitem{yang2017high}
Z.~Yang, Y.~Zhang, H.~Dong, et~al., High-order three-scale method for
  mechanical behavior analysis of composite structures with multiple periodic
  configurations, Compos. Sci. Technol. 152 (2017) 198--210.

\bibitem{ma2018multi}
Q.~Ma, Z.~Li, J.~Cui, Multi-scale asymptotic analysis and computation of the
  elliptic eigenvalue problems in curvilinear coordinates, Comput. Methods
  Appl. Mech. Engrg. 340 (2018) 340--365.

\bibitem{dong2022stochastic}
H.~Dong, Z.~Yang, X.~Guan, et~al., Stochastic higher-order three-scale strength
  prediction model for composite structures with micromechanical analysis, J.
  Comput. Phys. (2022) 111352.

\bibitem{canal2012intraply}
L.~P. Canal, C.~Gonz{\'a}lez, J.~Segurado, et~al., Intraply fracture of
  fiber-reinforced composites: Microscopic mechanisms and modeling, Compos.
  Sci. Technol. 72~(11) (2012) 1223--1232.

\bibitem{silling2000reformulation}
S.~A. Silling, Reformulation of elasticity theory for discontinuities and
  long-range forces, J. Mech. Phys. Solids. 48~(1) (2000) 175--209.

\bibitem{silling2005meshfree}
S.~A. Silling, E.~Askari, A meshfree method based on the peridynamic model of
  solid mechanics, Comput. Struct. 83~(17-18) (2005) 1526--1535.

\bibitem{hu2012peridynamic}
W.~Hu, Y.~D. Ha, F.~Bobaru, Peridynamic model for dynamic fracture in
  unidirectional fiber-reinforced composites, Comput. Methods Appl. Mech.
  Engrg. 217 (2012) 247--261.

\bibitem{zhou2017analyzing}
W.~Zhou, D.~Liu, N.~Liu, Analyzing dynamic fracture process in fiber-reinforced
  composite materials with a peridynamic model, Eng. Fract. Mech. 178 (2017)
  60--76.

\bibitem{sau2019peridynamic}
N.~Sau, J.~Medina-Mendoza, A.~C. Borbon-Almada, Peridynamic modelling of
  reinforced concrete structures, Eng. Fail. Anal. 103 (2019) 266--274.

\bibitem{li2018meso}
W.~Li, L.~Guo, Meso-fracture simulation of cracking process in concrete
  incorporating three-phase characteristics by peridynamic method, Constr.
  Build. Mater. 161 (2018) 665--675.

\bibitem{dong2021improved}
Y.~Dong, C.~Su, P.~Qiao, An improved mesoscale damage model for quasi-brittle
  fracture analysis of concrete with ordinary state-based peridynamics, Theor.
  Appl. Fract. Mec. 112 (2021) 102829.

\bibitem{peng2021application}
R.~Peng, W.~Qiu, M.~Jiang, Application of a micro-model for concrete to the
  simulation of crack propagation, Theor. Appl. Fract. Mec. 116 (2021) 103081.

\bibitem{mehrmashhadi2019stochastically}
J.~Mehrmashhadi, Z.~Chen, J.~Zhao, et~al., A stochastically homogenized
  peridynamic model for intraply fracture in fiber-reinforced composites,
  Compos. Sci. Technol. 182 (2019) 107770.

\bibitem{wu2021stochastically}
P.~Wu, F.~Yang, Z.~Chen, et~al., Stochastically homogenized peridynamic model
  for dynamic fracture analysis of concrete, Eng. Fract. Mech. 253 (2021)
  107863.

\bibitem{han2010statistical}
F.~Han, J.~Cui, Y.~Yu, The statistical second-order two-scale method for
  mechanical properties of statistically inhomogeneous materials, Int. J.
  Numer. Meth. Eng. 84~(8) (2010) 972--988.

\bibitem{jikov2012homogenization}
V.~V. Jikov, S.~M. Kozlov, O.~A. Oleinik, Homogenization of differential
  operators and integral functionals, Springer Science \& Business Media, 2012.

\bibitem{YANG2022CICP}
Z.~Yang, S.~Zheng, F.~Han, et~al., An improved peridynamic model with
  energy-based micromodulus correction method for fracture in particle
  reinforced composites, Commun. Comput. Phys. 32 (2022) 424--449.

\bibitem{madenci2020peridynamic}
E.~Madenci, A.~Yaghoobi, A.~Barut, et~al., Peridynamic unit cell for effective
  properties of complex microstructures with and without defects, Theor. Appl.
  Fract. Mec. 110 (2020) 102835.

\bibitem{azdoud2013morphing}
Y.~Azdoud, F.~Han, G.~Lubineau, A morphing framework to couple non-local and
  local anisotropic continua, Int. J. Solids. Struct. 50~(9) (2013) 1332--1341.

\bibitem{li2023peridynamics}
Z.~Li, F.~Han, The peridynamics-based finite element method (perifem) with
  adaptive continuous/discrete element implementation for fracture simulation,
  Eng. Anal. Bound. Elem. 146 (2023) 56--65.

\bibitem{nooru1993experimental}
M.~Nooru-Mohamed, E.~Schlangen, J.~G. van Mier, Experimental and numerical
  study on the behavior of concrete subjected to biaxial tension and shear,
  Adv. Cement. Based. Mater. 1~(1) (1993) 22--37.

\bibitem{winkler2001experimental}
B.~Winkler, G.~Hofstetter, G.~Niederwanger, Experimental verification of a
  constitutive model for concrete cracking, P. I. Mech. Eng. L-J. Mat. 215~(2)
  (2001) 75--86.

\bibitem{narayan2019gradient}
S.~Narayan, L.~Anand, A gradient-damage theory for fracture of quasi-brittle
  materials, J. Mech. Phys. Solids. 129 (2019) 119--146.

\bibitem{jamshaid2022natural}
H.~Jamshaid, R.~K. Mishra, A.~Raza, et~al., Natural cellulosic fiber reinforced
  concrete: Influence of fiber type and loading percentage on mechanical and
  water absorption performance, Materials 15~(3) (2022) 874.

\bibitem{ali2014seismic}
M.~Ali, Seismic performance of coconut-fibre-reinforced-concrete columns with
  different reinforcement configurations of coconut-fibre ropes, Constr. Build.
  Mater. 70 (2014) 226--230.

\bibitem{duxiaoqi}
X.~Du, Z.~Si, L.~Huang, et~al., Study on the influence of coarse aggregate
  particle shape and volume fraction on the uniaxial compressive properties of
  concrete, Chinese Journal of Applied Mechanics 37 (2020) 1828--1834.

\bibitem{azdoud2014morphing}
Y.~Azdoud, F.~Han, G.~Lubineau, The morphing method as a flexible tool for
  adaptive local/non-local simulation of static fracture, Comput. Mech. 54~(3)
  (2014) 711--722.

\end{thebibliography}

\end{document}